\documentclass[3p]{elsarticle}

\usepackage{lineno,hyperref}
\usepackage{amsfonts}
\usepackage{amsmath}
\usepackage{amssymb}
\usepackage{amsthm}
\usepackage{mathrsfs}
\usepackage{subcaption}
\usepackage{float}
\usepackage{xfrac}
\usepackage[ruled,vlined,linesnumbered]{algorithm2e}
\usepackage{multicol}
\usepackage{fonts,defs}
\usepackage{caption}
\usepackage{subcaption}
\usepackage{cleveref}
\usepackage{bm}

\modulolinenumbers[5]
\hypersetup{
  bookmarksnumbered = true,
  bookmarksopen=false,
  pdfborder=0 0 0,         
  pdffitwindow=true,      
  pdfnewwindow=true, 
  colorlinks=true,           
  linkcolor=blue,            
  citecolor=magenta,    
  filecolor=magenta,     
  urlcolor=cyan              
}

\journal{Journal of Computational Physics}

\newtheorem{defn}{Definition}

\newtheorem{prop}{Proposition}

\DeclareMathOperator*{\spann}{span}

\counterwithin{figure}{subsection}
\numberwithin{equation}{section}

\newcommand{\asgard}{{\rm ASGarD}}
\renewcommand{\dx}{{\,\mathrm{d}}}

\newcommand{\lavg}{\{\!\!\{}
\newcommand{\ravg}{\}\!\!\}}
\newcommand{\ljmp}{[\![}
\newcommand{\rjmp}{]\!]}
\newcommand{\ledge}{\big<\!\!\big<}
\newcommand{\redge}{\big>\!\!\big>}
\newcommand{\mT}{\mathcal{T}}

\newcommand{\W}{\Omega}
\renewcommand{\grad}{\nabla}

\newcommand{\mA}{\mathcal{A}}

\begin{document}

\begin{frontmatter}

\title{Sparse-grid Discontinuous Galerkin Methods for the Vlasov--Poisson--Lenard--Bernstein Model \tnoteref{copyright}\tnoteref{support}}
\tnotetext[copyright]{
This manuscript has been authored by UT-Battelle, LLC under Contract No. DE-AC05-00OR22725 with the U.S. Department of Energy. The United States Government retains and the publisher, by accepting the article for publication, acknowledges that the United States Government retains a non-exclusive, paid-up, irrevocable, world-wide license to publish or reproduce the published form of this manuscript, or allow others to do so, for United States Government purposes. The Department of Energy will provide public access to these results of federally sponsored research in accordance with the DOE Public Access Plan(\url{http://energy.gov/downloads/doe-public-access-plan}).}

\tnotetext[support]{This material is based upon work partially supported by the U.S. Department of Energy, Office of Science, Office of Advanced Scientific Computing Research, as part of their Applied Mathematics Research Program; the U.S. Department of Energy, Office of Science, Office of Fusion Energy Science as part of their Fusion Research Energy Program; and the Laboratory Directed Research and Development Program of Oak Ridge National Laboratory (ORNL), managed by UT-Battelle, LLC for the U.S. Department of Energy under Contract No. De-AC05-00OR22725.  This research used resources of the Oak Ridge Leadership Computing Facility at the Oak Ridge National Laboratory.
This work was performed under the auspices of the U.S. Department of Energy by Lawrence Livermore National Laboratory under Contract DE-AC52-07NA27344.}


\author[ornl-csmd]{Stefan Schnake\texorpdfstring{\corref{cor}}{*}}
\ead{schnakesr@ornl.gov}
\author[llnl]{Coleman Kendrick}
\ead{kendrick6@llnl.gov}
\author[ornl-csmd,utk-phys]{Eirik Endeve}
\ead{endevee@ornl.gov}
\author[ornl-csmd]{Miroslav Stoyanov}
\ead{stoyanovmk@ornl.gov}
\author[ornl-csmd]{Steven Hahn}
\ead{hahnse@ornl.gov}
\author[ornl-csmd,utk-math]{Cory D. Hauck}
\ead{hauckc@ornl.gov}
\author[csiro]{David L. Green}
\ead{d.green@csiro.au}
\author[ornl-fed]{Phil Snyder}
\ead{snyderpb@ornl.gov}
\author[type-one]{John Canik}
\ead{john.canik@typeoneenergy.com}

\cortext[cor]{Corresponding author.  Address: PO Box 2008, MS6013, Oak Ridge, TN 37831-6013}

\address[ornl-csmd]{Computer Science and Mathematics Division, Oak Ridge National Laboratory, Oak Ridge, TN 37831 USA }
\address[ornl-fed]{Fusion Energy Division, Oak Ridge National Laboratory, Oak Ridge, TN 37831 USA }
\address[llnl]{Center for Applied Scientific Computing, Lawrence Livermore National Laboratory, Livermore, CA 94550 USA}
\address[utk-phys]{Department of Physics and Astronomy, University of Tennessee Knoxville, TN 37996 USA}
\address[utk-math]{Department of Mathematics, University of Tennessee Knoxville, TN 37996 USA}
\address[type-one]{Type One Energy Group, Oak Ridge, TN 37830 USA}
\address[csiro]{Commonwealth Scientific and Industrial Research Organisation, Mayfield West, NSW, 2304, Australia}

\begin{abstract}
Sparse-grid methods have recently gained interest in reducing the computational cost of solving high-dimensional kinetic equations.
In this paper, we construct adaptive and hybrid sparse-grid methods for the Vlasov--Poisson--Lenard--Bernstein (VPLB) model.
This model has applications to plasma physics and is simulated in two reduced geometries: a $0x3v$ space homogeneous geometry and a $1x3v$ slab geometry.
We use the discontinuous Galerkin (DG) method as a base discretization due to its high-order accuracy and ability to preserve important structural properties of partial differential equations.  We utilize a multiwavelet basis expansion to determine the sparse-grid basis and the adaptive mesh criteria. 
We analyze the proposed sparse-grid methods on a suite of three test problems by computing the savings afforded by sparse-grids in comparison to standard solutions of the DG method.  
The results are obtained using the adaptive sparse-grid discretization library \asgard.
\end{abstract}


\end{frontmatter}

\clearpage


\nolinenumbers

\section{Introduction}
\label{sec:intro}

In this paper, we investigate sparse-grid solutions to kinetic equations with applications in plasma physics.  In a general setting, equations of this type are defined in terms of a kinetic distribution $f$ that evolves over a six-dimensional phase space (three position and three velocity variables).
To discretize $f$ in phase-space, we work with sparse-grid approximations based on the discontinuous Galerkin (DG) method.  First introduced for kinetic models of radiation transport \cite{reedHill_1973},
the DG method is a finite element method that offers high-order accurate solutions to elliptic \cite{riviere2008discontinuous} and hyperbolic partial differential equations (PDEs) \cite{cockburnShu_2001} with compact stencils.  In addition to being locally conservative, DG methods are efficient at preserving important physical constraints and structural properties inherited from the underlying PDEs that they are used to simulate.  Such properties include positivity \cite{zhangShu_2011}, energy conservation \cite{fu2019optimal,zhang2022energy,xing2013energy}, asymptotic limits \cite{larsen1989asymptotic,adams2001discontinuous,guermond2010asymptotic,sheng2021uniform,xiong2015high,tang2017asymptotic}, entropy stability \cite{chen2020review,friedrich2019entropy,chan2019efficient,yan2023entropy,barth2006discontinuous}, and invariant domains \cite{pazner2021sparse,chu2019realizability,olbrant2012realizability}. For these reasons, the DG method has become a popular tool in the simulation of kinetic equations \cite{hakim_etal_2020,endeveHauck_2022,alekseenko2012application,cheng2014energy,cheng2011brief,decaria2022asymptotic,endeve2015bound,garrett2018fast,hong2022generalized,qiu2011positivity,abdelmalik2016entropy}.

When applied to high-dimensional PDEs, Eulerian grid-based methods, including DG, suffer from the \textit{curse of dimensionality} \cite{bellman1959adaptive}, where the cost to approximate a general measurable function scales like $\mathcal{O}(N^d)$, with $d$ the dimension of the domain and $N$ the degrees of freedom in a single dimension.  Such a scaling in six dimensions makes the standard DG method intractable for approximating general kinetic equations, even on leadership class computing facilities \cite{guo2016sparse}.

Particle-based methods, e.g.~particle in cell, attempt to mitigate the curse of dimensionality using a Lagrangian approach \cite{birdsall2018plasma,chen2011energy,khaziev2018hpic,hu1994generalized,hirvijoki2018metriplectic,carrillo2020particle}. 
However, there has been recent interest in reducing the computational and memory footprint of Eulerian methods by compressing the full-resolution distribution.  
One popular technique is low-rank approximations where the discretized kinetic distribution is treated as a $d$-mode tensor and compressed using a low-rank factorization.
The low-rank decomposition is evolved through time using methods such as step-truncation \cite{rodgers2020step,guo2022low} or dynamical low-rank approximation \cite{einkemmer2018low,einkemmer2021efficient}.
Another popular avenue is the sparse-grid method \cite{guo2017adaptive,tao2019sparse,huang2023adaptive} which is the focus of this paper. 

The sparse-grid method \cite{bungartz2004sparse} is a general technique used for the approximation of high-dimensional functions.
These methods replace the $\mathcal{O}(N^d)$ scaling of tensor-based discretizations to $\mathcal{O}(N (\log N)^{d-1})$.
First developed for the integration of high-dimensional functions \cite{smolyak1963quadrature,gerstner1998numerical}, current flavors of the sparse-grid method are far reaching.
Sparse-grid interpolation has been successfully employed in the construction of surrogate models \cite{Bungartz_Griebel_2004, Chkifa_Cohen_Schwab_2014}
including addressing challenges of adaptivity for basis with local support \cite{pfluger2010spatially, pflueger10spatially, pflueger12spatially, ma2009adaptive, jakeman2012local, stoyanov2018adaptive},
global support
\cite{Stoyanov_Webster_2016, narayan2014adaptive, Nobile_Tamellini_Tempone_2016, morrow2020method}, 
and even discontinuous response surfaces
\cite{jakeman2011characterization, jakeman2013minimal, stoyanov2017predicting}.

Additionally, sparse-grids have gained favor in the approximation of high-dimensional PDEs with examples in finite differences \cite{leentvaar2006pricing}, finite volumes \cite{hemker1995sparse,mishra2012sparse}, conforming finite element methods \cite{schwab2008sparse, shen2010sparse, shen2010efficient, balder1996solution} as well as the DG method \cite{wang2016sparse,guo2016sparse}. 
In the DG context, functions are decomposed in a multiwavelet basis \cite{alpert1993class} with specific basis functions discarded via a sparse-grid selection rule.
This multiwavelet decomposition induces a decay in the coefficient's magnitude over finer levels.
This decay is utilized to build model-independent criteria for adaptively choosing whether to keep or discard basis functions.
This is referred to as the adaptive sparse-grid DG method and has shown promise in the modelling of kinetic equations \cite{kormann2016sparse,guo2017adaptive,tao2019sparse}.

The main goal of this work is to study the computational savings provided by sparse-grids on the Vlasov--Poisson--Lenard--Bernstein (VPLB) model.  We measure the savings by the reduction of the total degrees of freedom required to accurately represent the solution.
The grids of choice are the adaptive sparse-grid DG method and a hybrid sparse-grid method, called the mixed-grid method, which is a standard DG grid in position space tensored with a sparse-grid in velocity space.  
Similar hybrid splittings have been studied in the context of collisionless kinetic problems \cite{kormann2016sparse}.
The methods are tested on the following three problems: a simple relaxation to a Maxwellian equilibrium, the Sod shock tube problem \cite{sod_1978}, and an example of collisional Landau damping \cite{crestetto_etal_2012,hakim_etal_2020,francisquez_etal_2020a}.
In each problem, we present the computational savings achieved as well as general qualitative performance, such as capturing desired physical features, of the methods presented. In general, the adaptive sparse-grid method significantly reduces the storage cost of the distribution while the mixed-grid method only  provides favorable savings in determining lower-order moments of the distribution. 

We work with the VPLB model on a slab geometry which reduces the problem to a four-dimensional $1x3v$ problem (one position dimension, three velocity dimensions). 
This is done so that the problem size is sufficiently small to be run on a single node machine; the Chu reduction method \cite{chu1965kinetic} can be utilized to further reduce the problem to a moment model in $1x1v$ which then allows the creation of fine-resolution reference solutions; and, in a slab geometry, the VPLB model can be written as a sum of terms which are each tensor products of one-dimensional PDE operators. 
The latter property, commnoly referred to as separability, allows efficient evaluations of the model on sparse-grids.   We note that recent endeavors, e.g.~interpolatory wavelets \cite{huang2023adaptive}, are being researched to bypass the separability condition.
Nonetheless, DG approximations to the slab problem with sufficient resolution to accurately capture fine-scale features will still tax a single node machine and thus will require some sort of compression.
Additionally, we are keeping three velocity dimensions which we expect to provide the main source of savings captured by the adaptive sparse-grid DG method; this is because, locally in space, the collision operator mollifies the distribution in the velocity domain, which will in turn cause a rapid decay in the coefficient size when the distribution is represented in the multiwavelet DG basis.  
Moreover, in regimes of high collisionality, where the distribution typically approaches a local thermal equilibrium that is very smooth in the velocity variable, the distribution in physical space is close to a fluid model which can form shocks and other non-smooth features.  Therefore, we believe the $1x3v$ geometry provides an indication of the savings expected in the full $3x3v$ model, as well as clues for constructing hybrid approaches.

Complementing this work is the development of the adaptive sparse-grid DG codebase \asgard~(Adaptive Sparse-Grid Discretization) \cite{asgard2024}.  
The goal of this open-source project is to facilitate and promote the use of adaptive sparse-grid methods for the approximation of kinetic models by providing a robust yet flexible adaptive sparse-grid library.  
All sparse-grid results of this work were computed using \asgard.  
The algorithmic specifics of how \asgard~evaluates PDE operators will be delayed for a future work.
This manuscript only focuses on the mathematics of the adaptive sparse-grid method and the memory reduction realized via its utilization. 

The rest of the paper is organized as follows.
In \Cref{sec:model}, we present the VPLB model, as well as the Chu reduction method for generating reference solutions and the geometric reductions used to formulate the aforementioned test problems.
In \Cref{sec:method}, we present the DG method for the VPLB model which we refer to as the full-grid method.
\Cref{sec:sparsegrid} provides an overview to  the standard and adaptive sparse-grid methods and details the specifics implemented in \asgard.
In \Cref{sec:numerical}, we analyze the results of the adaptive sparse-grid and mixed-grid methods, compared against the full-grid method, for the chosen suite of test problems.
Finally, \Cref{sec:conclusions} gives our conclusions and future plans.
\section{The Vlasov--Poisson--Lenard--Bernstein Model}
\label{sec:model}

The Vlasov--Poisson--Lenard--Bernstein (VPLB) model describes the dynamics of charged particles influenced by a self-consistent electric field and collisional dynamics.  
It couples a kinetic equation for the phase-space distribution function of charged particles with a Poisson equation for the electrostatic potential.  
  
Assuming ions of unit mass and charge, the governing kinetic equation is
\begin{equation}
    \p_{t}f(\bx,\bv,t) + \bv\cdot\nabla_{\bx}f(\bx,\bv,t) + \bE(\bx,t)\cdot\nabla_{\bv}f(\bx,\bv,t) = \cC_{\lb}(f)(\bx,\bv,t),
    \label{eqn:kinetic}
\end{equation}
where the phase-space distribution function $f$ depends on position $\bx = (x_1,x_2,x_3)^\top \in \Omega_{\bx}\subseteq\bbR^{3}$, velocity $\bv = (v_1,v_2,v_3)^\top \in\bbR^{3}$, and time $t \geq 0$.  
The electric field $\bE=-\nabla_{\bx}\Phi$ is obtained from the electrostatic potential $\Phi$ by solving the Poisson equation
\begin{equation}
    -\nabla_{\bx}\cdot\nabla_{\bx}\Phi(\bx,t)
    =n_{f}(\bx,t)-n_{\rm e},
    \label{eqn:poisson}
\end{equation}
where $n_{f}=\vint{f}\equiv\int_{\bbR^{3}}f\dx\bv$ is the ion density, and $n_{\rm e}$ is a constant background electron density chosen to enforce global charge neutrality: $n_e = \int_{\Omega_\bx} n_f(\bx,t)\dx{\bx}$ for all $t \geq 0$.
The collision operator $\cC_{\lb}$ on the right-hand side of \eqref{eqn:kinetic} is the Lenard--Bernstein (LB) operator \cite{lenardBernstein_1958}.  It takes the form \cite{hakim_etal_2020,francisquez_etal_2020a}
\begin{equation}
   \cC_{\lb}[\bsrho_{f}](f)(\bx,\bv,t)
    =\nu\nabla_{\bv}\cdot\big(\,(\bv-\bu_{f})\,f+\theta_{f}\nabla_{\bv}f\,\big),
    \label{eqn:lenardBernstein}
\end{equation}
where $\nu \geq 0$ is a collision frequency that is assumed to be a constant independent of $\bv$, $\bx$, and $t$.  
The moments of $f$, 
\begin{equation}
    \bsrho_{f}=\vint{\be f}
    ,\quad\text{where}\quad
    \be=(e_{0},\be_{1},e_{2})^{\top}
    \equiv (1,\bv,\f{1}{2}|\bv|^{2})^{\top},
\end{equation}
represent the number, momentum, and energy densities, respectively, and the bulk velocity and temperature are defined from $\bsrho_{f}$ by
\begin{equation}\label{eqn:moments_to_fluid_variables}
    \bu_{f}=\f{1}{n_{f}}\vint{f\bv}
    \quad\text{and}\quad
    \theta_{f}=\f{1}{3n_{f}}\vint{f|\bv-\bu_{f}|^{2}}.  
\end{equation}
Direct calculations show that
\begin{equation}
    \bsrho_{f}=\big(\,n_{f},\,n_{f}\bu_{f},\,n_{f}(\frac12|\bu_{f}|^{2}+\frac32\theta_{f})\,\big)^{\top}.
\end{equation}

\begin{prop}[\cite{endeveHauck_2022}]\label{prop:lenardBernsteinOperator}
     The LB operator satisfies the following properties
    \begin{enumerate}
        \item Conservation of number, momentum, and energy:
        \begin{equation}
            \vint{\cC_{\lb}[\bsrho_{w}](w)\be}=0,
            \qquad
            \forall w\in\operatorname{Dom}(\cC_{\lb}).
            \label{eqn:lbConservation}
        \end{equation}
        \item Dissipation of entropy:
        \begin{equation}
            \vint{\cC_{\lb}[\bsrho_{w}](w)\log w}\le0,
            \qquad
            \forall w\in\operatorname{Dom}(\cC_{\lb}).
            \label{eqn:lbDissipation}
        \end{equation}
        \item Characterization of equilibria: For any $w\in\operatorname{Dom}(\cC_{\lb})$,
        \begin{equation}
            \vint{\cC_{\lb}[\bsrho_{w}](w)\log w}=0
            \label{eqn:lbEquilibrium}
        \end{equation}
        if and only if $w$ is a Maxwellian distribution, i.e.,
        \begin{equation}
            w=M_{w}:=\f{n_{w}}{(2\pi\theta_{w})^{3/2}}\exp\Big\{\,\f{-|\bv-\bu_{w}|^{2}}{2\theta_{w}}\,\Big\}.
            \label{eqn:maxwellian}
        \end{equation}
    \end{enumerate}
\end{prop}

\begin{prop}
[\cite{endeveHauck_2022}]
    On a periodic spatial domain $\Omega_{\bx}\subseteq\bbR^{3}$, the VPLB model satisfies the following global conservation laws:
    \begin{enumerate}
        \item Conservation of number:
        \begin{equation}
            \p_{t}\int_{\Omega_{\bx}}\vint{f}\dx\bx=0.
        \end{equation}
        \item Conservation of momentum:
        \begin{equation}
            \p_{t}\int_{\Omega_{\bx}}\vint{\be_{1}f}\dx\bx=0.
        \end{equation}
        \item Conservation of energy:
        \begin{equation}
            \p_{t}\int_{\Omega_{\bx}}(\vint{e_{2}f}+\frac{1}{2}|\bE|^{2})\dx\bx=0.
        \end{equation}
    \end{enumerate}
\end{prop}

\subsection{Geometric reductions}

\subsubsection{Space homogeneous problem} To investigate the relaxation induced by the LB collision operator of a velocity distribution to a Maxwellian, we consider the equation \eqref{eqn:lenardBernstein} under the assumption that $f$ does not depend on $\bx$.  In this case, the PDE is given by 
\begin{align}\label{eqn:relax:pde}
    \partial_t f(\bv,t) = \nu\,\cC_{\lb}(f)(\bv,t).
\end{align} 

\subsubsection{Reduction to slab geometry}

Under the assumption that $\p_y f = \p_z f = 0$, the VPLB model \eqref{eqn:kinetic} reduces to 
\begin{equation}
\label{eqn:kinetic_slab}
    \p_{t}f(x,\bv,t) + v_x \p_x f(x,\bv,t) + E \p_{v_x} f(x,\bv,t) = \cC_{\lb}(f)(x,\bv,t),
\end{equation}
where $E:= E_x = -\p_x \Phi$ and  $\Phi$ satisfies
\begin{equation}
    -\p_{xx} \Phi(x) = n_{f}(x,t)-n_{\rm e}.
    \label{eqn:Poisson_slab}
\end{equation}  
Let $(v_r,\vartheta,\varphi)$ be a spherical-polar coordinate system in which the $x$-axis is aligned with the polar direction, so that
\begin{equation}
    v_x = v_r \cos \vartheta, 
    \quad
    v_y = v_r \sin \vartheta \cos \varphi,
    \quad\text{and}\quad
    v_z = v_r \sin \vartheta \sin \varphi,
\end{equation}
where $v_{r}=|\bv|$, $\vartheta$ is the polar angle, and $\varphi$ is the azimuthal angle.  
We assume further that $f$ is independent of $\varphi$; as a result $(u_f)_y=(u_f)_z=0$, and by abuse of notation we set 
\begin{equation}\label{eqn:u_f_def}
    (u_f)_x := u_f = \frac{\langle f v_x \rangle_v}{n_f}
\end{equation}
so that $\bu_f = [u_f,0,0]^\top$.  The equation \eqref{eqn:kinetic_slab} has a phase space with four total dimensions: one for physical space and three for velocity space, i.e., $1x3v$.

\subsection{\texorpdfstring{Reduction to $1x1v$}{Reduction of 1x1v}}\label{subsec:reduction}
The Chu reduction method is a tool for further reducing the slab geometry problem to $1x1v$, at the cost of solving an additional equation.  It was first developed in \cite{chu1965kinetic} for the Bhatnagar--Gross--Krook (BGK) equation and is used here to provide reference solutions in \Cref{sec:numerical} for sparse-grid simulations when exact solutions are not known and full-grid reference calculations are prohibitively expensive.

To derive the Chu reduction of \eqref{eqn:kinetic_slab}, let 
\begin{gather}\label{eqn:reduction_g_def}
    g_1(x,v_x) = \int_{\bbR^{2}} f(x,\bv) \dx{v_y} \dx{v_z}
    \quand 
    g_2(x,v_x) = \int_{\bbR^{2}} (v_y^2+v_z^2) f(x,\bv) \dx{v_y} \dx{v_z} ,
\end{gather} 
Testing \eqref{eqn:kinetic_slab} by 1 and by $v_y^2+v_z^2$, respectively and integrating over $\dx{v_y} \dx{v_z} $ yields the following coupled system in $(x,v_x)$:
\begin{subequations}
\label{eqn:reduction_system}
\begin{align}
\label{eqn:reduction_system_g1}
    \partial_t g_1 + v_x\partial_{x}g_1 + E\partial_{v_x}g_1 &= \nu\,\mathcal{C}_{1}(g_1;u_f,\theta_f), \\
    \label{eqn:reduction_system_g2}
    \partial_t g_2 + v_x\partial_{x}g_2 + E\partial_{v_x}g_2 &= \nu\,\mathcal{C}_{1}(g_2;u_f,\theta_f) + \nu\,(\,4\theta_f g_1-2g_2\,), 
\end{align}
\end{subequations}
where 
\begin{equation}
    \mathcal{C}_{1}(g;u,\theta) = \partial_{v_x}((v_x-u)g + \theta\partial_{v_x}g),
    \quad 
     E = -\partial_x \Phi,
     \quand
     -\partial_{xx}\Phi = n_f-n_e,
\end{equation}
and, importantly, the velocity moments of $f$ can be expressed in terms of $g_1$ and $g_2$:
\begin{align}\label{eqn:chu_moments}
\begin{split}
    n_f = \int_{\bbR} g_1\dx{v_x},
    \quad 
    u_f = \frac{\int_{\bbR} g_1v_x \,\dx{v_x}}{n_f},
    \quand 
    \theta_f = \frac{1}{3n_f}\int_{\bbR}\big[\, g_1 (v_x-u_f)^2 + g_2 \,\big]\dx{v_x}.
\end{split}
\end{align}
The conservation properties of \eqref{eqn:lbConservation} are preserved; namely,
\begin{subequations}
\begin{align}
    \int_\bbR \mathcal{C}_{1}(\tilde{g}_1;u_f,\theta_f) \dx{v_x} = \int_\bbR \mathcal{C}_{1}(\tilde{g}_1;u_f,\theta_f)v_x \dx{v_x} &= 0, \\
    \frac{1}{2}\int_\bbR \mathcal{C}_{1}(\tilde{g}_1;u_f,\theta_f)v_x^2 + \mathcal{C}_{1}(\tilde{g}_2;u_f,\theta_f) + (4\theta_f \tilde{g}_1-2\tilde{g}_2)\dx{v_x} &= 0, 
\end{align}
\end{subequations}
for any $\tilde{g}_1$ and $\tilde{g}_2$ such that the fluid variables $n_f$, $u_f$, and $\theta_f$ are built via \eqref{eqn:chu_moments} using $\tilde{g}_1$ and $\tilde{g}_2$.

Though not required for \eqref{eqn:reduction_system}, we will, for diagnostic purposes in \Cref{sec:numerical}, also consider the function
\begin{equation}
    g_3(x,v_x) = \int_{\bbR^{2}} (v_y^4+v_z^4) f(x,\bv) \dx{v_y} \dx{v_z} ,
\end{equation}
which satisfies
\begin{equation}
\label{eqn:g3}
      \partial_t g_3 + v_x\partial_{x}g_3 + E\partial_{v_x}g_3  
      = \nu\,\mathcal{C}_{1}(g_3;u_f,\theta_f) + \nu\,(\,12\theta_{f} g_2-4g_3\,). 
\end{equation}

\section{Notation and the Discontinuous Galerkin Method}
\label{sec:method}

\subsection{Notation}\label{subsec:notation}
Let $\ell_x\in\bbN_0 = \{0,1,2,\ldots\}$, $\W_x=(-L_x,L_x)$ be an interval in physical space, and $\cT_{x,\ell_x}$ be a uniform mesh on $\W_x$ with $2^{\ell_x}$ elements.  
Let $\cE_{x,\ell_x}$ be the skeleton of $\cT_{x,\ell_x}$.  

Similarly, let $\ell_v\in\bbN_0$, $\W_v=(-L_v,L_v)^3 \subset \bbR^3$, and $\cT_{v,\ell_v}$ be a uniform cubic mesh on $\W_v$ with $2^{\ell_v}$ elements in each dimension.  
Let $\cE_{v,\ell_v}^{\text{I}}$ be the interior (i.e., not including boundaries) skeleton on this mesh.  
We will often use $\langle\,\cdot\,\rangle_v$ and $\langle\,\cdot\,\rangle_{v_y,v_z}$ to denote integration in $\dx{\bv}$ and $\dx{v_y}\dx{v_z}$, respectively. 

We let $\W=\W_x\times\W_v\subset\bbR^4$, and denote $L^2(\W)$ and $H^s(\W)$ to be the standard Lebesgue and Sobolev spaces on $\W$.  
Let $(\cdot\,,\cdot)$ be the $L^2(\W)$-inner product with norm $\|\cdot\|_{L^2(\W)}$ and let $\|\cdot\|_{H^s(\W)}$ be the norm on $H^s(\W)$.  We denote by $L^2(D)$ and $(\cdot\,,\cdot)_D$ the $L^2$ space with standard inner product on some domain $D$ which is typically $\W_x$ or $\W_v$.  Any of the inner products mentioned above can be trivially extended to vector-valued functions with the standard Euclidean inner product.  

Denote the discontinuous Galerkin finite element spaces $V_{x,\ell_{x}}\subset L^2(\W_x)$ and $V_{v,\ell_{v}}\subset L^2(\W_v)$ by 
\begin{equation}\label{eqn:DG_spaces}
\begin{split}
    V_{x,\ell_x} &= \{g\in L^2(\W_x):g\big|_{K} = \bbQ_k(K)~\forall K\in\cT_{x,\ell_x}\} \\
    V_{v,\ell_v} &= \{g\in L^2(\W_v):g\big|_{K} = \bbQ_k(K)~\forall K\in\cT_{v,\ell_v}\}
\end{split}
\end{equation}
where $\bbQ_k(K)$ is the set of all polynomials of maximum degree $k$ in any direction on $K$.  We assume $k=2$ unless written otherwise.
Let $\mathcal{V}_\ell=V_{x,\ell_x}\otimes V_{v,\ell_v}$.  

Given $x_*\in\cE_{x,\ell_x}$, let $g$ be a function with traces $g^\pm(x_*):=\lim_{x\to x_*^\pm}g(x)$ are well defined. 
Define the average and jump of $g$ in $x$, respectively, by
\begin{equation}\label{eqn:avg_jmp_x}
    \avg{g} = \tfrac{1}{2}(g^++g^-)
    \qquand
    \ljmp{g}\rjmp = g^--g^+.
\end{equation}
We account for the periodic boundary in $\cE_{x,\ell_x}$ by defining the jumps and averages on the boundary using \eqref{eqn:avg_jmp_x} with $g^+ = g(L_x)$ and $g^- = g(-L_x)$.
We denote by $S_{x,\ell_x}$ the intersection of $V_{x,\ell_x}$ with $k=1$ and all continuous and periodic functions on the closure of $\W_x$, i.e.~$\overline{\W_x}$.  The space $S_{x,\ell_x}$ is used for the discretization of \eqref{eqn:Poisson_slab} and uses linear functions so that the electric field $E$ is constant on each element.

Similarly, consider the edge $e\in\cE_{v,\ell_v}^{\text{I}}$, where $e=\partial K^+\cap \partial K^-$ and $K^\pm\in\cT_{v,\ell_v}$ with normal outward vector $\mathbf{n}_v^\pm$. Given a scalar and vector valued function $g$ and $\bm{\sigma}$ respectively with well defined traces on $\partial K^{\pm}$, define the average and jump of $g$ and $\bm{\sigma}$ in $\bv$, respectively, by
\begin{equation}\label{eqn:avg_jmp_v}
\begin{aligned}
    \avg{g} &= \tfrac{1}{2}(g^++g^-)&\text{and}&&\ljmp{g}\rjmp &= g^-\mathbf{n}_v^- + g^+\mathbf{n}_v^+, \\
    \avg{\bm{\sigma}} &= \tfrac{1}{2}(\bm{\sigma}^++\bm{\sigma}^-)&\text{and}&&\ljmp{\bm{\sigma}}\rjmp &= \bm{\sigma}^-\cdot\mathbf{n}_v^- + \bm{\sigma}^+\cdot\mathbf{n}_v^+,
\end{aligned}
\end{equation}
where for any $\bv_*\in e$,
\begin{equation}
g^\pm(\bv_*) = \lim_{\substack{\bv\to\bv_* \\ \bv\in K^\pm}}g(\bv)
\end{equation}
with analogous definition for $\sigma^{\pm}$.
While the same notation for average and jumps is used in the physical and velocity domains, the domain of integration of the DG formulation provides context to which case is used (see \eqref{eqn:VP_discrete_def}).
Let $\ledge\cdot\,,\cdot\redge_e$ be the $L^2$ inner product over an edge $e$ and denote
$\ledge\cdot\,,\cdot\redge_{\cE_{x,\ell_x}}=\sum_{e\in\cE_{x,\ell_x}}\ledge\cdot,\cdot\redge_e$ with an analogous definition for $\ledge\cdot\redge_{\cE_{v,\ell_v}^{\text{I}}}$.  For functions $g$ in $\mathcal{V}_\ell$, let $\partial_x$ and $\grad_{\bv}$ represent the piece-wise spatial derivative and velocity gradient $g$.

Finally, for time integration, let $\Delta t>0$ be the timestep, assumed for our purposes to be uniform.  For $\mfn\in\bbN_0$ define $t^\mfn=\mfn\Delta t$ and denote $f^\mfn$ to be an approximation to $f(t^\mfn)$.

\subsection{Discontinuous Galerkin Method}\label{subsec:discretization}

We first discretize \eqref{eqn:kinetic} in phase space on $\mathcal{V}_\ell$ by the following semi-discrete problem:  Find $f_h\in C([0,\infty];\mathcal{V}_\ell)$ such that

\begin{equation}\label{eqn:discrete_form}
    (\partial_t f_h,g_h) + \mA_{\text{VP}}(f_h,g_h) = \nu\mA_{\text{LB}}(f_h,\bm{\rho}_{f_h},g_h)
\end{equation}
holds for all $g_h\in \mathcal{V}_\ell$.  
The Vlasov--Poisson portion, $\mA_{\text{VP}}$, is discretized with upwind fluxes; specifically, 
\begin{align}\label{eqn:VP_discrete_def}
\begin{split}
    \AVP(w_h,g_h) &= -(v_xw_h,\partial_xg_h) + \ledge v_x\avg{w_h}+\tfrac{|v_x|}{2}\jmp{w_h},\jmp{g_h}\redge_{\cE_{x,\ell_x}\times\W_v} \\
    &\quad-(\tilde{\bm{E}}_hw_h,\grad_{\bv} g_h) + \ledge\avg{\tilde{\bm{E}}_hw_h}+\tfrac{|\tilde{\bm{E}}_h\cdot\textbf{n}_v|}{2}\jmp{w_h}, \jmp{g_h}\redge_{\W_x\times\cE_{v,\ell_v}^{\text{I}}}
\end{split}
\end{align}
for all $w_h,g_h\in \mathcal{V}_\ell$ where $\tilde{\bm{E}}_h:=(E_h,0,0)^\top$  and $E_h$ is given by $-\partial_x\Phi_h$ where $\Phi_h\in S_{x,\ell_x}$ satisfies
\begin{equation}\label{eqn:poisson_solve}
    -(\partial_x\Phi_h,\partial_x q_h)_{\W_x} = \big(\big<w_h\big>_v-n_e,q_h\big)_{\W_x}
\end{equation}
for every $q_h\in S_{x,\ell_x}$.
The boundary conditions are periodic in $x$ and we impose zero fluxes on the velocity boundaries.

The Lenard--Bernstein portion, $\mA_{\text{LB}}$, of \Cref{eqn:discrete_form} is discretized with the LDG method (e.g., \cite{cockburnShu_2001}), with central fluxes for the diffusion term and a local Lax--Friedrichs flux for the advection term; namely,
\begin{align}\label{eqn:LB_discrete_def}
\begin{split}
    \ALB(w_h,\bm{\rho}_h,g_h) &= -((\bv-\mathbf{u})w_h,\partial_v g_h) + \ledge \lavg \bv w_h\ravg-\tfrac{|\bv\cdot\textbf{n}_v|}{2}\jmp{w_h},\jmp{g_h}\redge_{\W_x\times\cE_{v,\ell_v}^{\text{I}}} \\
    &\quad -(\bm{\sigma}_h,\grad_{\bv} g_h) + \ledge\avg{\bm{\sigma}_h},\jmp{g_h}\redge_{\W_x\times\cE_{v,\ell_v}^{\text{I}}}
\end{split}
\end{align}
for every $w_h,g_h\in \mathcal{V}_\ell$, where $\bu$  is determined from $\bm{\rho}_h\in [V_{x,\ell_x}]^3$ via formulas in \eqref{eqn:moments_to_fluid_variables}.
Here $\bm{\sigma}_h\in [\mathcal{V}_\ell]^3$ is the approximation to the velocity gradient of $w_h$ and is defined by
\begin{equation}\label{eqn:LB_discrete_sigma_def}
    (\bm{\sigma}_h,\tau_h) = (\theta \grad_{\bv} w_h,\tau_h) - \ledge\theta\jmp{w_h},\avg{\tau_h}\redge_{\W_x\times\cE_{v,\ell_v}^{\text{I}}}
\end{equation}
for every $\tau_h\in[\mathcal{V}_\ell]^3$, where $\theta$ is determined by $\bm{\rho}_h$ the relevant formula in \eqref{eqn:moments_to_fluid_variables}.

If $w_h=0$ on $\partial\W_v$, then it can be shown that $\mA_{\text{LB}}(w_h,\bsrho_{w_h},\be q_h) = 0$ for all $q_h\in V_{x,\ell_x}$, which implies that the conservation properties in \Cref{eqn:lbConservation} hold.  

For brevity, we do not provide the discretization for the Chu reduction \eqref{eqn:reduction_system}, but we note it is similar to the discretizations given above for the slab problem.  

\subsection{Time Stepping Method}\label{subsec:time_stepping}

We discretize \eqref{eqn:discrete_form} in time via Implicit-Explicit (IMEX) Runge--Kutta (RK) methods \cite{ascher1997implicit}.  
Such methods are popular time steppers for evolving kinetic models that feature multiple time scales \cite{pareschiRusso_2005,chu_etal_2019,endeveHauck_2022}.  
In our case, the Vlasov--Poisson portion $\mA_{\text{VP}}$ will be evolved explicitly and the collision operator $\mA_{\text{LB}}$ will be evolved implicitly.  
We will use IMEX-RK method of \cite{chu_etal_2019} which is given by:
\begin{subequations}\label{eqn:IMEX_RK}
\begin{align}
(f_h^{(1,*)},g_h) &= (f_h^{\mathfrak{n}},g_h) - \Delta{t}\AVP(f_h^{\mathfrak{n}},g_h) \label{eqn:IMEX_RK:ex1}, \\
(f_h^{(1)},g_h)  &= (f_h^{(1,*)},g_h) + \Delta t\nu\ALB(f_h^{(1)},\bm{\rho}_{f_h^{(1*)}},g_h), \label{eqn:IMEX_RK:im1}\\ 
(f_h^{(2,*)},g_h) &= \tfrac{1}{2}(f_h^{\mathfrak{n}},g_h) - \tfrac{1}{2}\big((f_h^{(1)},g_h)-\Delta{t}\AVP(f_h^{(1)},g_h)\big), \label{eqn:IMEX_RK:ex2} \\
(f_h^{(2)},g_h) &= (f_h^{(2,*)},g_h) + \tfrac{1}{2}\Delta t\nu\ALB(f_h^{(2)},\bm{\rho}_{f_h^{(2*)}},g_h), \label{eqn:IMEX_RK:im2}
\end{align}
and $f_h^{\mathfrak{n}+1}:=f_h^{(2)}$
\end{subequations}

Assuming zero velocity-boundary data, the invariance of the discrete collision operator implies $\bm{\rho}_{f_h^{(s,*)}} = \bm{\rho}_{f_h^{(s)}}$, for $s \in \{1,2\}$.  
Therefore we plug the moments $f_h^{(s,*)}$ into the collision operator $\mA_{\text{LB}}$ in \eqref{eqn:IMEX_RK:im1} and \eqref{eqn:IMEX_RK:im2} for $s=1$ and $s=2$ respectively.  
This decouples the moments from the distribution and provides a linear solve for $f_h^{(s)}$.  Both \eqref{eqn:IMEX_RK:im1} and \eqref{eqn:IMEX_RK:im2} are solved iteratively using GMRES with the possible inclusion of a block-Jacobi preconditioner.
\section{Sparse-grid Method}
\label{sec:sparsegrid}

In this section, we describe the sparse-grid DG method and adaptivity procedure used in \asgard.  
The method, first introduced in \cite{wang2016sparse} (and from which some of the presentation of this section is based), is provided here for completeness.  
We first construct the wavelet basis in one dimension, then extend to multiple dimensions and introduce the sparse-grid selection rule, and finally discuss the adaptivity procedure.  

\subsection{Single Dimension Wavelet Basis}\label{subsec:1D_basis}

The one-dimensional wavelet basis is a hierarchical basis in which additional basis functions for resolving fine scale features are introduced using orthogonal complements to current functions in the basis.  
To simplify the presentation, we assume a domain $\W=[0,1]$. Given a level $\ell\in\{0,\ldots,N\}$, let $\mT_\ell$ be a uniform mesh of $\W$ with mesh size $h_\ell=2^{-\ell}$.  
The partition of $\mT_\ell$ is characterized by the union of disjoint intervals $I_{\ell,j}:=(2^{-\ell}j,2^{-\ell}(j+1))$ for $j=0,\ldots,2^\ell-1$.  
Given this mesh, define the corresponding DG finite element space $V_\ell$ by\footnote{We will often drop the polynomial degree superscript on $V_{\ell}^{k}$ for brevity.}
\begin{equation}\label{eqn:1D_DG_space}
    V_\ell := V_\ell^k = \left\{ g\in L^2(\W) : g\big|_{I_{\ell,j}} \in \mathbb{P}_k(I_{\ell,j})~\forall j=0,\ldots,2^\ell-1\right\},
\end{equation}
where $\mathbb{P}_k$ is the space of polynomials of degree up to $k$.  This space has dimension $\mathrm{dim}(V_\ell) = 2^\ell(k+1)$.  
Additionally, due to the uniform partitioning, 
\begin{equation}
    V_0 \subset V_1\subset V_2 \subset \cdots \subset V_N .
\end{equation}

Let $W_\ell$ to be the orthogonal complement of $V_{\ell-1}$ in $V_\ell$ with respect to the $L^2(\W)$ inner product; that is, $W_0 = V_0$, while for $\ell \geq 1$, 
\begin{equation}\label{eqn:W_def}
V_\ell = V_{\ell-1}\oplus W_\ell 
\qquand
W_\ell \perp V_{\ell-1},
\end{equation}
where $\oplus$ is the direct sum and $\mathrm{dim}(W_\ell) = \max\{0,2^{\ell-1}(k+1)\}$.  Then 
\begin{equation}\label{eqn:hierarchical_decomp_1d}
    V_N = \bigoplus_{\ell=0}^N W_\ell.
\end{equation}

The hierarchical decomposition in \eqref{eqn:hierarchical_decomp_1d} induces a natural decay in the coefficients for the approximation of smooth functions.  
Specifically, let $Q_\ell:L^2(\W)\to W_\ell$ be the orthogonal $L^2$ projection onto $W_\ell$.  
Then by standard polynomial approximation theory (see, e.g., \cite[Section 5.4.2]{canuto2007spectral} or \cite[Theorem 2.6]{riviere2008discontinuous}), there exists a constant $C>0$, independent of $\ell$, such that for any $g\in H^s(\W)$,
\begin{equation}\label{eqn:hier_approx_result_1d}
    \|Q_\ell g\|_{L^2(\W)} \leq C h_\ell^{\min\{s,k+1\}}\|g\|_{H^s(\W)}.
\end{equation}
This decay property motivates the adaptive strategy described in  \Cref{subsec:adaptive_sparse_grids}.

A standard choice for the basis of $W_\ell$ for $\ell\geq 1$ are \textit{wavelets} -- functions that are scaled and shifted to capture finer-scale features.  
The prototype wavelet is the piece-wise constant Haar basis \cite{haar1909theorie}. Here we use Alpert wavelets \cite{alpert1993class}.

\begin{defn}\label{defn:alpert}
The Alpert wavelets are a set of a functions $\{\phi_i(y):i=1,\ldots,k+1\}\subset L^2(\mathbb{R})$ with support in $[-1,1]$ and defined such that
\begin{enumerate}
    \item $\phi_i\big|_{(0,1)}\in\mathbb{P}_k(0,1)$.
    \item $\phi_i(y) = (-1)^{i+k}\phi_i(-y)$.
    \item $\int_{-1}^1\phi_i(y) y^j\dx{y} = 0$ for all $j=0,1,\ldots,i+k-1$. \label{enum:polynomial_orthogonality}
    \item $\int_{-1}^1 \phi_i(y)\phi_j(y) \dx{y} = \delta_{ij}$ for all $i,j=1,\ldots,k$ where $\delta_{ij}$ is the Kronecker delta.
\end{enumerate}
\end{defn}

For a given polynomial degree $k$, the Alpert wavelets satisfying \Cref{defn:alpert} are unique up to a sign.  
The Alpert basis is not hierarchical in the polynomial degree; thus each wavelet must be reconstructed when $k$ is changed.  
For $k=0$, Alpert's wavelets correspond to the Haar basis.  For $k=2$, the wavelets are given on the interval $(0,1)$ by 
\begin{equation}\label{eqn:quad_wavelet_basis}
    \phi_1(y) = \tfrac{1}{3}\sqrt{\tfrac{1}{2}}(1-24y+30y^2),\quad 
    \phi_2(y) = \tfrac{1}{2}\sqrt{\tfrac{3}{2}}(3-16y+15y^2),\quad
    \phi_3(y) = \tfrac{1}{3}\sqrt{\tfrac{5}{2}}(4-15y+12y^2)
\end{equation}
Construction of the wavelets and examples for other polynomial degrees can be found in \cite[Page 5]{alpert1993class}. 

\begin{figure}[!ht]
    \centering
    \begin{subfigure}[b]{0.24\textwidth}
        \centering
        \includegraphics[width=\textwidth]{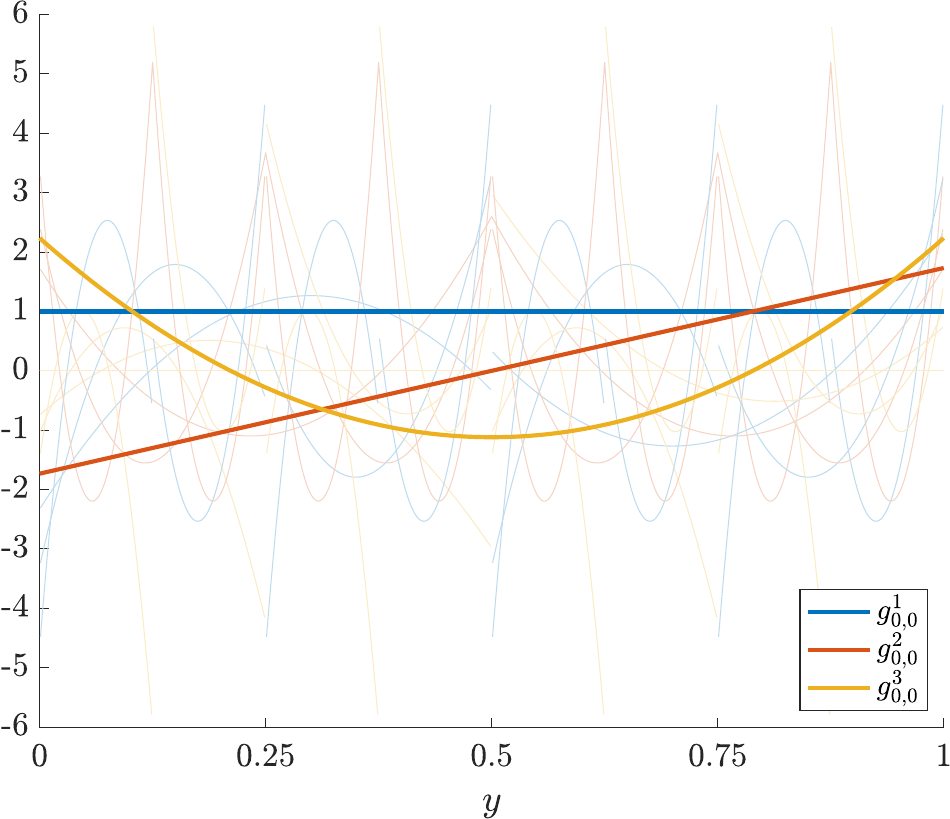}
        \caption{$\ell=0,j=0$}
    \end{subfigure}
    \begin{subfigure}[b]{0.24\textwidth}
        \centering
        \includegraphics[width=\textwidth]{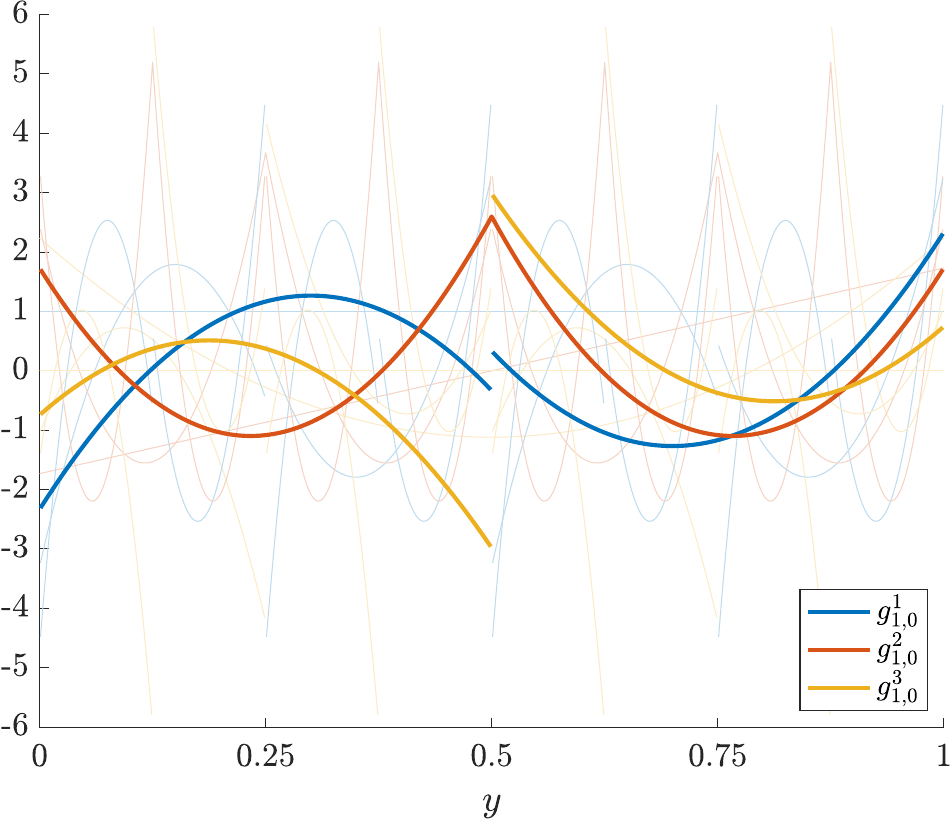}
        \caption{$\ell=1,j=0$}
    \end{subfigure}
    \begin{subfigure}[b]{0.24\textwidth}
        \centering
        \includegraphics[width=\textwidth]{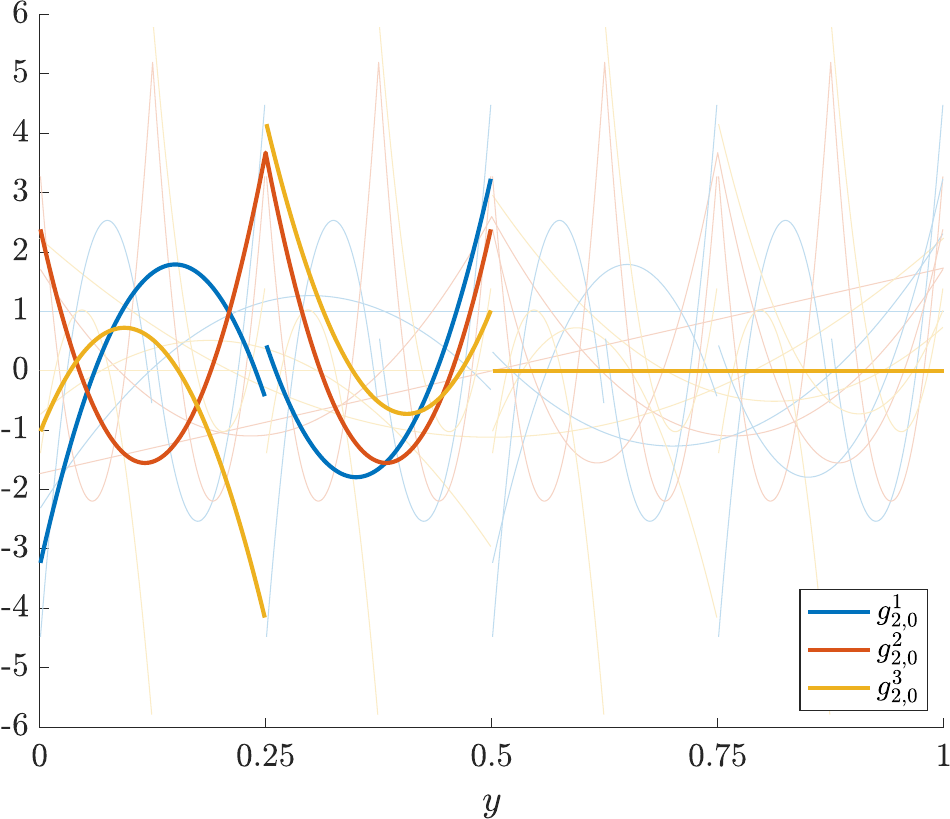}
        \caption{$\ell=2,j=0$}
    \end{subfigure}
    \begin{subfigure}[b]{0.24\textwidth}
        \centering
        \includegraphics[width=\textwidth]{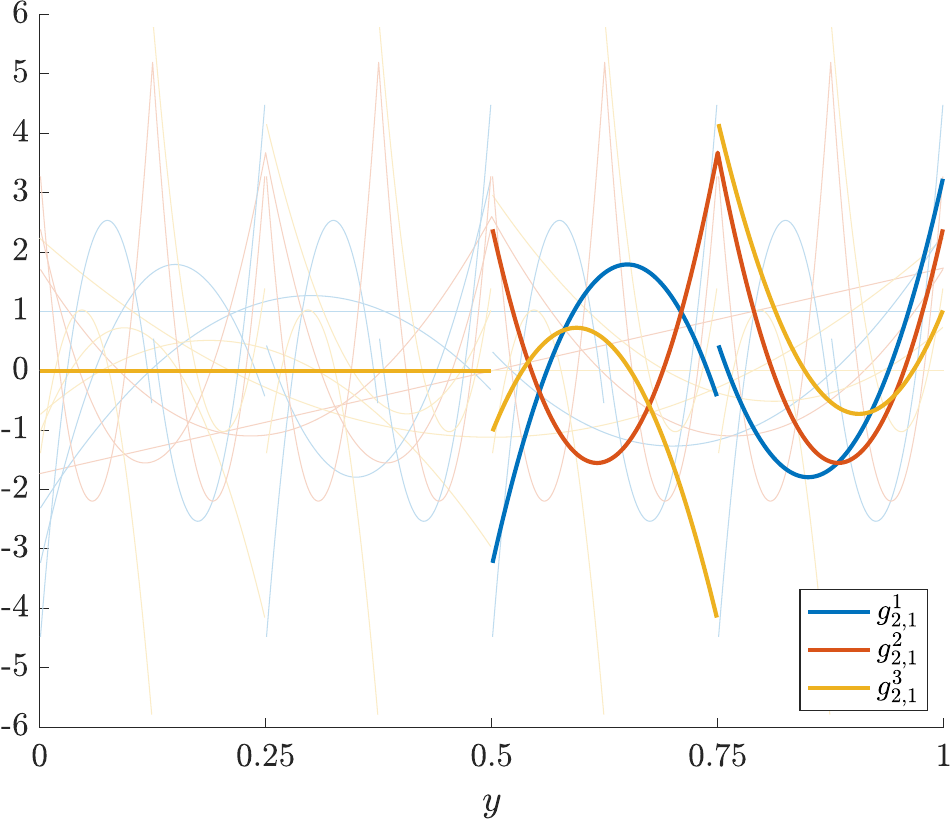}
        \caption{$\ell=2,j=1$}
    \end{subfigure}

    \begin{subfigure}[b]{0.24\textwidth}
        \centering
        \includegraphics[width=\textwidth]{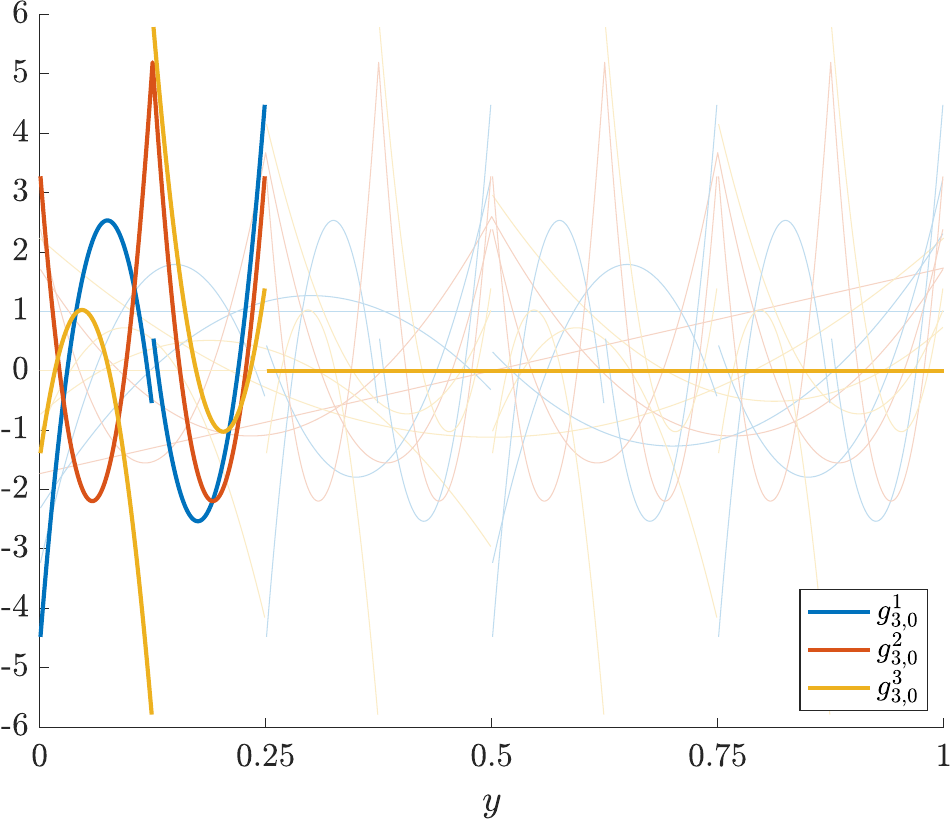}
        \caption{$\ell=3,j=0$}
    \end{subfigure}
    \begin{subfigure}[b]{0.24\textwidth}
        \centering
        \includegraphics[width=\textwidth]{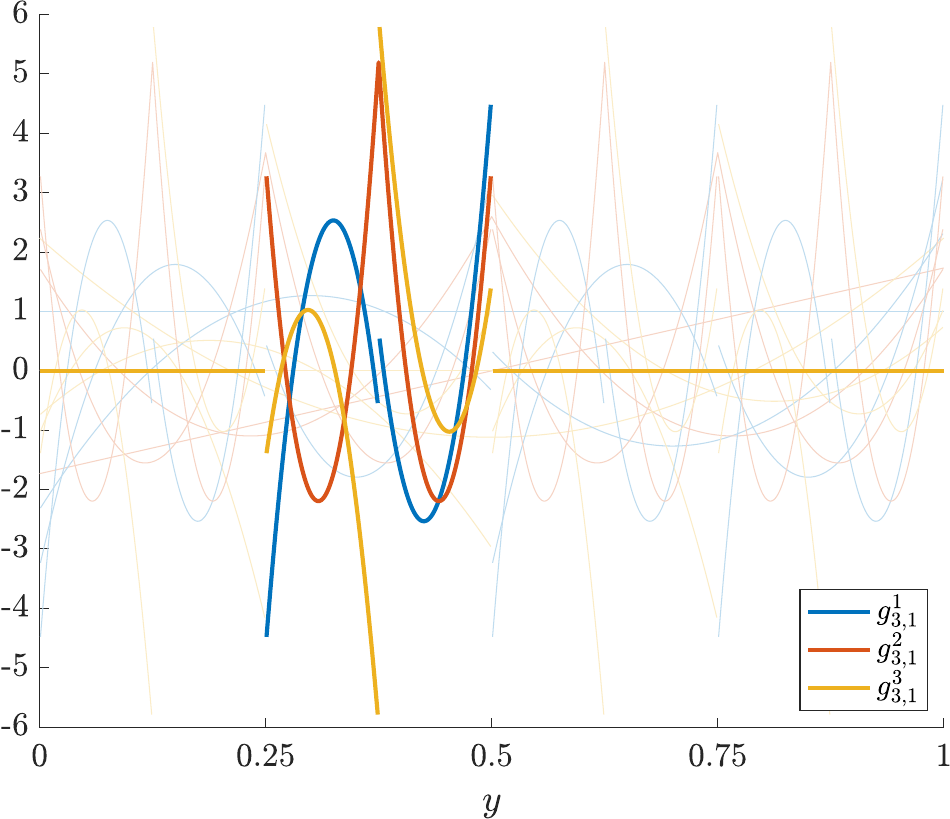}
        \caption{$\ell=3,j=1$}
    \end{subfigure}
    \begin{subfigure}[b]{0.24\textwidth}
        \centering
        \includegraphics[width=\textwidth]{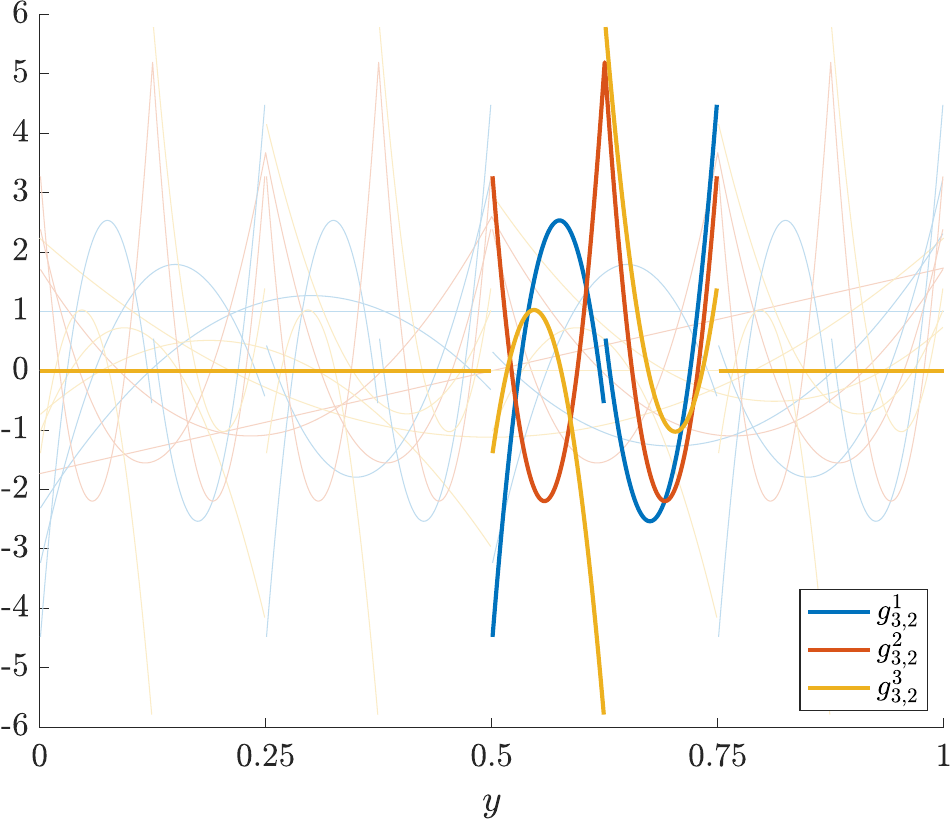}
        \caption{$\ell=3,j=2$}
    \end{subfigure}
    \begin{subfigure}[b]{0.24\textwidth}
        \centering
        \includegraphics[width=\textwidth]{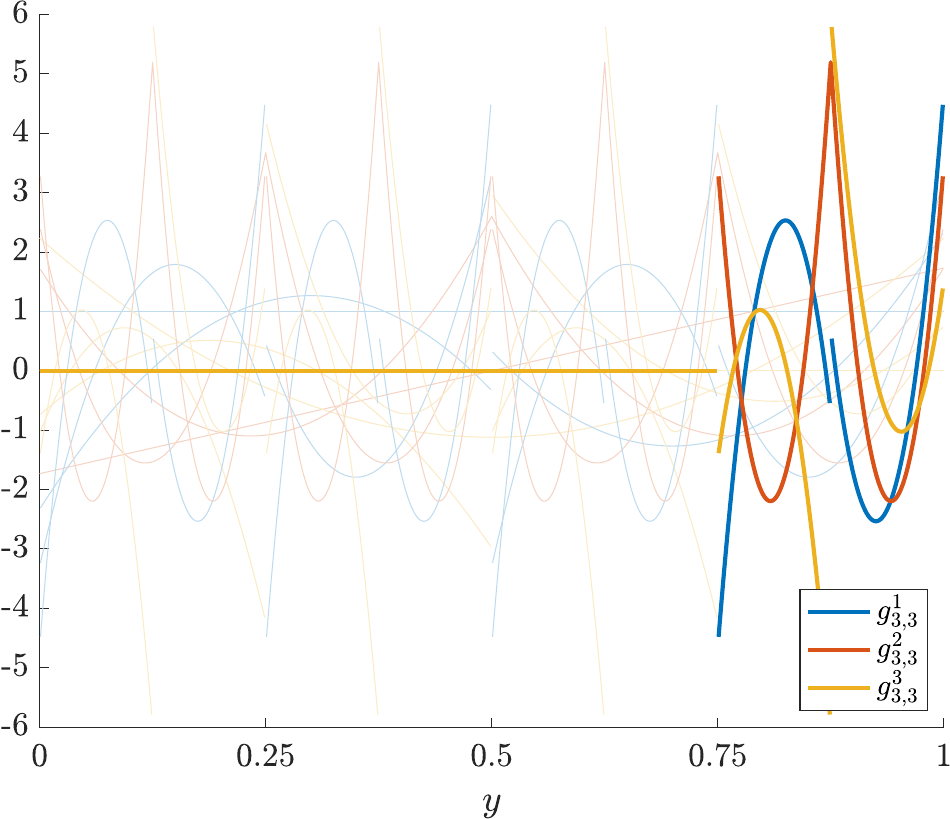}
        \caption{$\ell=3,j=3$}
    \end{subfigure}
    \caption{Plots of the wavelet basis $g_{\ell,j}^i$, given by \eqref{eqn:wavelet_basis_1D_def}, for $k=2$.  In each plot, the entire set of wavelet basis functions for level $\ell=3$ and lower are shown in each plot and are translucent.}
    \label{fig:wavelets_1d}
\end{figure}

For each $\ell \geq 0$, we use the Alpert wavelets to define a basis set $\{g_{\ell,j}^{i}\}$ of $W_\ell$. For $\ell=0$, we choose $g_{0,0}^i$ to be the shifted Legendre polynomials normalized on $L^2(\W)$.  For $\ell \geq 1$, we shift and rescale the Alpert wavelets so that for each $x\in(0,1)$,
\begin{equation}\label{eqn:wavelet_basis_1D_def}
    g_{\ell,j}^i(y) = 2^{(\ell-1)/2}\gamma_i(2^{\ell-1}y-j) 
    ,\quad \text{where}
    \quad
    \gamma_i(y):=\sqrt{2}\phi_i(2y-1).
\end{equation}
Here $\ell$ is the level, $j=0,\ldots,2^{\ell-1}-1$ is the level index, and $i=1,\ldots,k+1$ is the polynomial index.  
The support of $g_{\ell,j}^i$ is precisely $I_{\ell-1,\lfloor j/2\rfloor}$, where $\lfloor\cdot\rfloor$ is the floor function.
Additionally, since every wavelet $g_{\ell',j'}^{i'}$ for any $i',j',$ and $\ell'<\ell$ is a polynomial on $I_{\ell-1,\lfloor j/2\rfloor}$, 
\Cref{enum:polynomial_orthogonality} of \Cref{defn:alpert} ensures that the wavelet bases are all orthonormal; that is,
\begin{equation}\label{eqn:wavelet_orthonormality}
    \int_0^1 g_{\ell,j}^i(y)g_{\ell',j'}^{i'}(y) \dx{y} = \delta_{ii'}\delta_{\ell\ell'}\delta_{jj'}.
\end{equation}
Plots of the wavelets $g_{\ell,j}^i$ for $\ell=0,1,2,3$ and $k=2$ are given in \Cref{fig:wavelets_1d}.

\subsection{Multiwavelets}\label{subsec:nD_basis}

A $d$-dimensional basis is achieved through a tensor product extension.  
Let $\W^{d}=(0,1)^d$ with $\bm{y}=(y_1,\ldots,y_d)\in\W^{d}$.  
Given a multi-index $\bm{\alpha}=(\alpha_1,\ldots,\alpha_d)\in\mathbb{N}_0^d$, define the norms
\begin{equation}\label{eqn:multi_index_norms}
    |\bm{\alpha}|_1 = \sum_{m=1}^d\alpha_m
    \qquand 
    |\bm{\alpha}|_\infty = \max_{1\leq m\leq d}\alpha_m.
\end{equation}  

Let $\bm{\ell}=(\ell_1,\ldots,\ell_d)$ be a multi-index level set, where $\ell_d$ defines the level for dimension $d$, and let $\mathcal{T}_{\bm{\ell}}$ be a tensor product mesh with multi-dimensional mesh parameter $\bm{h}:=(2^{-\ell_1},\ldots,2^{-\ell_d})$.
We label all elements in $\mathcal{T}_{\bm{\ell}}$ by $I_{\bm{\ell},\bm{j}} = \{\bm{y}:y_m\in(2^{-\ell_m}j_m,2^{-\ell_m}(j_m+1)\}$ and define the tensor product finite element space by
\begin{equation}
    \bm{V_\ell}:=\bm{V_\ell}^k=\{g\in L^2(\W):g\big|_{I_{\bm{\ell},\bm{j}}}\in\mathbb{Q}_k(I_{\bm{\ell},\bm{j}}),~\forall~0\leq j_m\leq 2^{\ell_m}-1, m=1,\ldots,d\},
\end{equation}
where $\mathbb{Q}_k(I_{\bm{\ell},\bm{j}})$ represents the set of polynomials of degree up to $k$ in each dimension on $I_{\bm{\ell},\bm{j}}$.  If $\bm{\ell}=(N,\ldots,N)$, then we abbreviate $\bm{V_\ell}$ by $\bm{V}_N$.

Recall the one-dimensional hierarchical decomposition in \Cref{subsec:1D_basis}.  Given the complementary sets $W_{\ell_m}$ defined in \eqref{eqn:W_def}, let
\begin{equation}
    \bm{W_\ell} = W_{\ell_1}\otimes W_{\ell_2}\otimes\cdots\otimes W_{\ell_d}.
\end{equation}
Then \eqref{eqn:hierarchical_decomp_1d} extends to the multidimensional setting:
\begin{equation}\label{eqn:hierarchical_decomp_nd}
    \bm{V_\ell} = V_{\ell_1}\otimes\cdots\otimes V_{\ell_d} = \bigoplus_{\bm{0}\leq\bm{\ell'}\leq\bm{\ell}}\bm{W_{\ell'}}.
\end{equation}
An extension of the coefficient decay result \eqref{eqn:hier_approx_result_1d} also holds.  Let $\bm{Q_\ell}:L^2(\W)\to \bm{W_\ell}$ be the orthogonal $L^2$ projection onto $\bm{W_\ell}$, then
\begin{equation}\label{eqn:hier_approx_result_nd}
    \|\bm{Q_\ell}g\|_{L^2(\W)} = \mathcal{O}\left(\prod_{m=1}^d h_m^{\min\{s,k+1\}}\right),
\end{equation}
where $s$ is a regularity parameter tied to a Sobolev-like space including high-order mixed derivative control.  We refer the reader to \cite[(A.8)]{guo2016sparse} and \cite[Proposition 5.1]{schwab2008sparse} for specifics on \eqref{eqn:hier_approx_result_nd}.

The basis we choose for $\bm{W_\ell}$ are the \textit{multiwavelets} which are products of the 1D wavelets in \eqref{eqn:wavelet_basis_1D_def}:
\begin{equation}\label{eqn:wavelet_basis_nD_def}
    g_{\bm{\ell},\bm{j}}^{\bm{i}}(\bm{y}) := \prod_{m=1}^d g_{\ell_m,j_m}^{i_m}(y_m),~\text{where}~j_m=0,\ldots,\max\{0,2^{\ell_m-1}-1\},i_m=1,\ldots,k+1.
\end{equation}
It follows from repeated application of \eqref{eqn:wavelet_orthonormality} in each dimension that these multiwavelets are orthonormal in $L^2(\W)$ .

\subsection{The Sparse-grid Selection Rule}\label{subsec:sparse_grid_rule}

The spaces $\bm{W}_\ell$ are used to define the sparse grid.  
From \eqref{eqn:hierarchical_decomp_nd} we can rewrite the full-grid as
\begin{equation}\label{eqn:full_grid_def}
    \bm{V}_N = \bigoplus_{|\bm{\ell}|_\infty\leq N}\bm{W_\ell}.
\end{equation}
This space has dimension $\mathrm{dim}(\bm{V}_N) = (k+1)^d2^{Nd}$.
The sparse grid is defined via a selection rule that relaxes the index norm in \eqref{eqn:full_grid_def}.
\begin{defn}[\cite{wang2016sparse,Bungartz_Griebel_2004}]\label{defn:sparse_grid}
    The level $N$ sparse grid, $\hat{\bm{V}}_N \subseteq \bm{V}_N$, is defined by
    \begin{equation}\label{eqn:sparse_grid_def}
        \hat{\bm{V}}_N = \bigoplus_{|\bm{\ell}|_1\leq N}\bm{W_\ell}.
    \end{equation}    
\end{defn}

By definition, the sparse-grid only includes components $\bm{W_\ell}$ whose level indices $\bm{\ell}$ sum up to $N$, and throws away basis functions deemed too fine to include in multiple dimensions.  It was shown in \cite[Lemma 2.3]{wang2016sparse} that
\begin{equation}
\textrm{dim}(\hat{\bm{V}}_N)=\Theta((k+1)^d 2^N N^{d-1}),
\end{equation}
which avoids the costly $\mathcal{O}(2^{Nd})$ scaling of the full-grid in \eqref{eqn:full_grid_def} but still maintains exponential dependence on $k$ and on $\log(N)$.

\begin{figure}[ht]
    \centering
    \begin{subfigure}[b]{0.45\textwidth}
        \includegraphics[height=2.25in]{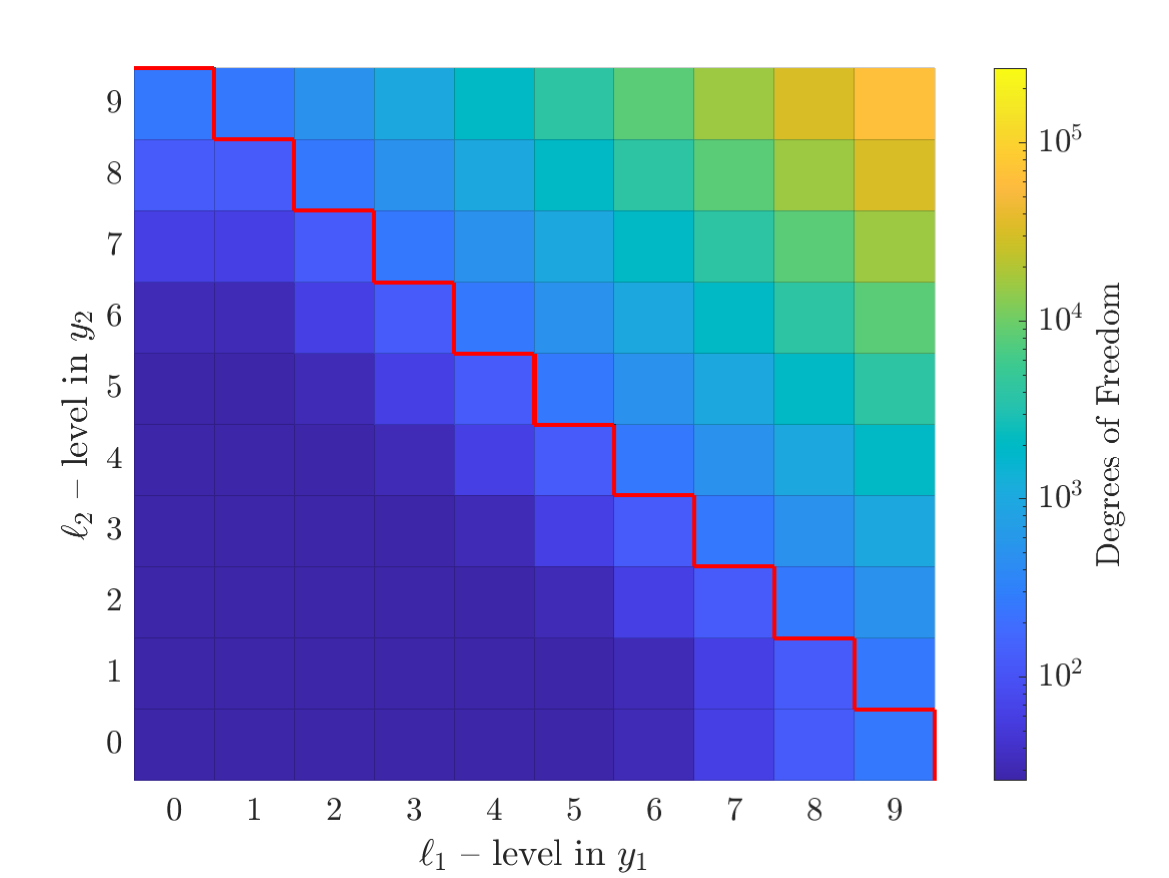}
        \caption{A heat map for the degrees of freedom of $\bm{W}_{\bm{\ell}}$ for a 2 dimensional problem. The whole rectangle corresponds to all degrees of freedom for the full-grid $\bm{V}_9$ while the sparse-grid $\hat{\bm{V}}_9$ only contains the spaces on the lower-left portion divided by the red line.}
    \label{fig:sg_vs_fg_dof_heatmap}    
    \end{subfigure}
    \begin{subfigure}[b]{0.45\textwidth}
        \includegraphics[height=2.15in]{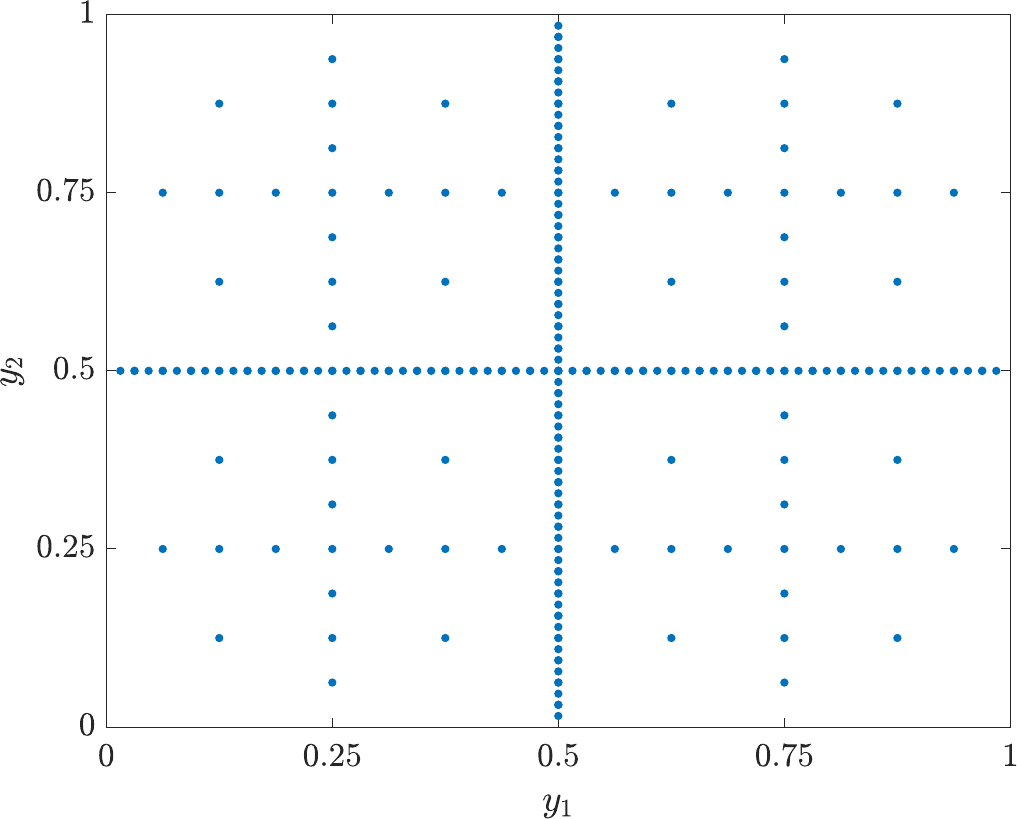}
        \caption{Plot showing the coverage of sparse-grids in two dimensions.  Each point represents the barycenter of the support of a wavelet that is in the level 7 sparse-grid.\\ \\}
    \label{fig:sg_vs_fg_dots}
    \end{subfigure}
    \caption{Sparse-grid illustrations.}
    \label{fig:sparse_grid_illustration} 
\end{figure}

\Cref{fig:sparse_grid_illustration} illustrates which basis functions are kept in the sparse-grid and the reduction in degrees of freedom that sparse-grids provide for the case with $d=2$.  
\Cref{fig:sg_vs_fg_dof_heatmap} shows that the dimension of the spaces $\bm{W}_{\bm{\ell}}$ being thrown away in the sparse-grid truncation are significantly larger on average than the dimension of the spaces that are kept.  
As a result, there is a reduction in degrees of freedom from the full-grid space $\bm{V}_9$ of size $2^{18}\approx 2.62\times 10^5$ to the sparse grid space $\hat{\bm{V}}_9$ of size 2816.\footnote{Here we use $k=0$ to calculate $\dim({{\bm{V}}_9})$ and $\dim({\hat{\bm{V}}_9})$.}
\Cref{fig:sg_vs_fg_dots} shows that the basis functions kept in the sparse-grid allow accurate approximations of derivatives in coordinate directions while throwing away mixed-derivative data which is assumed to be smaller than the components kept by the sparse-grid. 
It has been shown that $\hat{\bm{V}}_N$ shares similar approximation properties to $\bm{V}_N$ in $L^2$, which is $\mathcal{O}(h^{k+1})$, up to a poly-logarithmic factor of $|\log_2 h|^{d-1}$ (see \cite[Theorem 2.4]{wang2016sparse}). 
This result holds for functions with bounded mixed derivatives of sufficient order.

\subsection{Adaptive Sparse-grids}\label{subsec:adaptive_sparse_grids}

The adaptive sparse-grid method uses an adaptive algorithm based on the hierarchical framework of the sparse-grid method \cite{guo2017adaptive}. 
The first step is to further decompose the orthogonal complements $\bm{W_\ell}$ by their level $\bm{\ell}$ and position $\bm{j}$ within the level.
This position $\bm{j}$ in the level is based on the multiwavelet basis.
Given the basis in \eqref{eqn:wavelet_basis_nD_def}, we define the space $\bm{W_{\ell,j}}\subset \bm{W_\ell}$, called a \textit{hierarchical element}, by
\begin{equation}\label{eqn:adaptive_basis_def}
    \bm{W_{\ell,j}} = \spann_{\substack{1\leq i_m\leq k+1 \\ 1\leq m\leq d}} \{g_{\bm{\ell,j}}^{\bm{i}}\}.
\end{equation}
This space has dimension $\textrm{dim}(\bm{W_{\ell,j}})=(k+1)^d$ and
\begin{equation}
    \bm{W_\ell} = \bigoplus_{\bm{j}\in\mathcal{B}_{\bm{\ell}}} \bm{W_{\ell,j}}
\end{equation}
where
\begin{equation}
    \mathcal{B}_{\bm{\ell}} := \{\bm{j}=(j_1,\ldots,j_d):j_m=0,\ldots,\max\{0,2^{\ell_m-1}-1\},\forall m=1\ddd d\}.
\end{equation}
The spaces $\bm{W}_{\bm{\ell},\bm{j}}$ are deemed hierarchical because they carry a natural parent-child relationship of which the details will be postponed (see \Cref{defn:parents_and_children}).
The full- and sparse-grid spaces, \eqref{eqn:full_grid_def} and \eqref{eqn:sparse_grid_def} respectively, can be written as
\begin{equation}\label{eqn:full_and_sparse_grid_adapt_hier}
    \bm{V}_N = \bigoplus_{\substack{|\bm{\ell}|_\infty\leq N \\ 
 \bm{j}\in\mathcal{B}_{\bm{\ell}}}} \bm{W_{\ell,j}} \qquand \hat{\bm{V}}_N = \bigoplus_{\substack{|\bm{\ell}|_1\leq N \\ 
 \bm{j}\in\mathcal{B}_{\bm{\ell}}}} \bm{W_{\ell,j}}.
\end{equation}
For the adaptive sparse-grid algorithm, it is helpful to view the full- and sparse-grid spaces as direct sums of the hierarchical elements $\bm{W_{\ell,j}}$.  We can now define an adaptive sparse-grid which is an arbitrary collection of hierarchical elements.

\begin{defn}\label{defn:adapt_grid}
    Given a \textit{max level} $N_{\text{max}}\in\mathbb{N}_0$ and a \textit{level index set} $\{(\bm{\ell}^\iota, \bm{j}^\iota)\}_{\iota=1}^M$ such that for all $\iota=1\ddd M$, $|\bm{\ell}^\iota|_\infty \leq N_{\text{max}}$ and $\bm{j}^\iota\in\mathcal{B}_{\bm{\ell}}$, 
    the adaptive sparse-grid $\bm{V}\subseteq \bm{V}_{N_{\text{max}}}$ is defined as 
    \begin{equation}\label{eqn:adaptive_grid_def}
        \bm{V} = \bigoplus_{\iota} \bm{W}_{\bm{\ell}^\iota,\bm{j}^\iota}.
    \end{equation}
    Here $M$ is said to be the number of active elements of the adaptive sparse-grid $\bm{V}$.
\end{defn}
We will often drop the $\iota$ superscript in \eqref{eqn:adaptive_grid_def} and refer to the level index set as $\{(\bm{\ell},\bm{j})\}$.  
From \eqref{eqn:full_and_sparse_grid_adapt_hier}, the standard sparse-grid is a specific adaptive sparse-grid where we include all hierarchical elements $\bm{W}_{\bm{\ell},\bm{j}}$ such that $|\bm{\ell}|_1\leq N_{\text{max}}$ and $\bm{j}\in\mathcal{B}_{\bm{\ell}}$.   

\subsubsection{Adaptive Approximation of Initial Data}\label{subsec:adapt_initial_data}

Let $\mathcal{P}_{\bm{V}}$ be the $L^2$ projection from $L^2(\W)$ onto $\bm{V}$. The main idea of the adaptive sparse-grid is to choose a grid $\bm{V}\subseteq\bm{V}_{N_{\text{max}}}$, depending on the distribution $w$, such that 
\begin{enumerate}
    \item The relative projection error $\|w-\mathcal{P}_{\bm{V}}w\|_{L^2(\W)}/\|w\|_{L^2(\W)}$ is small;
    \item $\textrm{dim}(\bm{V})$ is approximately minimal. 
\end{enumerate}
We will first demonstrate this process for an initial condition, and then extend the result to functions formulated via a dynamical system.

For a fixed max level $N_{\text{max}}$, choosing $\bm{V}=\bm{V}_{N_{\text{max}}}$ would minimize the $L^2$ projection error over all possible adaptive sparse-grid spaces, but with significant costs in terms of the number of degrees of freedom.  Thus we assume $w\in\bm{V}_{N_{\text{max}}}$ is our target; then the coefficient expansion with respect to the multiwavelet basis of \eqref{eqn:wavelet_basis_nD_def} is given by
\begin{equation}\label{eqn:coeff_expansion_full_grid}
    w(\bm{y}) = \sum_{\Big\{\substack{(\bm{\ell},\bm{j}): \\ |\bm{\ell}|_\infty\leq N_{\text{max}},
 \bm{j}\in\mathcal{B}_{\bm{\ell}}}\Big\}} \sum_{\substack{1\leq i_m\leq k+1 \\ 1\leq m\leq d}} w_{\bm{\ell},\bm{j}}^{\bm{i}} g_{\bm{\ell},\bm{j}}^{\bm{i}}(\bm{y})~~\text{where}~~ w_{\bm{\ell},\bm{j}}^{\bm{i}} = \int_\W w(\bm{y})g_{\bm{\ell},\bm{j}}^{\bm{i}}(\bm{y})\dx{\bm{y}}.
\end{equation}
For simplification, we define $w_{\bm{\ell},\bm{j}}$ to be the multilinear rank-$d$ tensor with $k+1$ entries in each dimension, defined by
\begin{equation}
    [w_{\bm{\ell},\bm{j}}]_{\bm{i}} = w_{\bm{\ell},\bm{j}}^{\bm{i}}.
\end{equation}
When taking the norm of $w_{\bm{\ell},\bm{j}}$, we first flatten the tensor into a vector in $\mathbb{R}^{(k+1)^d}$ and apply the appropriate vector norm in $\ell^p$ where $1\leq p\leq \infty$.  
Then the $L^2$-norm of $w$ can be written as
\begin{equation}
    \|w\|_{L^2(\W)}^2 = \sum_{\Big\{\substack{(\bm{\ell},\bm{j}): \\ |\bm{\ell}|_\infty\leq N_{\text{max}},
 \bm{j}\in\mathcal{B}_{\bm{\ell}}}\Big\}} \|w_{\bm{\ell},\bm{j}}\|_2^2
\end{equation}
Additionally, for any adaptive sparse-grid space $\bm{V}$ with level index set $\{(\bm{\ell},\bm{j})\}$ we have
\begin{equation}\label{eqn:coeff_expansion_adapt_grid}
    \mathcal{P}_{\bm{V}}w = \sum_{(\bm{\ell},\bm{j})} \sum_{\substack{1\leq i_m\leq k+1 \\ 1\leq m\leq d}} w_{\bm{\ell},\bm{j}}^{\bm{i}} g_{\bm{\ell},\bm{j}}^{\bm{i}}~~\text{and}~~\|\mathcal{P}_{\bm{V}}w\|_{L^2(\W)}^2 = \sum_{(\bm{\ell},\bm{j})} \|w_{\bm{\ell},\bm{j}}\|_2^2.
\end{equation}
From \eqref{eqn:coeff_expansion_full_grid} and \eqref{eqn:coeff_expansion_adapt_grid}, it is clear that the relative projection error satisfies
\begin{equation}\label{eqn:projection_err_adapt}
    \frac{\|w-\mathcal{P}_{\bm{V}}w\|_{L^2(\W)}^2}{\|w\|_{L^2(\W)}^2} = \sum_{(\bm{\ell},\bm{j}):\bm{W}_{\bm{\ell},\bm{j}}\not\subseteq\bm{V}} \frac{\|w_{\bm{\ell},\bm{j}}\|_2^2}{\|w\|_{L^2(\W)}^2}.
\end{equation}
Therefore, given $\tau>0$, called the \textit{threshold}, we want to keep all hierarchical elements $\bm{W}_{\bm{\ell},\bm{j}}$ such that 
\begin{equation}\label{eqn:refine_requirement}
    \|w_{\bm{\ell},\bm{j}}\|_2 \geq \tau\|w\|_{L^2(\W)}
\end{equation}
lest they contribute to the error in \eqref{eqn:projection_err_adapt}.

\paragraph{Refinement}
We will now describe how hierarchical elements are added to the adaptive sparse-grid -- which we call refinement. The primary challenge in building a grid that contains all elements satisfying \eqref{eqn:refine_requirement} is to avoid checking all hierarchical elements in the full-grid -- an operation that naively would require $\mathcal{O}(2^{Nd})$ operations.   

The refinement process is iterative, where an initial grid is chosen and then added upon.  For adapting an initial condition, we choose our initial grid as the sparse-grid $\bm{V}=\hat{\bm{V}}_{N_{\text{max}}}$.  Given a current grid $\bm{V}$, the coefficients $w_{\bm{\ell},\bm{j}}$ are computed for every hierarchical element in the grid.  In order to determine what elements to add to the grid, we appeal to the hierarchical representation of the full-grid space which embeds the following parent-child relation.

\begin{defn}\label{defn:parents_and_children}
    Let $\bm{W}_{\bm{\ell},\bm{j}}$ for $(\bm{\ell},\bm{j})=\big( (\ell_1,\ldots,\ell_d),(j_1,\ldots,j_d) \big)$ be a hierarchical element with max level $N_{\text{max}}$.  The children of $\bm{W}_{\bm{\ell},\bm{j}}$, with up to two per dimension, are defined for each dimension $m=1,\ldots,d$ by the following:
    \begin{itemize}
        \item If $\ell_m=0$, then $\bm{W}_{\bm{\ell}',\bm{j}'}$, where \begin{equation}
            (\bm{\ell}',\bm{j}')=\big( (\ell_1,\ldots,\ell_{m-1},1,\ell_{m+1},\ldots,\ell_d),(j_1,\ldots,j_{m-1},0,j_{m+1},\ldots,j_d) \big),
        \end{equation} is a child of $\bm{W}_{\bm{\ell},\bm{j}}$.
        \item If $0<\ell_m<N_{\text{max}}$, then $\bm{W}_{\bm{\ell}',\bm{j}'}$, where
        \begin{subequations}
          \begin{align}
            (\bm{\ell}',\bm{j}')&=\big( (\ell_1,\ldots,\ell_{m-1},\ell_m+1,\ell_{m+1},\ldots,\ell_d),(j_1,\ldots,j_{m-1},2j_m,j_{m+1},\ldots,j_d) \big)  \quad\text{and}\\
            (\bm{\ell}',\bm{j}') &=\big( (\ell_1,\ldots,\ell_{m-1},\ell_m+1,\ell_{m+1},\ldots,\ell_d),(j_1,\ldots,j_{m-1},2j_m+1,j_{m+1},\ldots,j_d) \big),
        \end{align}   
        \end{subequations}
        are children of $\bm{W}_{\bm{\ell},\bm{j}}$.
        \item If $\ell_m=N_{\text{max}}$, then there are no children of $\bm{W}_{\bm{\ell},\bm{j}}$ in dimension $m$.
    \end{itemize}
    The parents of an element $\bm{W}_{\bm{\ell},\bm{j}}$ are all elements $\bm{W}_{\bm{\ell}',\bm{j}'}$ such that $\bm{W}_{\bm{\ell},\bm{j}}$ is a child of $\bm{W}_{\bm{\ell}',\bm{j}'}$.
\end{defn}
It is clear from \Cref{defn:parents_and_children} that each hierarchical element can have up to $2d$ children and up to $d$ parents.
To tie \Cref{defn:parents_and_children} to the wavelet representation, for a fixed dimension $m \in \{1,\dots,d\}$, the children of a wavelet given in \eqref{eqn:wavelet_basis_1D_def} are the up to two wavelets of one greater level whose support is contained in the parent.
Furthermore, based on the coefficient decay estimate \eqref{eqn:hier_approx_result_nd}, if $w$ is sufficiently smooth, then it is reasonable to assume that if $\bm{W}_{\bm{\ell}',\bm{j}'}$ is a child of $\bm{W}_{\bm{\ell},\bm{j}}$, then $\|w_{\bm{\ell}',\bm{j}'}\|\leq\|w_{\bm{\ell},\bm{j}}\|$.
Therefore, if the size of a hierarchical element in the grid is small, we assume the size of the children are also small, and we do not need to search further along this path.
This assumption leads to a stopping mechanism for the refinement strategy:  Given a grid $\bm{V}$ with level index set $\{(\bm{\ell},\bm{j})\}$, if 
\begin{equation}\label{eqn:refine_stopping_criterion}
    \|w_{\bm{\ell},\bm{j}}\|_2 \geq \tau\Big(\sum_{({\bm{\ell}',\bm{j}}')} \|w_{\bm{\ell}',\bm{j}'}\|_2^2\Big)^{\frac{1}{2}} = \tau\|\mathcal{P}_{\bm{V}}w\|_{L^2(\W)},
\end{equation} then we add all children of $\bm{W}_{\bm{\ell},\bm{j}}$ to the grid.  We repeat this process iteratively until no new children are added.

\paragraph{Coarsening} The process of removing active elements from the current grid, i.e., coarsening, is achieved by simple thresholding of the coefficients.
Let $0<\mu<1$ be the \textit{coarsening factor}.
For a given grid $\bm{V}$, if $\bm{W}_{\bm{\ell},\bm{j}}$ is a hierarchical element such that 
\begin{equation}\label{eqn:coarsen_stopping_criterion}
    \|w_{\bm{\ell},\bm{j}}\|_2 \leq \mu\tau \Big(\sum_{({\bm{\ell}',\bm{j}}')}
    \|w_{\bm{\ell}',\bm{j}'}\|_2^2\Big)^{\frac{1}{2}} = \mu\tau\|\mathcal{P}_{\bm{V}}w\|_{L^2(\W)},
\end{equation} then it is removed from the grid.
We acknowledge this coarsening strategy does not preserve structural properties like parent completeness, i.e., requiring that all parents of an active element are active (see \cite[Section 3]{stoyanov2018adaptive}), but still yields stable and accurate approximations as evidenced in \Cref{sec:numerical}.

While the criteria for refinement \eqref{eqn:refine_stopping_criterion} and coarsening \eqref{eqn:coarsen_stopping_criterion} are based on $\ell^2$-type norms, other discrete norms can be used.
For instance, the $\ell^\infty$ norm can also be used:
\begin{subequations}\label{eqn:inf_norm_criteria}
\begin{alignat}{2}
    \|w_{\bm{\ell},\bm{j}}\|_\infty &\geq \tau \max_{({\bm{\ell}',\bm{j}}')} \|w_{\bm{\ell}',\bm{j}'}\|_\infty &&
    \qquad \text{(for refinement),}
    \label{eqn:inf_norm_criteria:refine}\\
    \|w_{\bm{\ell},\bm{j}}\|_\infty &\leq \mu\tau \max_{({\bm{\ell}',\bm{j}}')} \|w_{\bm{\ell}',\bm{j}'}\|_\infty &&
    \qquad \text{(for coarsening).} \label{eqn:inf_norm_criteria:coarsen}
\end{alignat}
\end{subequations}
\asgard~uses \eqref{eqn:inf_norm_criteria} as its refinement/coarsening criteria.

\subsubsection{Adaptive sparse-grids of a dynamical system} 

Unlike adapting initial conditions, where the coefficients are drawn from analytic or quadrature data, the adaptive strategy can also be utilized to create temporally varying grids that dynamically capture features of the solution in time.  
To extend our adaptive strategy to dynamical systems, consider the abstract problem
\begin{equation}\label{eqn:abstract_problem}
(\partial_t w,g) = \mathcal{A}(w,g)\quad\forall g\in \bm{V}_{N_{\text{max}}}
\end{equation}
where $\mathcal{A}:\bm{V}_{N_{\text{max}}}\times\bm{V}_{N_{\text{max}}}\to\mathbb{R}$ (c.f.~\eqref{eqn:discrete_form}). Here $\cA$ is one of the discretizations in \eqref{eqn:VP_discrete_def} or \eqref{eqn:LB_discrete_def}.  For a given adaptive sparse-grid $\bm{V}\subseteq\bm{V}_N$ define the operator $\mathcal{R}:\bm{V}\to\bm{V}$ by
\begin{equation}\label{eqn:ritg_galerkin_R_def}
    (\mathcal{R}_{\bm{V}}w,g) = \mathcal{A}(w,g) \quad\forall g\in \bm{V}.
\end{equation}
Then \eqref{eqn:abstract_problem} can be succinctly written as $\partial_tw=\mathcal{R}_{\bm{V}_{N_{\text{max}}}} w$.

Consider a solution $w^{\mfn}$ at timestep $t^\mfn$ defined on an adaptive sparse-grid $\bm{V}^\mfn$.  
To refine, we first set $\bm{V}=\bm{V}^\mfn$ and advance the abstract problem $\partial_tw = \mathcal{R}_{\bm{V}}w$ from $t^\mfn$ to $t^{\mfn+1}$ via a IMEX Runge--Kutta method \eqref{eqn:IMEX_RK} to produce $w^{\mfn+1}\in \bm{V}$.  
We then check for elements $\bm{W}_{\bm{\ell},\bm{j}}$ of $\bm{V}$ that satisfy the same refinement requirement as the initial condition case, namely,
\eqref{eqn:refine_stopping_criterion} for a $\ell^2$-norm refinement or  \eqref{eqn:inf_norm_criteria:refine} for a $\ell^\infty$-norm refinement.
If there are elements satisfying the refinement criterion, then their children are added to $\bm{V}$.  
We then go back to time $t^\mfn$ and advance $\partial_tw = \mathcal{R}_{\bm{V}}w$ from $t^\mfn$ to $t^{\mfn+1}$ with the updated space $\bm{V}$.  
Since $\bm{V}^\mfn\subseteq \bm{V}$, the coefficients of the state $w^\mfn$ can be extended into $\bm{V}$ by setting $w_{\bm{\ell},\bm{j}}^\mfn=0$ if $\bm{W}_{\bm{\ell},\bm{j}}\subseteq\bm{V}$ but not if $\bm{W}_{\bm{\ell},\bm{j}}\not\subseteq\bm{V}^\mfn$. 
This process is repeated until no new children are added into the grid $\bm{V}$ -- in which case we set $\bm{V}^{\mfn+1}=\bm{V}$. Typically, only one or two refinements are needed per timestep, but more may be needed for the first few timesteps due to initial layers.
Coarsening after refinement is done in a manner analogous to the initial condition case.  The procedure for refining and coarsening are summed up in \Cref{alg:adapt_refine} and \Cref{alg:adapt_coarsen} respectively.

\begin{algorithm}
\caption{Adaptive refinement using $\ell^\infty$-norm}\label{alg:adapt_refine}

\SetKwInOut{Input}{Input}
\SetKwInOut{Output}{Output}
\SetKwRepeat{Do}{do}{while}
\SetKwComment{Comment}{/* }{ */}

\Input{Adaptive sparse-grid $\bm{V}$, threshold $\tau>0$}
\Input{Distribution $w$ or dynamical system $(w,\bm{U})\to\partial_tw=\mathcal{R}_{\bm{U}}w$ defined in \eqref{eqn:ritg_galerkin_R_def} with coefficients $w_{\bm{\ell},\bm{j}}^\mfn$ computed for all $\bm{W}_{\bm{\ell},\bm{j}}\subset \bm{V}$}
\Output{Adaptive sparse-grid $\bm{V}^*$}
\Output{Coefficients $w_{\bm{\ell},\bm{j}}$ for all $\bm{W}_{\bm{\ell},\bm{j}}\subset \bm{V}^*$}
$\bm{V}^*:=\bm{V}$\;
\Do{$\bm{N}\neq\{0\}$}{
    $\bm{N}:=\{0\}$\;
    Compute $w_{\bm{\ell},\bm{j}}$ for all $\bm{W}_{\bm{\ell},\bm{j}}\subset \bm{V}^*$ via $w$ or dynamical system $\partial_t f = \mathcal{R}_{\bm{V}^*}f$\;
    \For{$\bm{W}_{\bm{\ell},\bm{j}}\subset \bm{V}^*$}{
        \If(\Comment*[f]{Check if element needs refining}){$\|w_{\bm{\ell},\bm{j}}\|_\infty \geq \tau \max_{({\bm{\ell}',\bm{j}}')} \|w_{\bm{\ell}',\bm{j}'}\|_\infty$}{
            Compute all children $\bm{W}_{\bm{\ell}',\bm{j}'}$ of $\bm{W}_{\bm{\ell},\bm{j}}$ using \Cref{defn:parents_and_children}\;
            \For{children $\bm{W}_{\bm{\ell}',\bm{j}'}$ of $\bm{W}_{\bm{\ell},\bm{j}}$}{
                \If{$\bm{W}_{\bm{\ell}',\bm{j}'}\not\subset\bm{V}^*$ and $\bm{W}_{\bm{\ell}',\bm{j}'}\not\subset\bm{N}$}{
                    $\bm{N}:=\bm{N}\oplus\bm{W}_{\bm{\ell}',\bm{j}'}$\Comment*[r]{Add element}
                    $w_{\bm{\ell}',\bm{j}'}^\mfn=0$\Comment*[r]{Zero out new element at $t^\mfn$}
                }
            }
        }
    }
    $\bm{V}^*:=\bm{V}^*\oplus\bm{N}$\;
}(\Comment*[f]{Repeat until no children are added.})
\end{algorithm}

\begin{algorithm}
\caption{Adaptive coarsening using $\ell^\infty$-norm}\label{alg:adapt_coarsen}

\SetKwInOut{Input}{Input}
\SetKwInOut{Output}{Output}
\SetKwRepeat{Do}{do}{while}
\SetKwComment{Comment}{/* }{ */}

\Input{Adaptive sparse-grid $\bm{V}$, threshold $\tau>0$, coarsening factor $0<\mu<1$}
\Input{Coefficients $w_{\bm{\ell},\bm{j}}$ for all $\bm{W}_{\bm{\ell},\bm{j}}\subset \bm{V}$}
\Output{Adaptive sparse-grid $\bm{V}^*$}
$\bm{V}^*:=\{0\}$\;
\For{$\bm{W}_{\bm{\ell},\bm{j}}\subset \bm{V}$}{
    \If(\Comment*[f]{Check if element needs to be removed}){$\|w_{\bm{\ell},\bm{j}}\|_\infty > \mu\tau \max_{({\bm{\ell}',\bm{j}}')} \|w_{\bm{\ell}',\bm{j}'}\|_\infty$}{
        $\bm{V}^*:=\bm{V}^*\oplus\bm{W}_{\bm{\ell},\bm{j}}$\;
    }
}
\end{algorithm}

\begin{figure}[!ht]
    \centering
    \begin{subfigure}[b]{0.45\textwidth}
        \includegraphics[width=\textwidth]{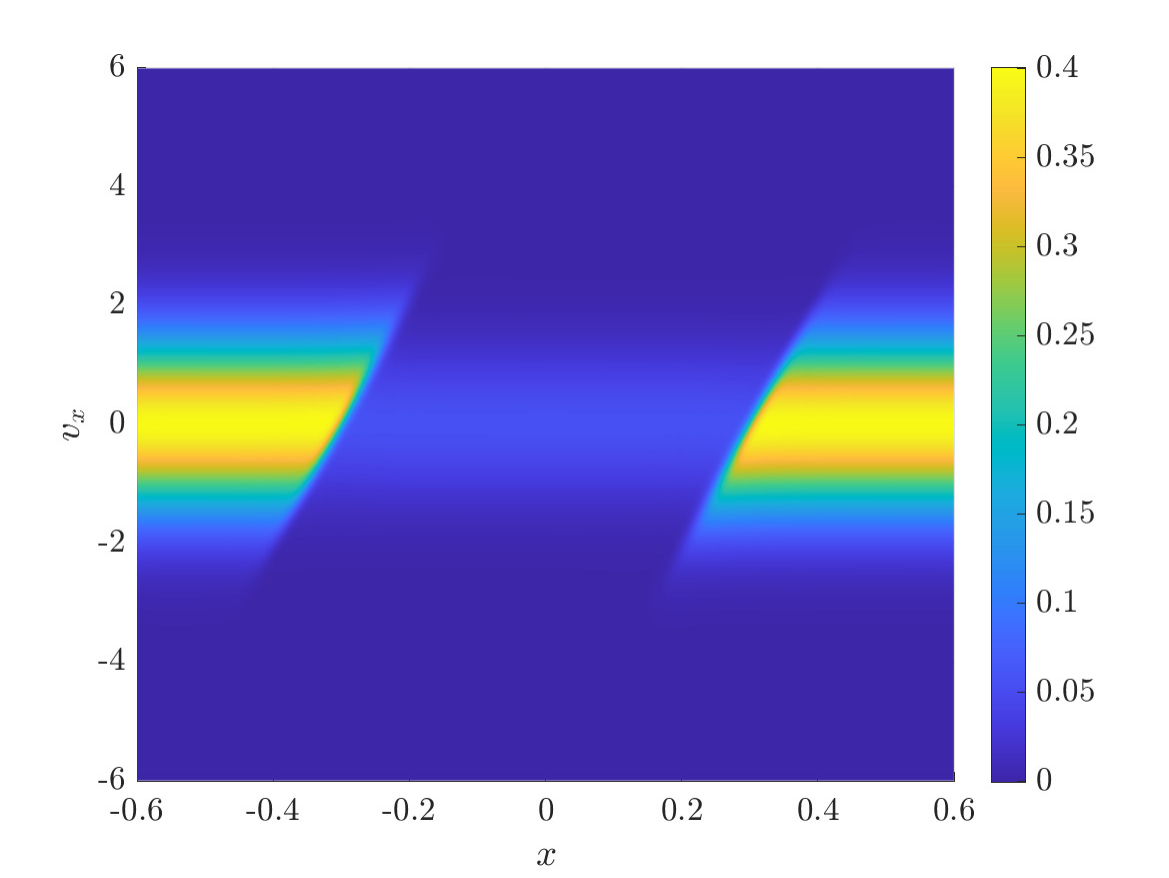}
        \caption{Phase space plot of $\langle f\rangle_{v_y,v_z}$ in $(x,v_x)$ where $f$ is the adaptive sparse-grid distribution.\newline }
    \label{fig:adapt_riemann:dist}    
    \end{subfigure}
    \begin{subfigure}[b]{0.45\textwidth}
        \includegraphics[width=0.85\textwidth]{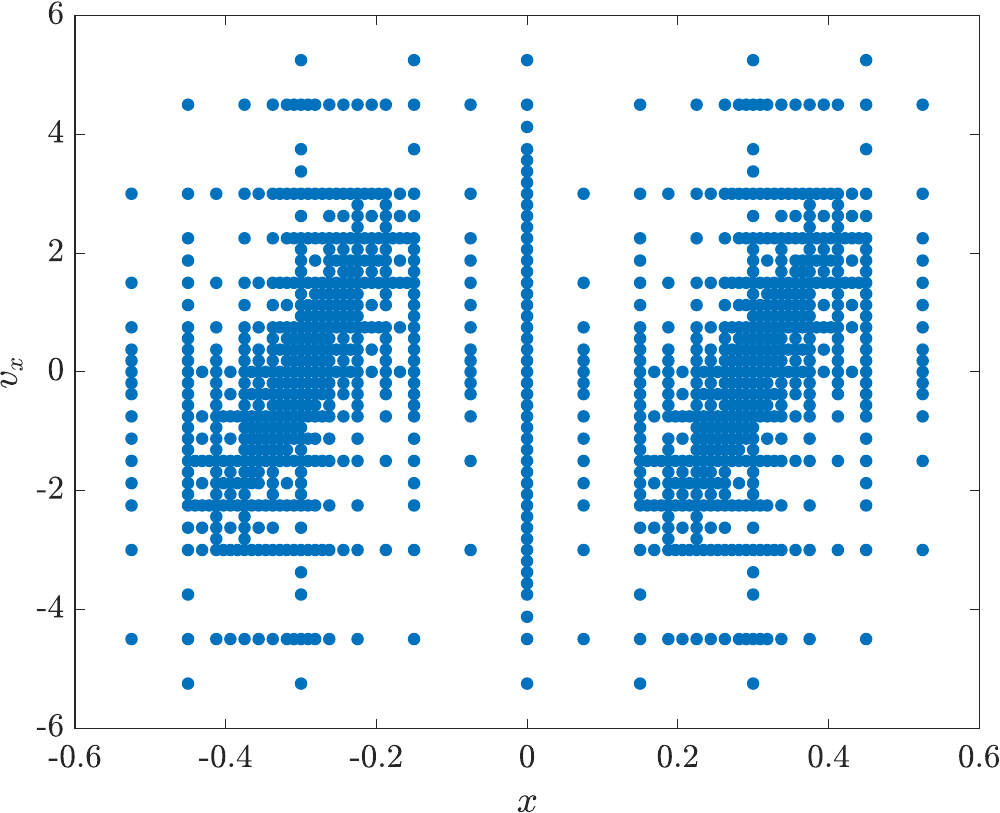}
        \caption{Plot showing which hierarchical elements are active in \Cref{fig:adapt_riemann:dist}.  The points represent the $(x,v_x)$ coordinates of the barycenter of the support of each active wavelet.}
    \label{fig:adapt_riemann:dots}    
    \end{subfigure}
    \caption{Riemann problem -- \Cref{subsec:riemann} -- $\nu=1$: Adaptive Sparse-grid Method at $t=0.04918$.  The threshold is $\tau=10^{-4}$ and the adaptive sparse-grid cannot refine past $\bm{\ell}=(7,6,6,6)$. }
    \label{fig:adapt_riemann} 
\end{figure}

As visual illustration of the adaptive sparse-grid method is shown in \Cref{fig:adapt_riemann}, where it is applied to $1x3v$ Riemann problem in \Cref{subsec:riemann}.  
As seen in \Cref{fig:adapt_riemann:dots}, the adaptive algorithm focuses on refinement around the discontinuity in the distribution, plotted in \Cref{fig:adapt_riemann:dist}, while coarsening occurs near the velocity boundaries.  

\section{Numerical Experiments}
\label{sec:numerical}

In this section, we present results from various test problems relevant to plasma physics.   Our goals are to demonstrate the performance of the adaptive sparse-grid and mixed-grid DG methods with IMEX time stepping implemented in \asgard, and investigate the computational benefit of the adaptive sparse-grid and mixed-grid methods over the full-grid methods (see \Cref{subsec:choice_of_grids} for definitions).  
In increasing degree of complexity, we consider: (i) relaxation to a Maxwellian velocity distribution (\Cref{subsec:relax}); (ii) a Riemann problem for two different values of the collision frequency $\nu$ (\Cref{subsec:riemann}); (iii) and the collisional Landau damping problem (\Cref{subsec:landau}), also for two different values of the collision frequency.   All the results presented in this section were obtained with quadratic polynomials, i.e., $k=2$.   This choice of $k$ natural considering that the velocity moments with respect to $1$, $\bv$, and $|\bv|^2$ are the important fluid variables.

\subsection{Choice of Grids}\label{subsec:choice_of_grids}

In the simulations presented below we choose to compare results obtained with three types of grids: full-grid, mixed-grid, and adaptive sparse-grid.  
We provide the specifics of each grid in this section.

Our first choice is the standard full-grid $\bm{V}_{\bm{\ell}}$, where $\bm{\ell}=(\ell_x,\ell_v,\ell_v,\ell_v)$.
We use the Chu reduction method of \Cref{subsec:reduction} in \Cref{subsec:riemann,subsec:landau} with $\bm{\ell}=(\ell_x,\ell_v)$ to build reference solutions and numerical approximations with the full-grid.  
This is because the full-grid space is too large in comparison to the other two grids and can easily fill the memory of a single-node machine.  When using the Chu reduction, the discretization is performed using a local Legendre polynomial basis instead of the multiwavelets. 
When determining the degrees of freedom or number of active elements for a full-grid run, we will always assume that the underlying run is 4D, even if the Chu reduction method is used.  

We have found that standard 4D sparse-grids such as $\hat{\bm{V}}_{\bm{\ell}}$ are unstable for the VPLB model in \eqref{eqn:kinetic_slab}.  
This is due to both the lack of resolution in $x$ and the lack of regularity of the distribution function in physical space.  
Specifically, the temperature $\theta_f$ becomes negative which causes the solution to blow up.
As we expect savings to come from the smoothness in velocity space, induced by the LB collision operator, we propose a mixed-grid approach for our second choice.  
The mixed-grid of level $\bm{\ell}=(\ell_x,\ell_v,\ell_v,\ell_v)$ is defined by 
\begin{equation}\label{eqn:mixed_grid_def}
    \tilde{\bm{V}}_{\bm{\ell}}=\bigoplus_{\substack{\bm{\ell}': \ell_1'\leq \ell_x, \\|(\ell_2',\ell_3',\ell_4')|_1\leq \ell_v}}\bm{W}_{\bm{\ell}'}=V_{\ell_x}\otimes \hat{\bm{V}}_{(\ell_v,\ell_v,\ell_v)}.
\end{equation}
The mixed-grid space is a tensor product of a full-grid in physical space and a sparse-grid in velocity space.  This can be viewed as a sparse-grid in velocity space attached to each degree of freedom in $x$, and thus provides computational savings relative to the full-grid (without the Chu reduction method).  
The dimension of $\tilde{\bm{V}}_{\bm{\ell}}$ is $\mathcal{O}((k+1)^42^{\ell_x+\ell_v}\ell_v^2)$.  
We find this space is sufficient to maintain stability of the DG method for the problems considered here.  

Additionally, since $1$ and $v_y^2+v_z^2$ are admissible DG functions that live on level $(0,0)$ in $(v_y,v_z)$ when $k\ge2$, a 4D full-grid of level $(\ell_x,\ell_v,0,0)$ is sufficient to recover $g_1$ and $g_2$ in a 2D full-grid of level $(\ell_x,\ell_v)$.  
Since a full-grid of level $(\ell_x,\ell_v,0,0)$ is a subgrid of a mixed-grid with level $(\ell_x,\ell_v,\ell_v,\ell_v)$, the reduced moments $g_1$ and $g_2$ created by the mixed-grid solution will be similar to the full-grid.  
However, for $k\le3$, the function $v_y^4+v_z^4$ is not a DG function, and its projection onto the DG space will excite finer level coefficients that are better captured by the full-grid than by the mixed-grid for a certain level.
We therefore evolve $g_3$ in the Chu reduction method in order to better understand differences in accuracy between the mixed-grid and full-grid methods.

Our last grid is the adaptive sparse-grid, $\bm{V}$, that is coarsened and refined as detailed in \Cref{alg:adapt_coarsen,alg:adapt_refine}, using $\ell^\infty$-thresholding. 
The refinement threshold $\tau$ will be problem dependent, but we use the coarsening factor $\mu=0.1$ for all our examples.  Instead of a max level $N_{\text{max}}$ used in \Cref{sec:sparsegrid}, we will not allow the adaptive sparse-grid to refine above a full-grid of specified level $\bm{\ell}=(\ell_x,\ell_v,\ell_v,\ell_v)$.
The number of degrees of freedom, or active elements, presented in the results below will be of the adaptive sparse-grid solution after the refinement step but before coarsening.

It is useful to view each of these grids as a velocity grid attached to each spatial degree of freedom.  
The full-grid attaches a three-dimensional full-velocity grid to every spatial degree of freedom while the mixed-grid attaches a sparse-velocity grid.  
The adaptive sparse-grid attaches a variable velocity grid, with possibly zero elements, to each spatial degree of freedom.

Finally, we will track the number of active elements, see \Cref{defn:adapt_grid}, as opposed to degrees of freedom in order to more clearly present the advantages of the mixed-grid and adaptive sparse-grid methods.

\subsection{Relaxation Problem}\label{subsec:relax}

We first consider the $0x3v$ problem in \eqref{eqn:relax:pde} in order to test the relaxation to equilibrium induced by the LB collision operator.  In this case $f = f(\bv,t)$ and the computational domain is truncated so that $\bv\in (-8,12)^3$.  
The initial condition is given by the sum of three Maxwellians, each sharing $n_{f}=1/3$, $\theta_{f}=1/2$, but differing in the bulk velocities, which are given by $[3,0,0]$, $[0,3,0]$, and $[0,0,3]$, respectively. 
This initial condition induces the following velocity moments: $n_{f}=1$, $\bu_{f}=[1,1,1]^{\intercal}$, and $\theta_{f}=2.5$.  
By the properties of the LB collision operator \Cref{prop:lenardBernsteinOperator}, these moments are expected to remain constant in time and the velocity distribution to relax to the Maxwellian defined by the initial moments.  

\begin{figure}[!ht]
    \centering
    \begin{subfigure}[b]{0.48\textwidth}
        \centering
        \includegraphics[width=\textwidth]{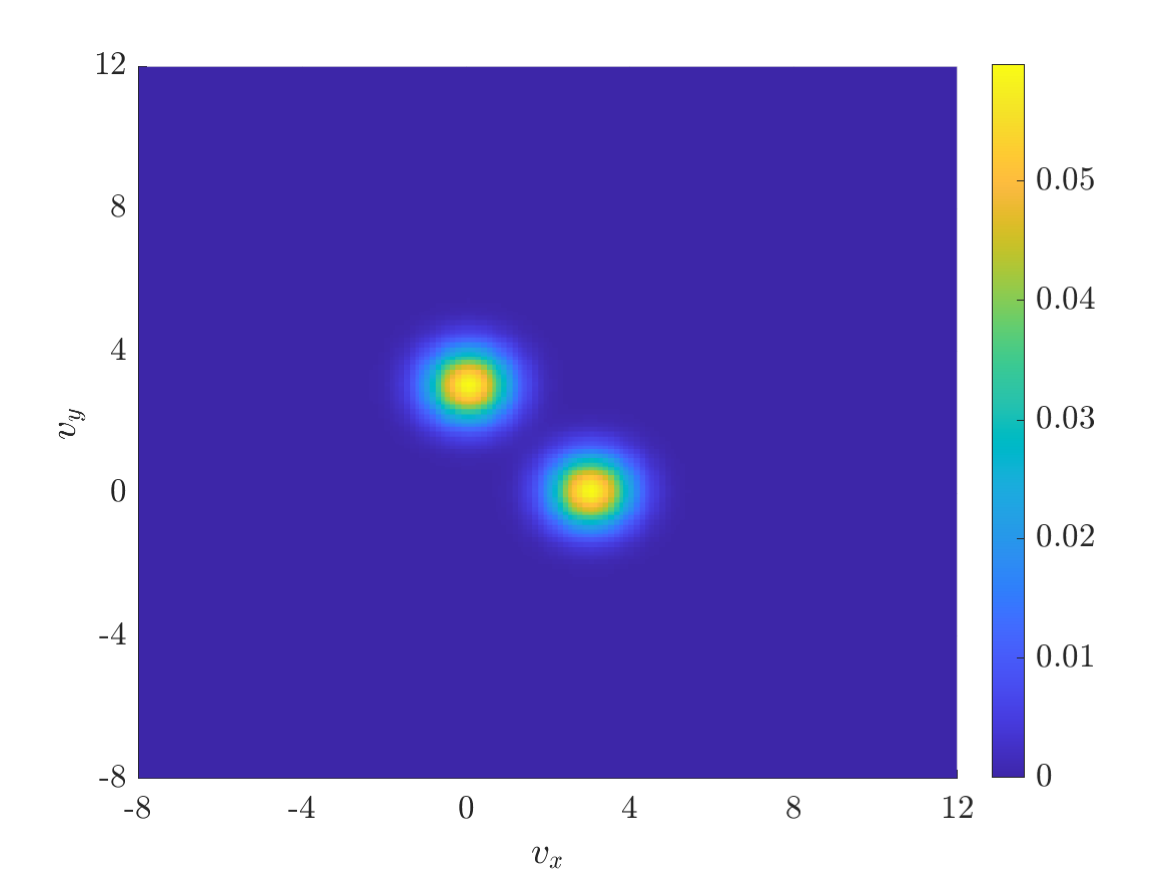}
        \caption{$\nu t=0$}
    \end{subfigure}
    \begin{subfigure}[b]{0.48\textwidth}
        \centering
        \includegraphics[width=\textwidth]{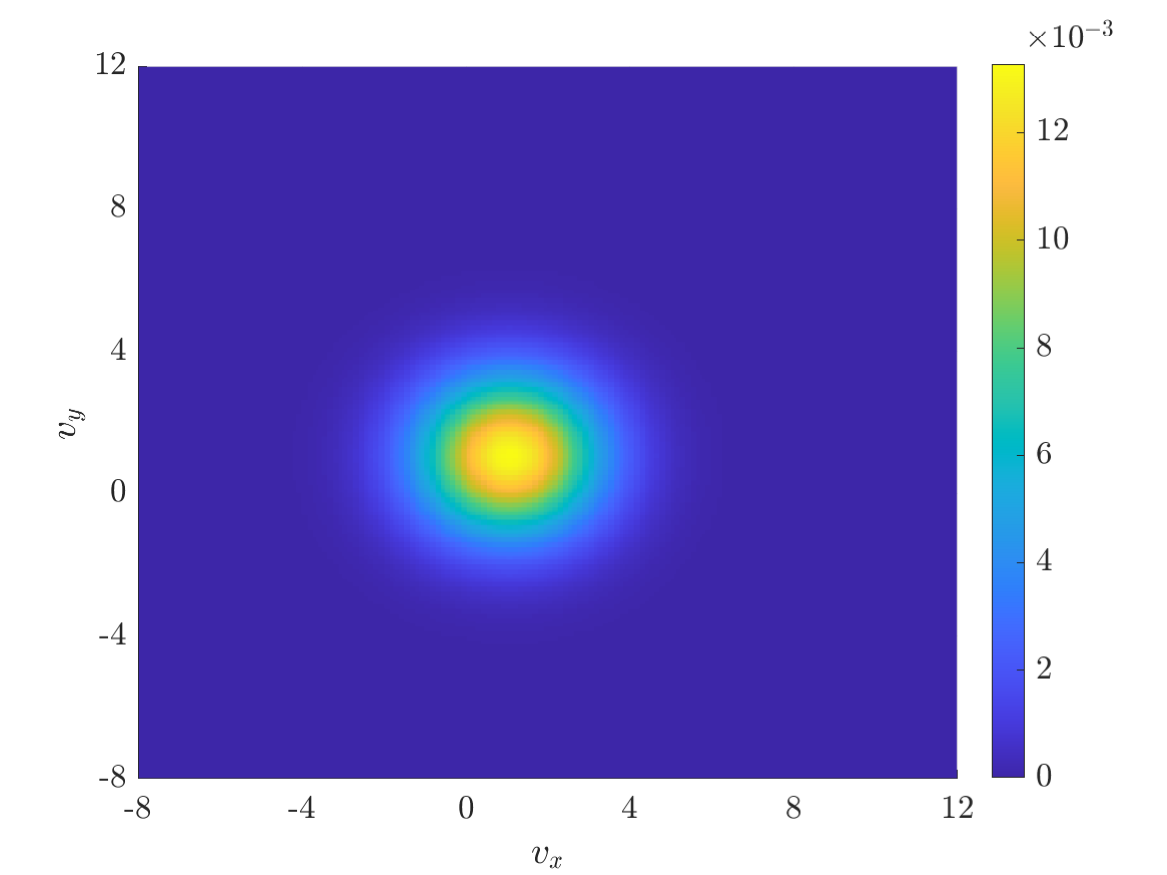}
        \caption{$\nu t=20$}
    \end{subfigure}
    \caption{Relaxation Problem -- \Cref{subsec:relax}: 2D plot of the velocity distribution $f_h(v_x,v_y,v_z=0.019)$ at the start (left) and end (right) of a relaxation simulation.  These results were obtained with a full-grid run with $\bm{\ell}=(5,5,5)$.}
    \label{fig:relaxation:2d_slice}
\end{figure}
  
For this test, we will use a 3D sparse-grid of level $(\ell_v,\ell_v,\ell_v)$ as a substitute for the mixed-grid.  The 4D ($1x3v$) definitions of the full-grid and adaptive sparse-grid naturally carry to the 3D ($0x3v$) case.
We set $\nu=10^{3}$, $\Delta t=5\times10^{-4}$, and use backward Euler time stepping for this problem, with a tolerance of $10^{-8}$ for the GMRES implicit solve.  
\Cref{fig:relaxation:2d_slice} illustrates the initial and final (equilibrium) distributions in the $(v_x,v_y)$-plane for a full-grid model.  

\begin{figure}[!ht]
    \centering
    \begin{subfigure}[b]{0.48\textwidth}
    \includegraphics[width=\textwidth]{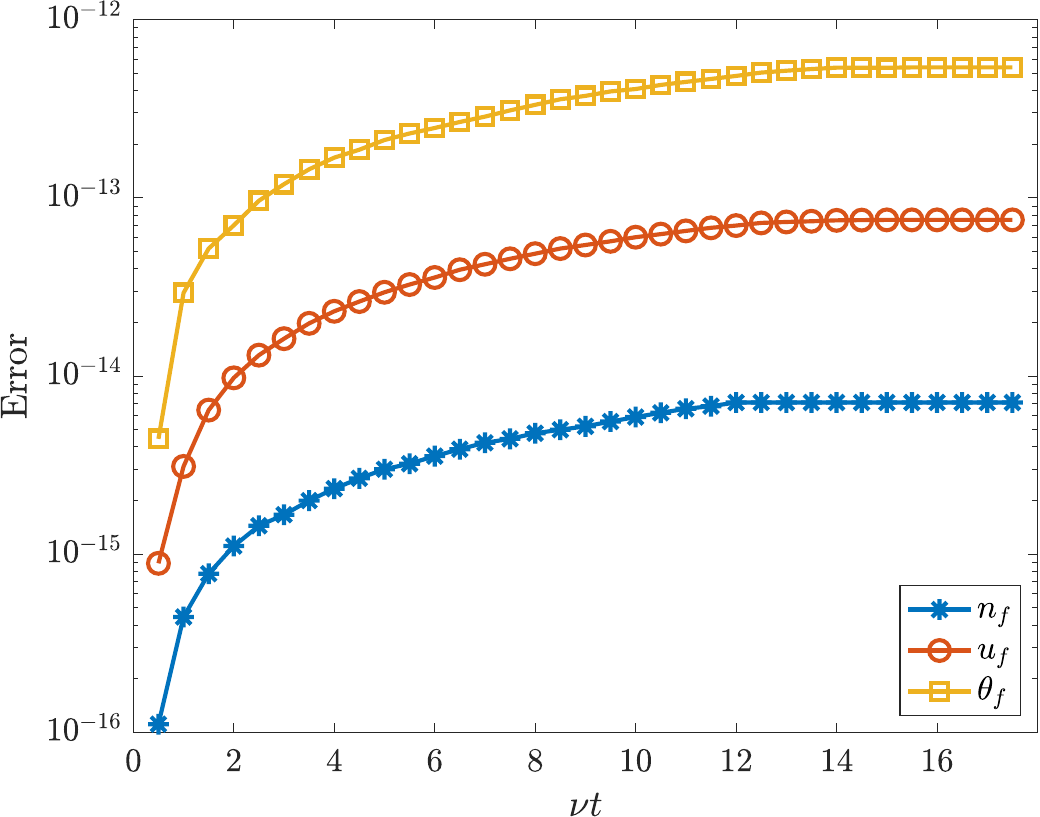}
    \captionsetup{margin=1em}
    \caption{Plots of $|n_f(t)-n_f(0)|$, $|u_f(t)-u_f(0)|$, and $|\theta_f(t)-\theta_f(0)|$ versus time where the fluid variables are approximations with $\bm{\ell} = (4,4,4)$.\newline\newline}
    \label{fig:relaxation:conservation}
    \end{subfigure}\hfill
    \begin{subfigure}[b]{0.48\textwidth}
    \centering
    \includegraphics[width=\textwidth]{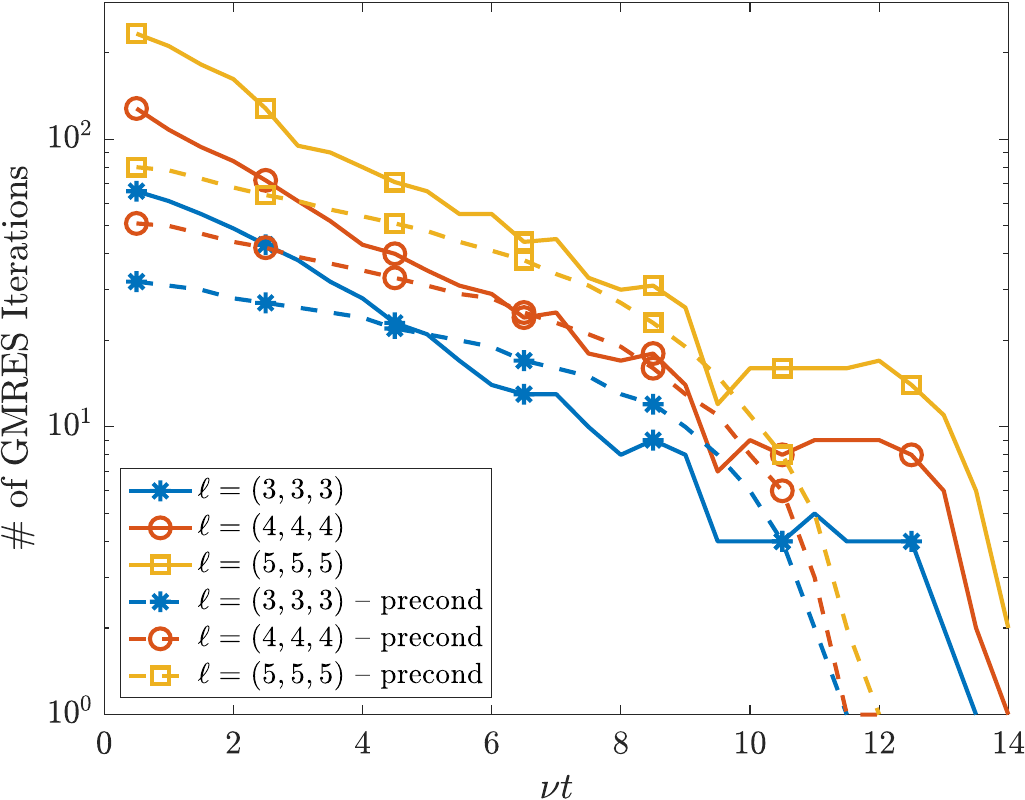}
    \captionsetup{margin=1em}
    \caption{Number of GMRES iterations per timestep with (dashed) and without (solid) a block-Jacobi preconditioner.  GMRES was restarted every 100 iterations and exited when the residual norm was less than 1e-8.  Instances when GMRES exited in zero iterations are not plotted.}
    \label{fig:relaxation:gmres_iters}
    \end{subfigure}
    \caption{Relaxation Problem -- \Cref{subsec:relax}: Plots of interest for full-grid runs with varying levels.}
\end{figure}

\Cref{fig:relaxation:conservation} plots the change in the fluid variables $n_{f}$, $u_{f}$, and $\theta_{f}$ from their initial values as a function of $\nu t$, when using the full grid; the figure clearly shows that the loss in conservation of the moments is well below the GMRES tolerance.  The error profiles for the mixed-grid and adaptive sparse-grid runs are similar, but not shown.

\Cref{fig:relaxation:gmres_iters} shows the number of GMRES iterations for each timestep for varying full-grid levels.  
The block-Jacobi preconditioner reduces the number of GMRES iterations for each simulation (dashed lines) and overall smoothly decays the iteration count as a function of timestep.  
However, the constant jump of the iteration count, in logarithmic scale, between velocity levels in both the standard GMRES and precondioned version shows that the preconditioner does not asymptotically lower the $\mathcal{O}(4^{\ell_v})$ conditioning of the diffusion term in the LB operator.    
We found that the sparse-grid's iteration count was roughly two-thirds of the full-grid for the same level.  
Additionally, we found that the adaptive sparse-grid method often included elements from level 9 grids which caused a significant increase in the number of GMRES iterations in the adaptive sparse-grid over full-grid runs with a similar number of active elements.

\begin{figure}[!ht]
    \centering
    \includegraphics[width=0.6\textwidth]{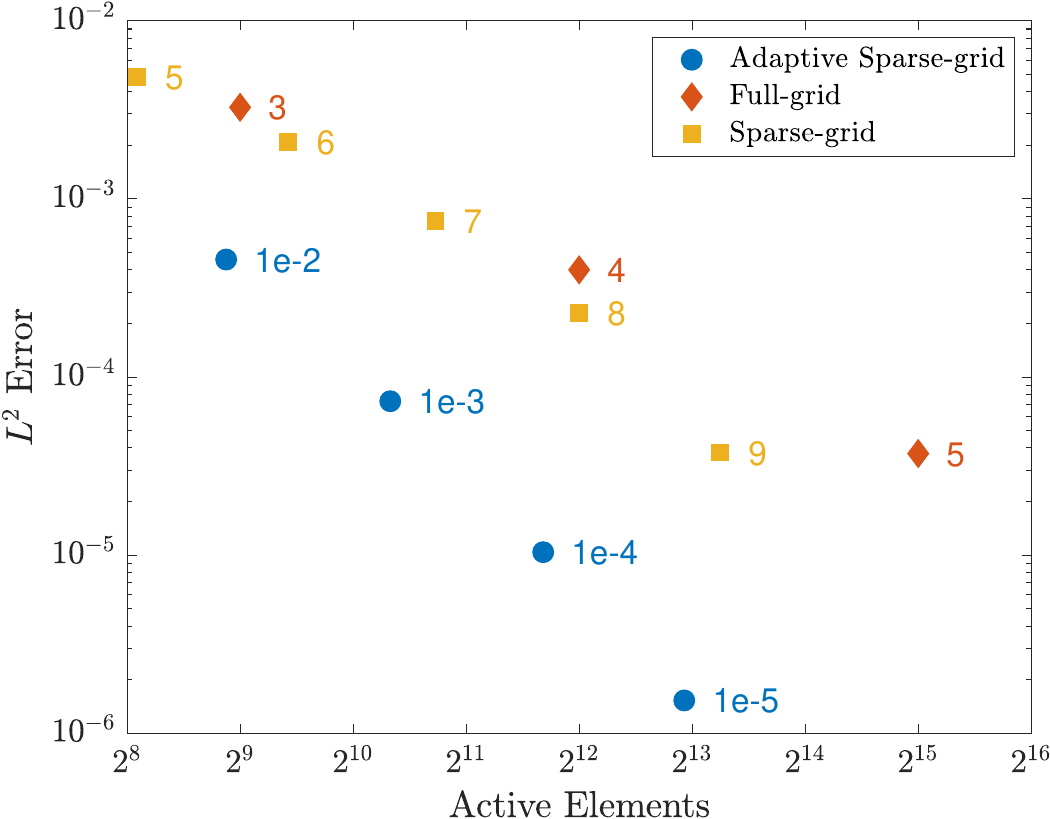}
    \caption{Relaxation Problem -- \Cref{subsec:relax}:  The $L^2$ error of solution versus the number of active elements used for the full-grid and adaptive sparse-grid runs.  The error is calculated against the analytic equilibrium in \eqref{eqn:maxwellian}.  The full- and mixed-grid runs were set at $\bm{\ell}=(\ell_v,\ell_v,\ell_v)$ where $\ell_v$ is the number by the marker. The marker next to the adaptive sparse-grid runs is the tolerance $\tau$ at which the run was set, and the adaptive run was not allowed to exceed a level of $\bm{\ell}=(9,9,9)$.}
    \label{fig:relaxation:l2_error}
\end{figure}

\Cref{fig:relaxation:l2_error} illustrates the advantages of adaptive sparse-grids over the full- and mixed-grid methods for the relaxation problem.  
The $L^{2}$ error of the relaxed distribution, relative to the analytic Maxwellian, is plotted versus the number of active elements.  
When plotted against the number of active elements, adaptive sparse grids are more accurate and asymptotically superior when compared against the other formulations. 
Additionally, the mixed-grid is comparable to the full-grid with the mixed-grid only gaining an advantage when a large number of active elements is used.  
This is not surprising as the Maxwellian, being radially symmetric, has large mixed derivatives and the coefficients to capture mixed derivative information are thrown away in the standard sparse-grid construction.  
However, adaptive sparse-grids are able to capture these mixed-derivative coefficients.

\subsection{Riemann Problem}\label{subsec:riemann}

Next, we consider a problem that includes both phase-space advection and collisions.  
The Sod shock tube problem \cite{sod_1978} is a standard test for numerical simulations of kinetic models with collisions (e.g., \cite{bennoune_etal_2008,filbetJin_2010}).  
For this test, the PDE is given by \eqref{eqn:kinetic_slab}
We consider two regimes of collisionality: The first is an intermediate regime with $\nu = 1$, and the second is a collisional regime with $\nu=10^{3}$.  
For both problems we fix $\bv\in (-6,6)^3$ and set the initial condition to a Maxwellian with moments given by:
\begin{align}\label{eqn:riemann:initial_fluid_vars}
    \begin{bmatrix}
        n_{f} \\ \bu_f \\ \theta_{f}
    \end{bmatrix} 
    =
    \begin{bmatrix}
        1 \\ \bm{0} \\ 1
    \end{bmatrix}
    ~~\text{if}~|x| \geq s_{\text{initial}};
    \qquad
    \begin{bmatrix}
        n_{f} \\ \bu_f \\ \theta_{f}
    \end{bmatrix} 
    =
    \begin{bmatrix}
        0.125 \\ \bm{0} \\ 0.8
    \end{bmatrix}
    ~~\text{if}~|x| < s_{\text{initial}}
\end{align}  
where $s_{\text{initial}}$ is the location of the initial discontinuity. 
We set the GMRES tolerance to $10^{-8}$.  

\begin{figure}[!ht]
    \centering
    \begin{subfigure}[b]{0.49\textwidth}
        \centering
        \includegraphics[width=1.1\textwidth]{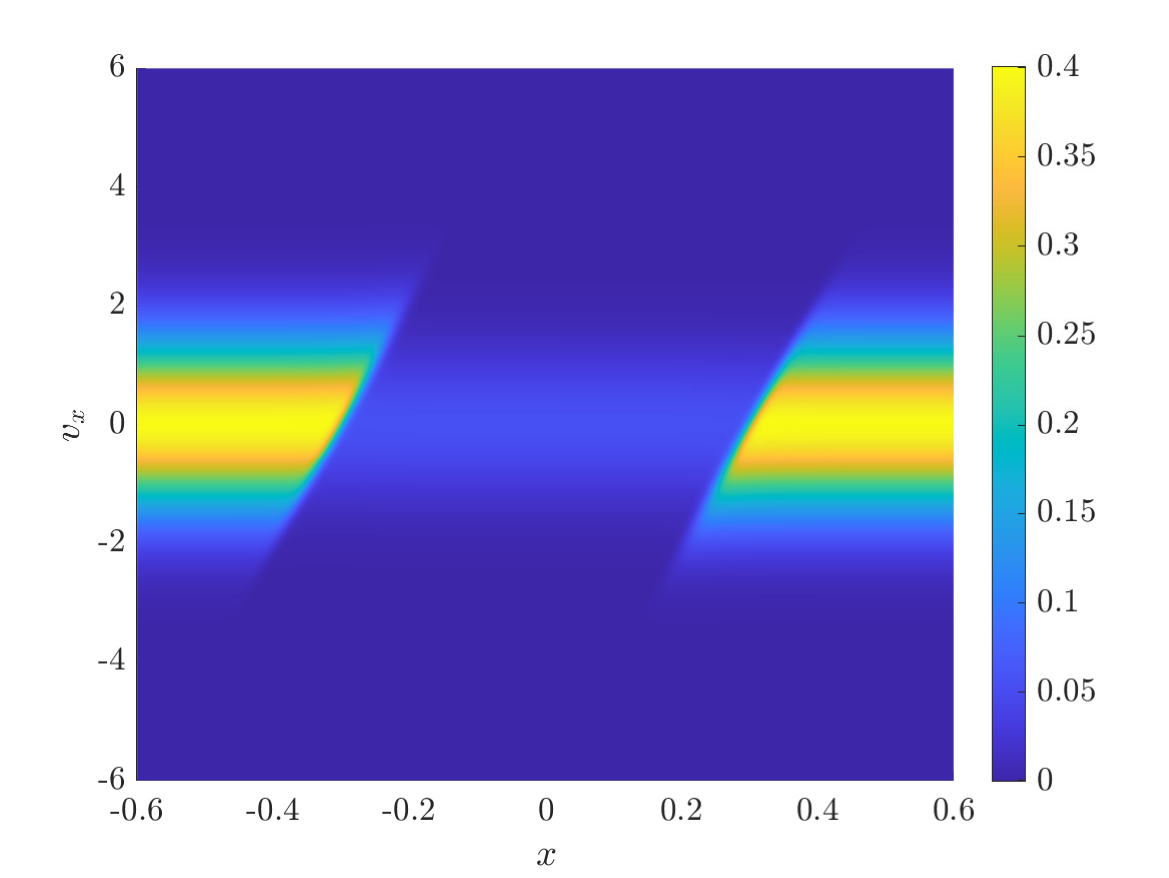}
        \caption{$\langle f\rangle_{v_y,v_z}$ -- $\nu = 1$}
        \label{fig:riemann:reference:a}
    \end{subfigure}
    \begin{subfigure}[b]{0.49\textwidth}
        \centering
        \includegraphics[width=1.1\textwidth]{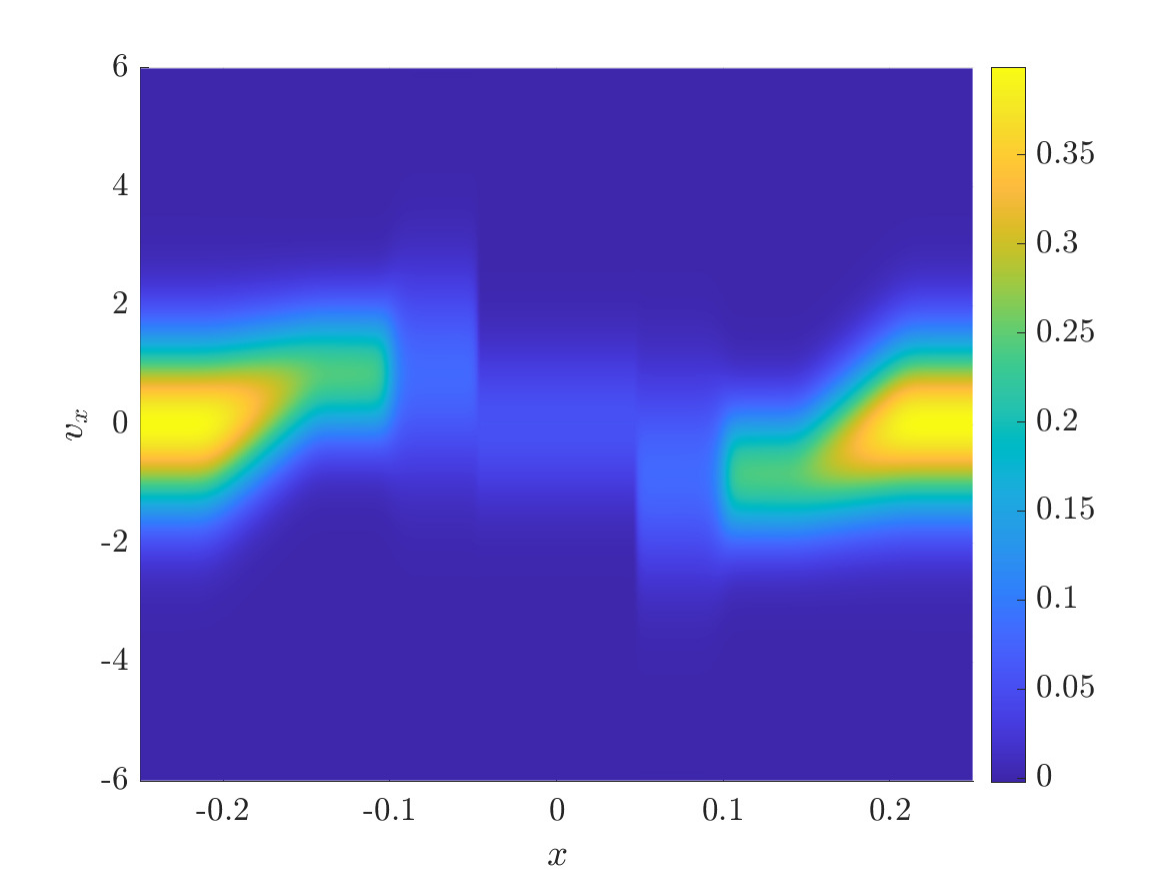}
        \caption{$\langle f \rangle_{v_y,v_z}$ -- $\nu = 10^{3}$}
        \label{fig:riemann:reference:b}
    \end{subfigure}

    \centering
    \begin{subfigure}[b]{0.49\textwidth}
        \centering
        \includegraphics[height=2.6in]{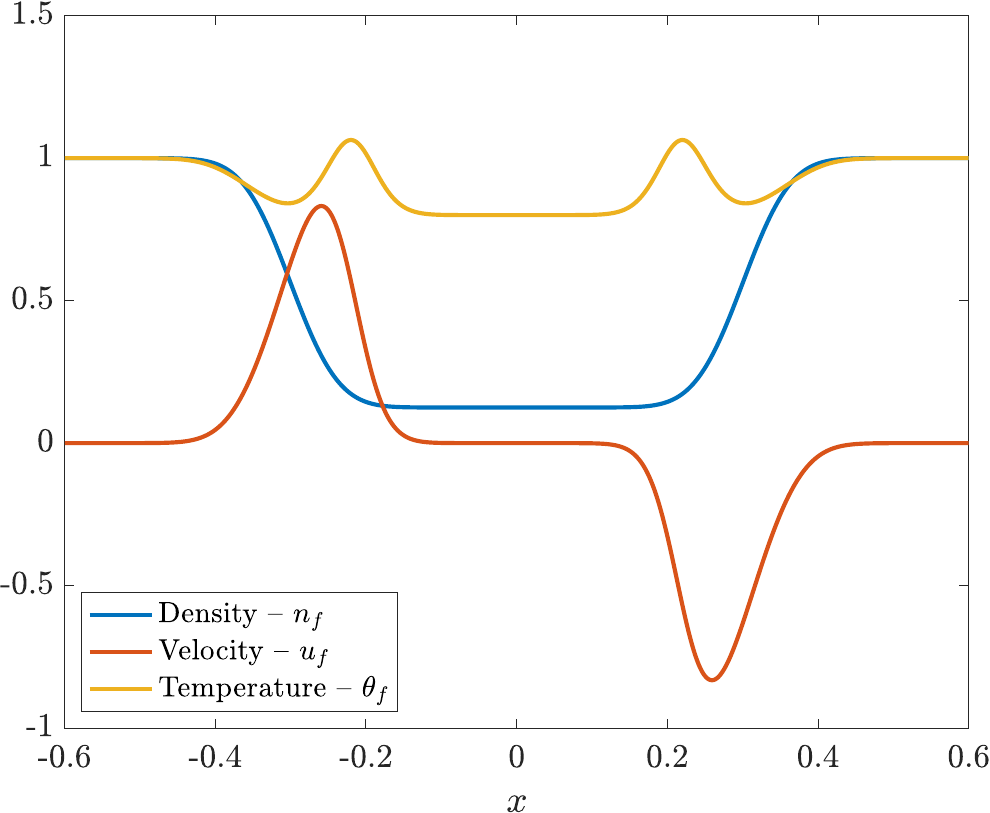}
        \caption{Fluid variables -- $\nu = 1$}
        \label{fig:riemann:reference:c}
    \end{subfigure}
    \begin{subfigure}[b]{0.49\textwidth}
        \centering
        \includegraphics[height=2.6in]{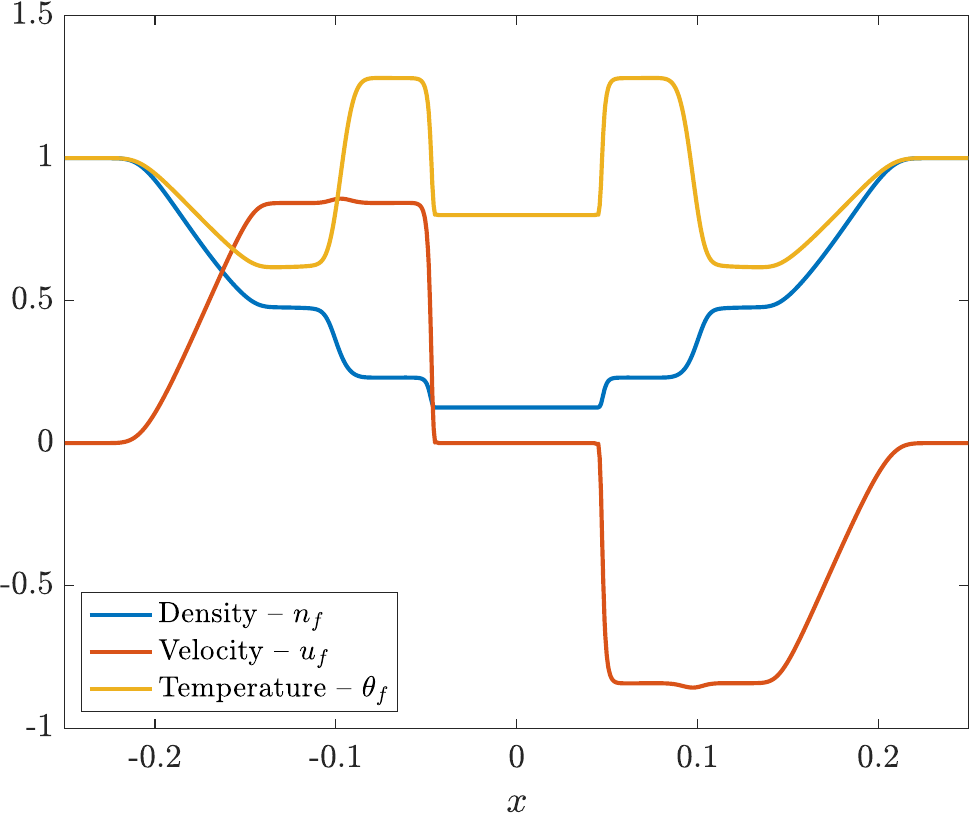}
        \caption{Fluid variables -- $\nu = 10^{3}$}
        \label{fig:riemann:reference:d}
    \end{subfigure}
    \caption{Riemann Problem -- \Cref{subsec:riemann}: Plots of the distribution and fluid variables for the which is computed using the Chu reduction model in \Cref{subsec:reduction} with $(\ell_x,\ell_v)=(9,8)$ for $\nu=1$ and $(\ell_x,\ell_v)=(8,8)$ for $\nu=10^3$.  The $\nu=1$ and $\nu=10^{3}$ plots are taken at time $t=0.04918$ and $t=0.05$ and with $s_{\text{initial}} = 0.3$ and $s_{\text{initial}}= \sfrac{9}{64}$ respectively. }
    \label{fig:riemann:reference}
\end{figure}

\Cref{fig:riemann:reference} shows plots of the distribution in the $(x,v_x)$-plane and plots of the velocity moments versus position, as obtained with the full-grid using the Chu reduction technique.  
We will use these as reference solutions when evaluating the performance of the adaptive sparse grid method.  
For moderate collisionality, i.e.~$\nu=1$, the distribution, as shown in \Cref{fig:riemann:reference:a}, deviates from the Maxwellian due to the streaming and features a discontinuity in the $(x,v_x)$ space.
Additionally, as seen in \Cref{fig:riemann:reference:c}, the streaming effect smooths out features of the fluid variables.  
In the collision dominated regime ($\nu=10^3$), the distribution, as  seen in \Cref{fig:riemann:reference:b}, remains close to a local Maxwellian parameterized by the local fluid variables in \Cref{fig:riemann:reference:d}.

\begin{figure}[!ht]
    \centering
    \begin{subfigure}[b]{0.49\textwidth}
        \centering
        \includegraphics[width=\textwidth]{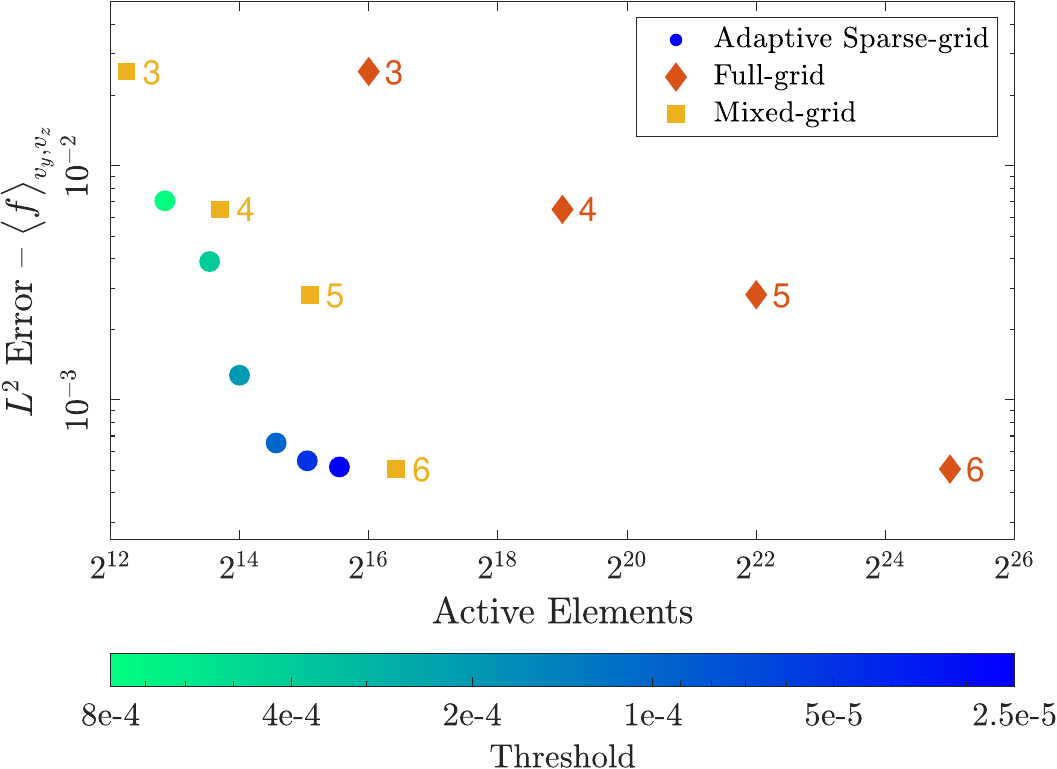}
        \caption{$\Big\|\big<f\big>_{v_y,v_z}-g_1\Big\|_{L^2}$}
        \label{fig:riemann:nu1e0_dof_vs_err:1}
    \end{subfigure}
    \begin{subfigure}[b]{0.49\textwidth}
        \centering
        \includegraphics[width=\textwidth]{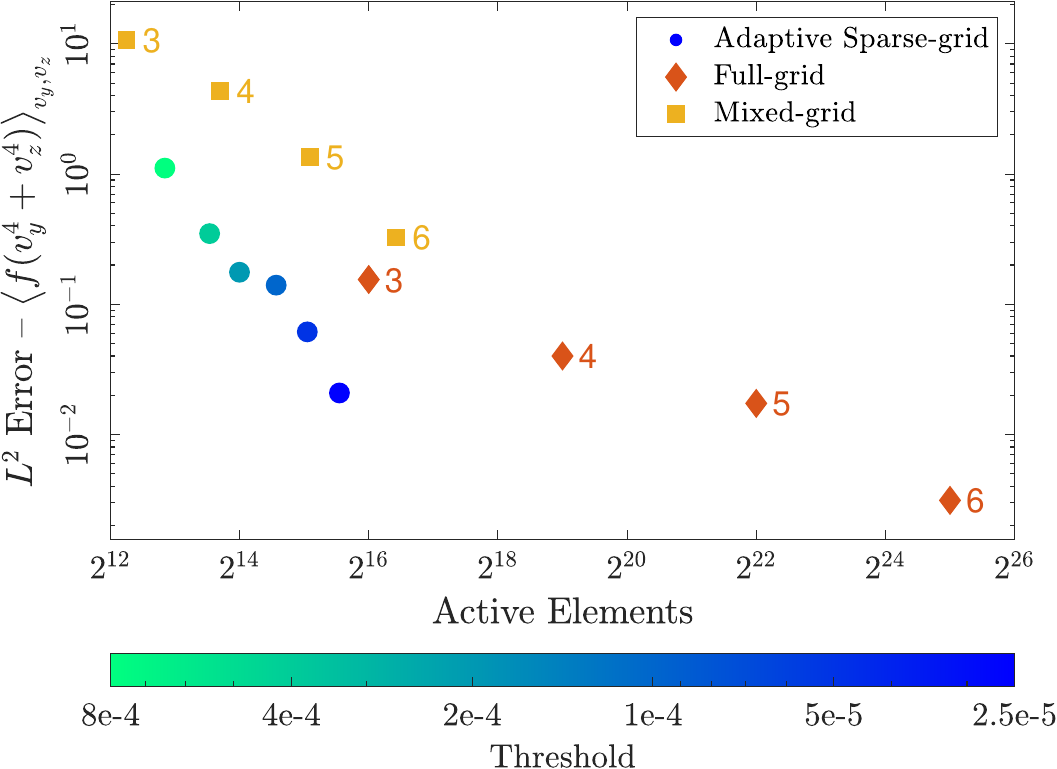} \caption{$\Big\|\big<f(v_y^4+v_z^4)\big>_{v_y,v_z}-g_3\Big\|_{L^2}$}
        \label{fig:riemann:nu1e0_dof_vs_err:2}
    \end{subfigure}
    \caption{Riemann Problem -- \Cref{subsec:riemann} -- $\nu=1$: Errors of the distribution and fluid variables  at $t=0.01$ for the $1x3v$ Riemann problem in \Cref{subsec:riemann} with $\nu=1$.  All errors are measured against the full-grid solution at level $\bm{\ell}=(9,8,8,8)$ (see \Cref{fig:riemann:reference:a}).  All adaptive sparse-grid runs are capped at $\bm{\ell}=(7,6,6,6)$.  The full- and mixed-grid runs use $\bm{\ell}=(7,\ell_v,\ell_v,\ell_v)$ where $\ell_v$ is the symbol by each marker. The adaptive sparse-grid method performs well in both cases while the mixed-grid method is accurate only in the low-order moment.}
    \label{fig:riemann:nu1e0_dof_vs_err}
\end{figure}

When comparing results obtained with different grids, we first consider the case of $\nu=1$, and we set $x\in(-0.6,0.6)$, $s_{\text{initial}}=0.3$, final time $T=0.04918$, and time step $\Delta t = 2.3419\times 10^{-4}$.  
Our reference solution is the full-grid solution of level $\bm{\ell}=(9,8,8,8)$, displayed in the left panels in \Cref{fig:riemann:reference}.
\Cref{fig:riemann:nu1e0_dof_vs_err} shows the error versus the number of active elements for $g_1$ and $g_3$ (defined in \Cref{subsec:reduction}).  
It is shown in \Cref{fig:riemann:nu1e0_dof_vs_err:1} that the mixed-grid yields the same error as the full-grid -- for the same velocity space resolution level $\ell_{v}$.  
This is because $g_1$ is embedded in the mixed-grid as mentioned in \Cref{subsec:choice_of_grids}. 
The adaptive sparse-grid error saturates at the level of the mixed-grid error when $\ell_{v}=6$, but with about 50\% fewer active elements.  
The saturation is because the adaptive grid is not allowed to refine  past level $\bm{\ell}=(7,6,6,6)$ in the hierarchy (see \Cref{subsec:choice_of_grids}) and therefore the associated error will not be significantly lower than the full-grid of level $\bm{\ell}=(7,6,6,6)$.
When viewing the same plot for the higher-order moment $g_3$ in \Cref{fig:riemann:nu1e0_dof_vs_err:2}, we see the degradation in the mixed-grid method when compared to the full-grid and adaptive sparse-grid methods.  
While the slope in the error from the mixed-grid method is steeper than the full-grid method, its error constant is significantly larger.
Additionally, the adaptive sparse-grid method is significantly better than both the mixed-grid and full-grid methods with respect to both the slope and error constant.  

\begin{figure}[!ht]
    \centering
    \begin{subfigure}[b]{0.32\textwidth}
        \centering
        \includegraphics[width=0.95\textwidth]{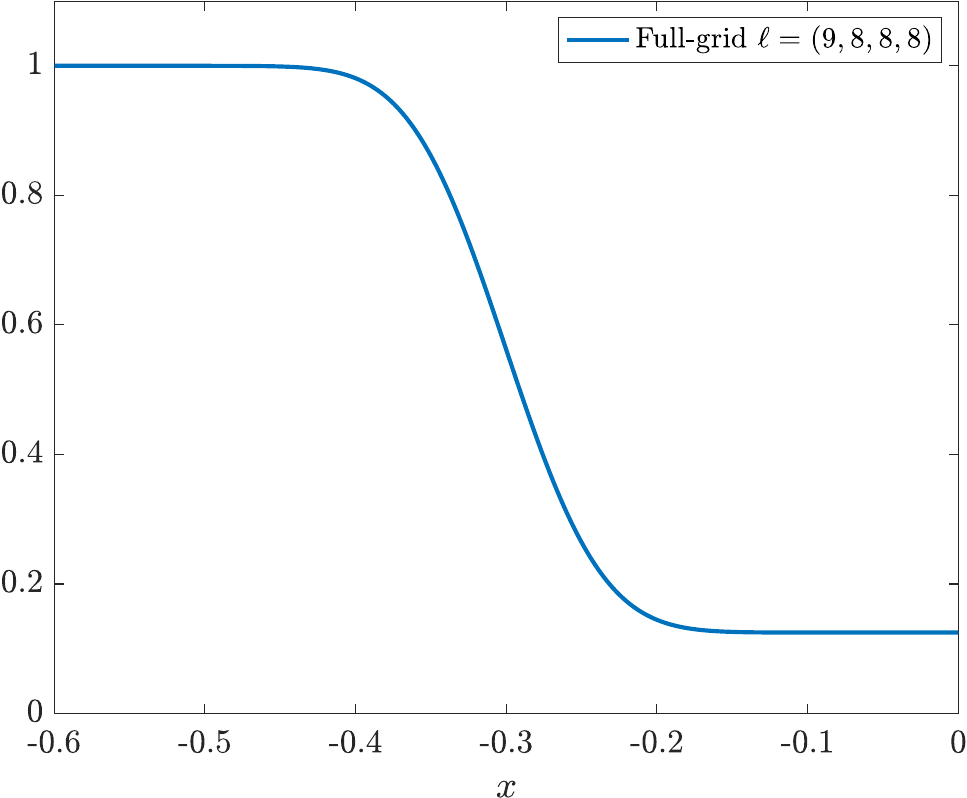}
        \caption{$n_f$}
        \label{fig:riemann:nu1e0_moment_error:a}
    \end{subfigure}
    \begin{subfigure}[b]{0.32\textwidth}
        \centering
        \includegraphics[width=\textwidth]{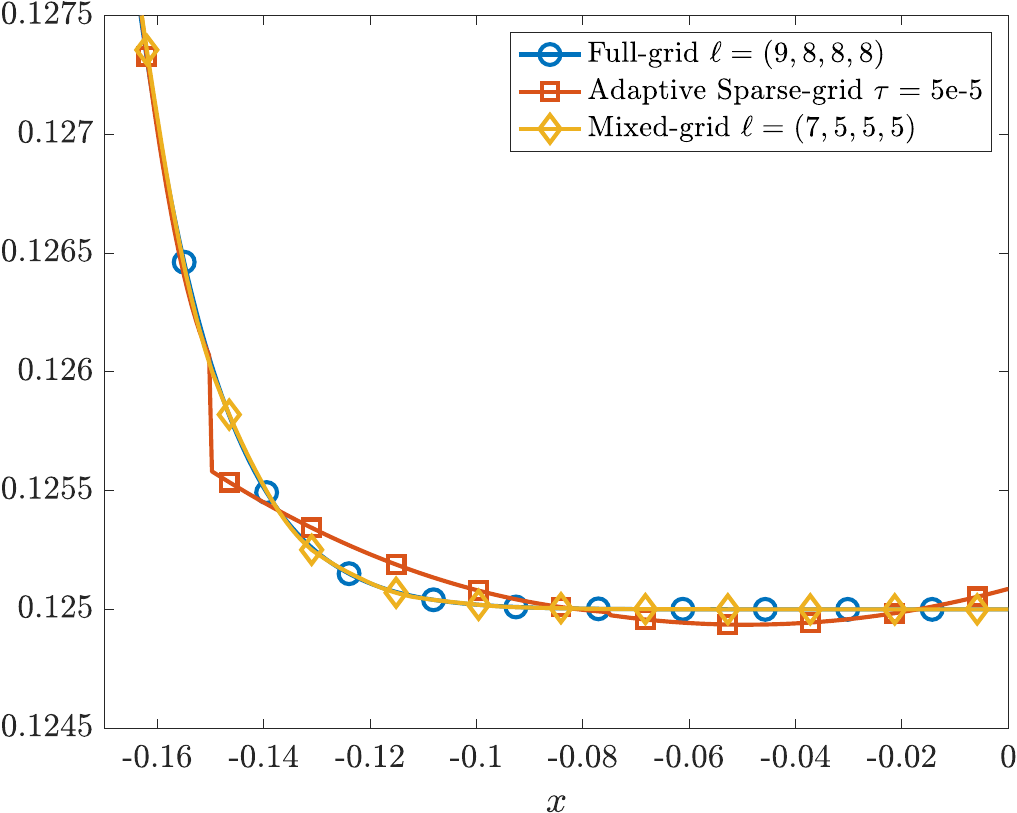}
        \caption{$n_f$ -- slice for $x\in(-0.17,0)$}
        \label{fig:riemann:nu1e0_moment_error:b}
    \end{subfigure}
    \begin{subfigure}[b]{0.32\textwidth}
        \centering
        \includegraphics[width=1.025\textwidth]{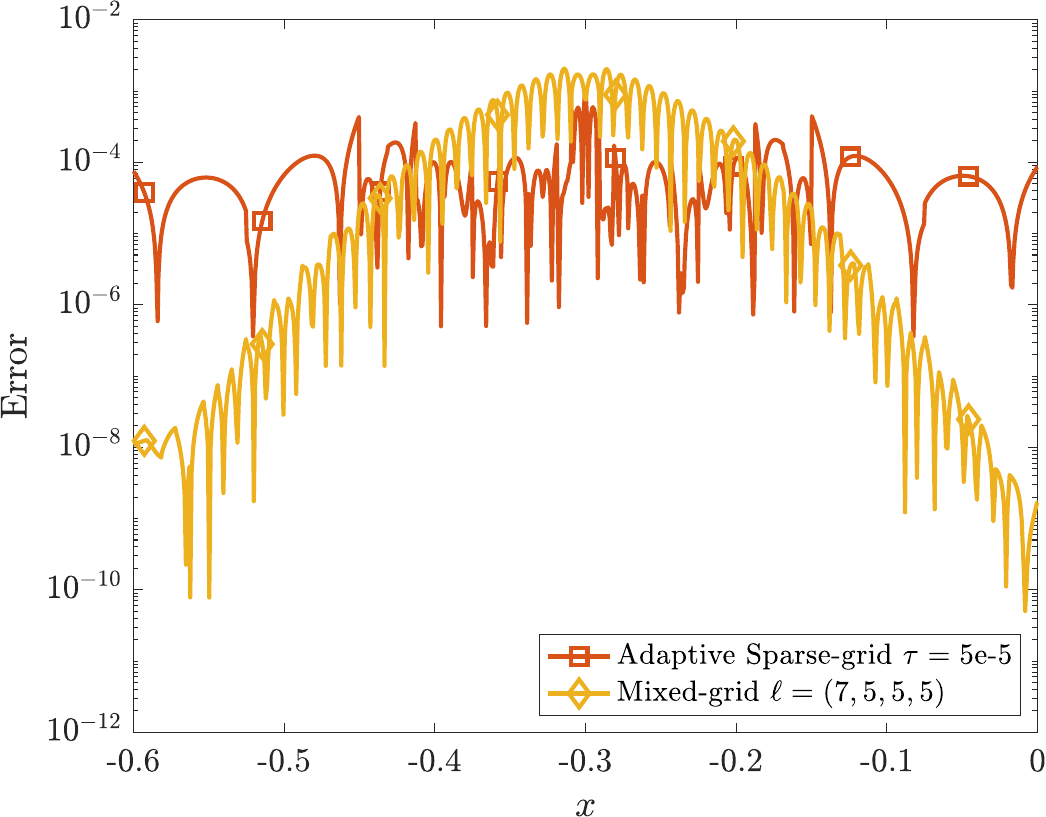}
        \caption{Error -- $n_f$}
        \label{fig:riemann:nu1e0_moment_error:d}
    \end{subfigure}
    \caption{Riemann Problem -- \Cref{subsec:riemann} -- $\nu=1$: Plots of the density and error to the reference density for the $1x3v$ Riemann problem with $\nu=1$ and $t=0.04918$. The reference density is calculated with the full-grid method at level $\bm{\ell}=(9,8,8,8)$ (see \Cref{fig:riemann:reference:c}).  The adaptive sparse-grid solution is not allowed to be refined beyond level $\bm{\ell}=(7,6,6,6)$.  The adaptive sparse-grid method equally spaces out the error in physical space while the mixed-grid is only accurate in the constant regions of the density.}
    \label{fig:riemann:nu1e0_moment_error}
\end{figure}

\Cref{fig:riemann:nu1e0_moment_error} shows the particle density $n_{f}$ (left and middle panels) and the pointwise error of the particle density (right panel) for a mixed-grid and an adaptive sparse-grid model with a similar number of active elements.  
\Cref{fig:riemann:nu1e0_moment_error:a} shows that the density appears to be relatively constant in $x$ toward the edges of the plot.  
When zooming in on a smaller $x$-range near the right edge, see \Cref{fig:riemann:nu1e0_moment_error:b}, it becomes clear that the density obtained with the adaptive sparse-grid features a discontinuity (around $x=-0.15$) and exhibits more spatial variation when compared to the full-grid and mixed-grid solutions.  
This is primarily caused by the adaptive method uniformly distributing the error across the spatial domain, and this is further evidenced in the error plot (see~\Cref{fig:riemann:nu1e0_moment_error:d}), where the error in $x$ is much more uniform across the spatial domain for the adaptive sparse-grid than it is with the mixed-grid method.  
In the mixed-grid method, where each DOF in $x$ is attached with the same sparse-grid in $v$, the moment errors are much smaller away from the wave regions, i.e., the regions where the moments are constant.

\begin{figure}[!ht]
    \centering
    \begin{subfigure}[b]{0.49\textwidth}
        \centering
        \includegraphics[width=\textwidth]{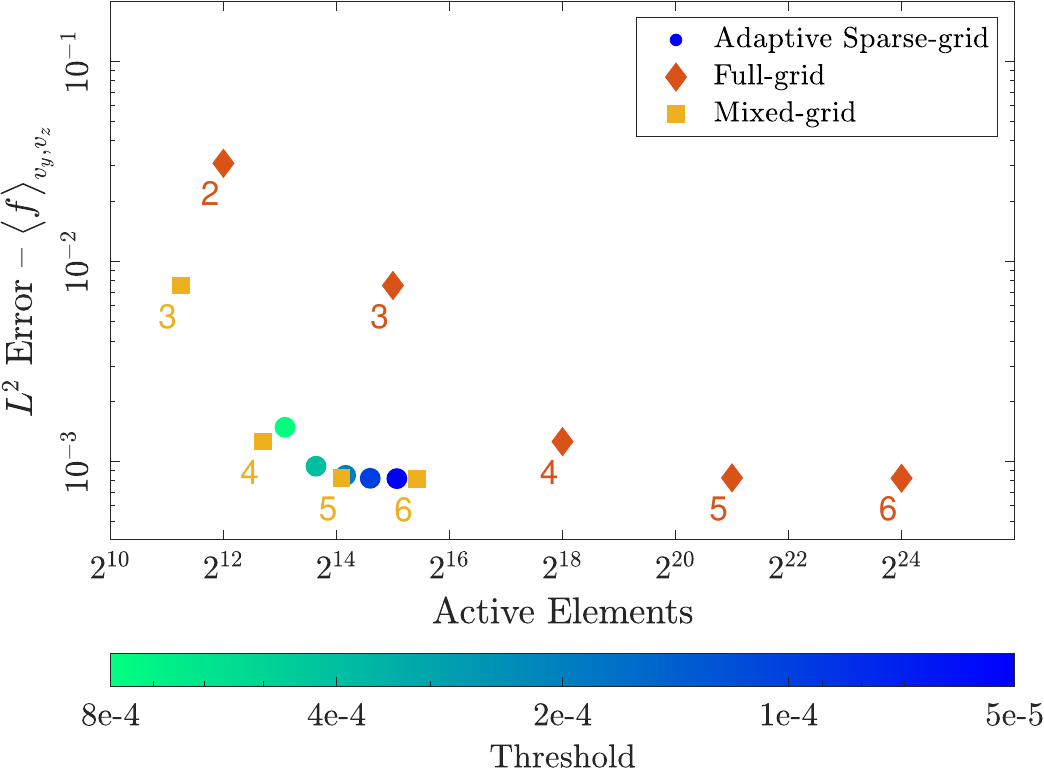}
        \caption{$\Big\|\big<f\big>_{v_y,v_z}-g_1\Big\|_{L^2}$}
        \label{fig:riemann:nu1e3_dof_vs_err:a}
    \end{subfigure}
    \begin{subfigure}[b]{0.49\textwidth}
        \centering
        \includegraphics[width=\textwidth]{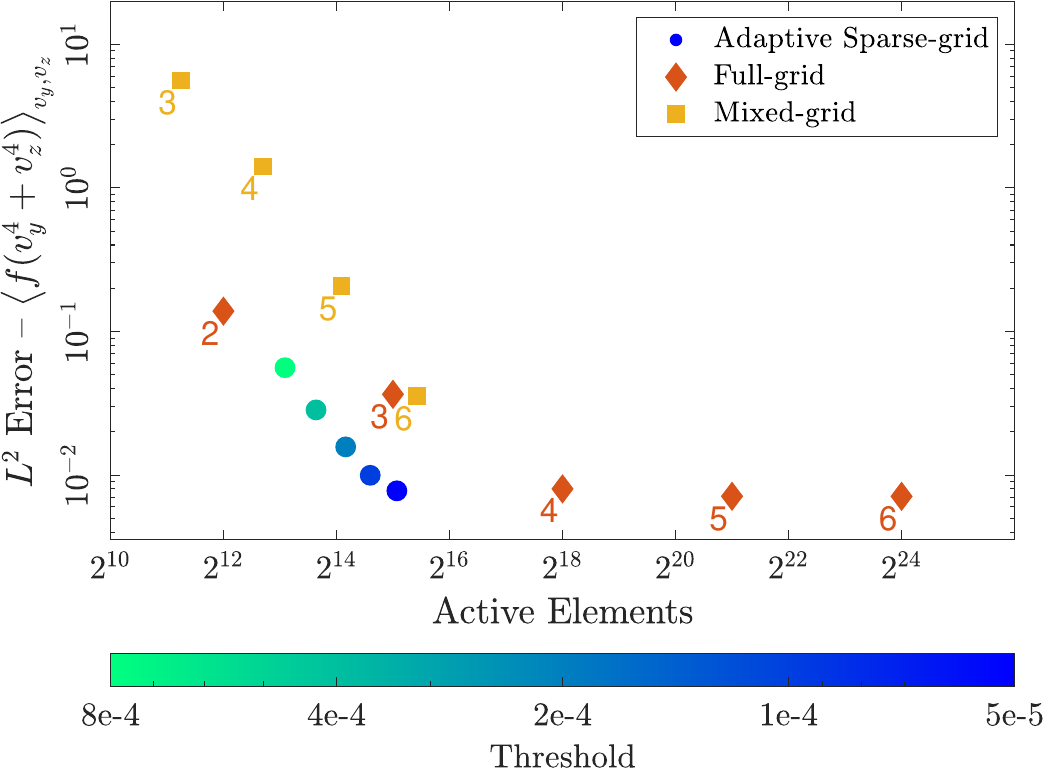}
        \caption{$\Big\|\big<f(v_y^4+v_z^4)\big>_{v_y,v_z}-g_3\Big\|_{L^2}$}
        \label{fig:riemann:nu1e3_dof_vs_err:b}
    \end{subfigure}
    \caption{Riemann Problem -- \Cref{subsec:riemann} -- $\nu=10^3$: Errors of the distribution at $t=0.05$ for the $1x3v$ Riemann problem in \Cref{subsec:riemann} with $\nu=10^3$.  All errors are measured against the full-grid solution at level $\bm{\ell}=(8,8,8,8)$ (see \Cref{fig:riemann:reference:b}).  All adaptive sparse-grids are capped at level $\bm{\ell}=(6,6,6,6)$.  The full- and mixed-grid levels are given by $\bm{\ell}=(6,\ell_v,\ell_v,\ell_v)$ where $\ell_v$ is the symbol to the lower left of the marker.  The quick saturation of the error is due to smoothness in velocity and the discontinuities in the fluid variables (see \Cref{fig:riemann:reference:d}).  The adaptive sparse-grid method performs well in both cases while the mixed-grid method is accurate only in the low-order moment.}
    \label{fig:riemann:nu1e3_dof_vs_err}
\end{figure}

Next, we consider the case with $\nu=10^{3}$.  
Here we set $x\in(-0.25,0.25)$, $s_{\text{initial}}=\sfrac{9}{64}$, $T=0.05$, and $\Delta t=2\times 10^{-4}$.  
\Cref{fig:riemann:nu1e3_dof_vs_err} shows the error of $g_1$ and $g_3$ against the number of active elements.  
In this higher-collisional regime, the distribution is much smoother in velocity, and the $L^2$ error saturates sooner than when $\nu=1$.  This saturation is due to the dominant error that appears near the discontinuities in the $x$-domain (see \Cref{fig:riemann:reference:b,fig:riemann:reference:d}). 
In \Cref{fig:riemann:nu1e3_dof_vs_err:a}, the mixed-grid and adaptive sparse-grid methods are very similar.
At saturation, the number of active elements for the mixed-grid and adaptive sparse grid, around $2^{14}$ are approximately 128 times fewer than the number of active elements in the full-grid, which is $2^{21}$.

When looking at the error in the higher-order moments, \Cref{fig:riemann:nu1e3_dof_vs_err:b}, we observe a separation in the performance of the mixed-grid and adaptive sparse-grid that is similar to the $\nu=1$ case.
However, in this case, the adaptive sparse-grid method has nearly hit saturation while the mixed-grid with $\ell_v=6$ is still not at saturation.  In particular, the grouping of the errors for mixed-grid $\ell_v=6$ and the full-grid $\ell_v=3$ is similar to the grouping in the relaxation case (see \Cref{fig:relaxation:l2_error}).  This shows that the dominant error in the mixed-grid method is the lack of velocity resolution sufficient to capture the local Maxwellian behavior of the distribution.

\begin{figure}[!ht]
    \centering
    \begin{subfigure}[b]{0.32\textwidth}
        \centering
        \includegraphics[width=\textwidth]{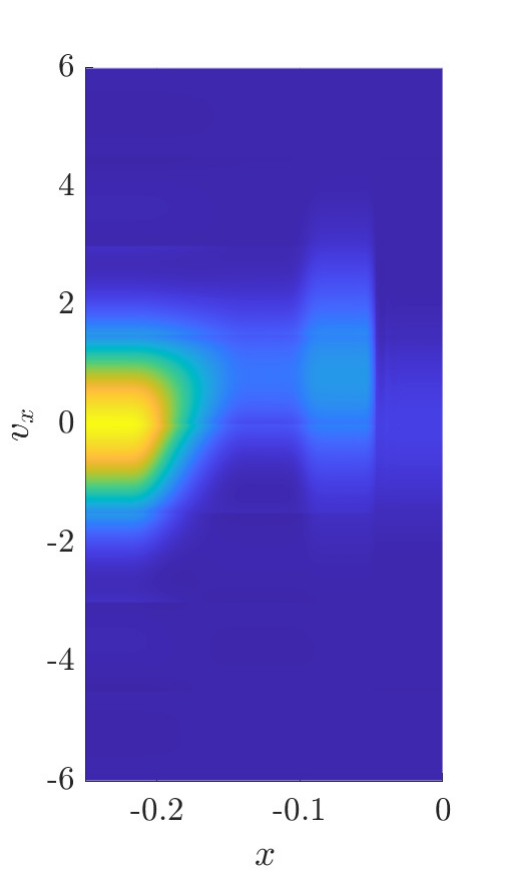}
        \caption{Full-grid $\bm{\ell}=(6,3,3,3)$}
        \label{fig:riemann:nu1e3_2d_errs:a}
    \end{subfigure}
    \begin{subfigure}[b]{0.32\textwidth}
        \centering
        \includegraphics[width=\textwidth]{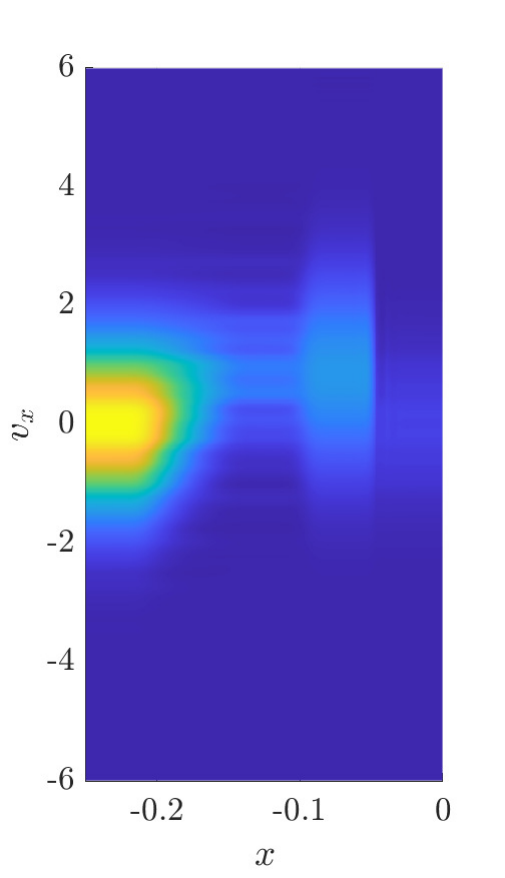}
        \caption{Mixed-grid $\bm{\ell}=(6,6,6,6)$}
        \label{fig:riemann:nu1e3_2d_errs:b}
    \end{subfigure}
    \begin{subfigure}[b]{0.32\textwidth}
        \centering
        \includegraphics[width=\textwidth]{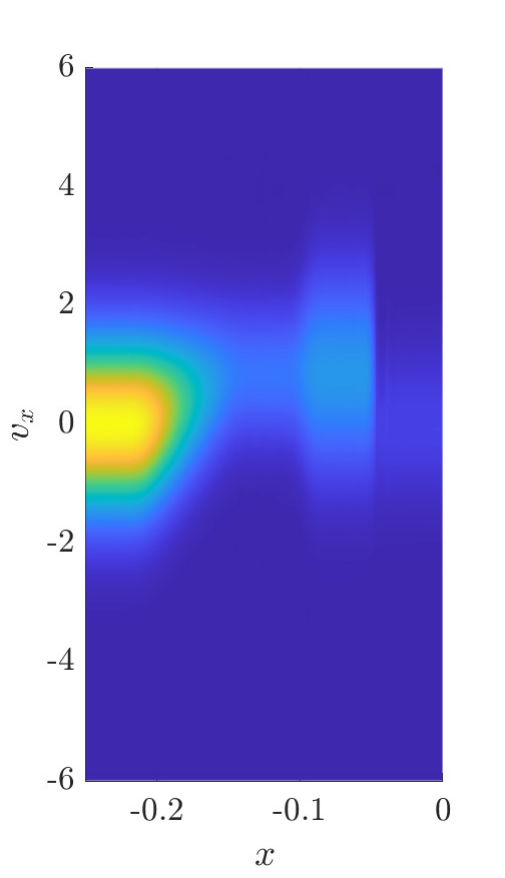}
        \caption{Adaptive Sparse-Grid $\tau=5 \times 10^{-5}$}
        \label{fig:riemann:nu1e3_2d_errs:c}
    \end{subfigure}
    \caption{Riemann Problem -- \Cref{subsec:riemann} -- $\nu=10^3$: Plots of $\langle f(v_y^4+v_z^4) \rangle_{v_y,v_z}$ for the $1x3v$ Riemann Problem with $\nu=10^3$ and $t=0.05$ in the $(x,v_x)$ plane for $x\in(-0.25,0)$ and $v_x\in(-6,6)$.  The adaptive sparse-grid was not allowed to refine past level $\bm{\ell}=(6,6,6,6)$.  The artifacts seen in the full-grid and mixed-grid solutions are not found in the adaptive- sparse grid solution.}
    \label{fig:riemann:nu1e3_2d_errs}
\end{figure}

We include a plot of the fourth-order moment $g_3$ in the $(x,v_x)$-plane for each grid type, each having a similar number of degrees of freedom, in \Cref{fig:riemann:nu1e3_2d_errs}.  
The full-grid solution, \Cref{fig:riemann:nu1e3_2d_errs:a}, exhibits discontinuities on element interfaces in the velocity dimension (due the discontinuous basis) while the mixed-grid moment, \Cref{fig:riemann:nu1e3_2d_errs:b}, is oscillatory in the region immediately left of the contact line, i.e.~$x\in(-0.15,-0.1)$.  
The solution obtained with the adaptive sparse-grid, \Cref{fig:riemann:nu1e3_2d_errs:c}, is the most accurate of the three and does not suffer from either of the previously mentioned artifacts.  

\subsection{Collisional Landau Damping}\label{subsec:landau}

Finally, we consider a version of the collisional Landau damping test (e.g., \cite{crestetto_etal_2012,hakim_etal_2020,francisquez_etal_2020a}), which involves phase-space advection of charged particles, influenced by a self-consistent electric field and particle collisions.  
The PDEs solved in this test are given by the VPLB system in \eqref{eqn:kinetic_slab} and \eqref{eqn:Poisson_slab}.

The $1x3v$ phase-space domain is given by $x\in(-2\pi,2\pi)$ and $\bv\in(-6,6)^3$, and the model is evolved to the final time $T=50$.  
The initial condition is set as Maxwellian with a small spatial perturbation so that the velocity moments are $n_{f} = 1+10^{-4}\cos(\tfrac{x}{2})$, $u_{f}=0$, $\theta_{f} = 1$.  
The timestep taken depends on the spatial resolution and $\max|v_{x}|$, and is taken as $\Delta t = \tfrac{0.75}{30}\Delta x$, where $\Delta x = \tfrac{4\pi}{2^{\ell_x}}$.

\begin{figure}[!ht]
    \centering
    \begin{subfigure}[b]{0.32\textwidth}
        \centering
        \includegraphics[width=\textwidth]{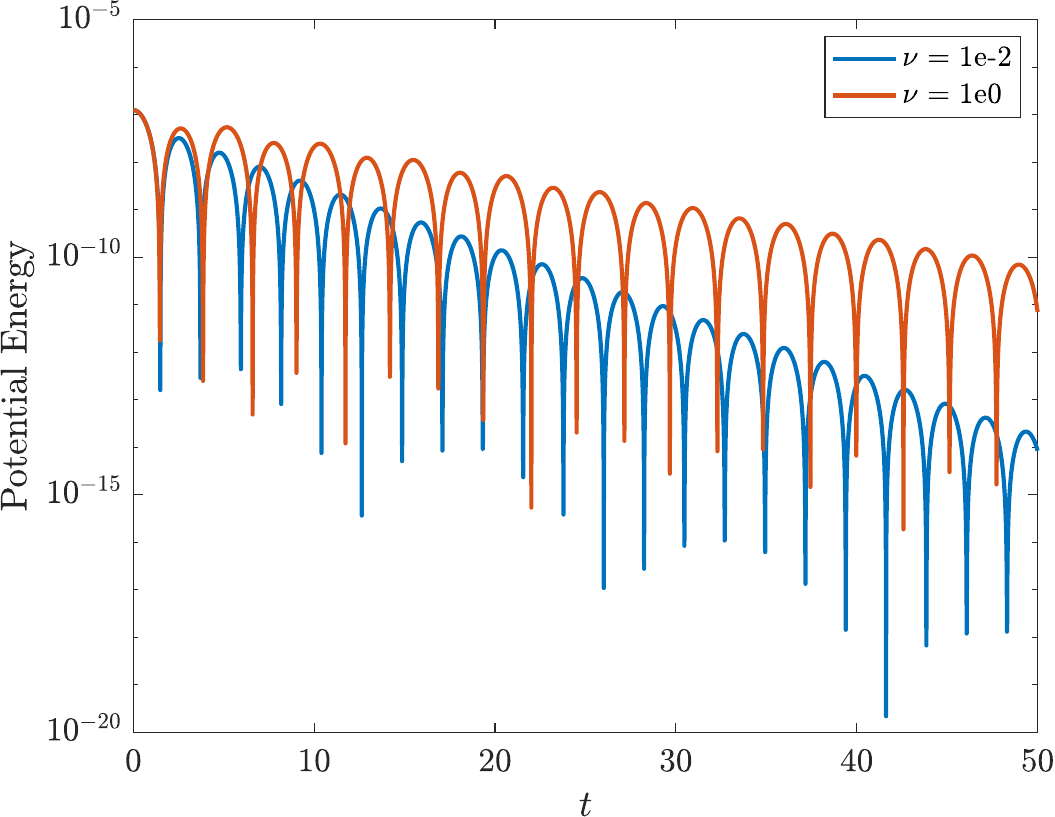}
        \caption{The electric potential energy versus time for two collision frequencies $\nu$. \\}
        \label{fig:landau:nu_plot:FG_runs}
    \end{subfigure}
        \begin{subfigure}[b]{0.32\textwidth}
        \centering
        \includegraphics[width=\textwidth]{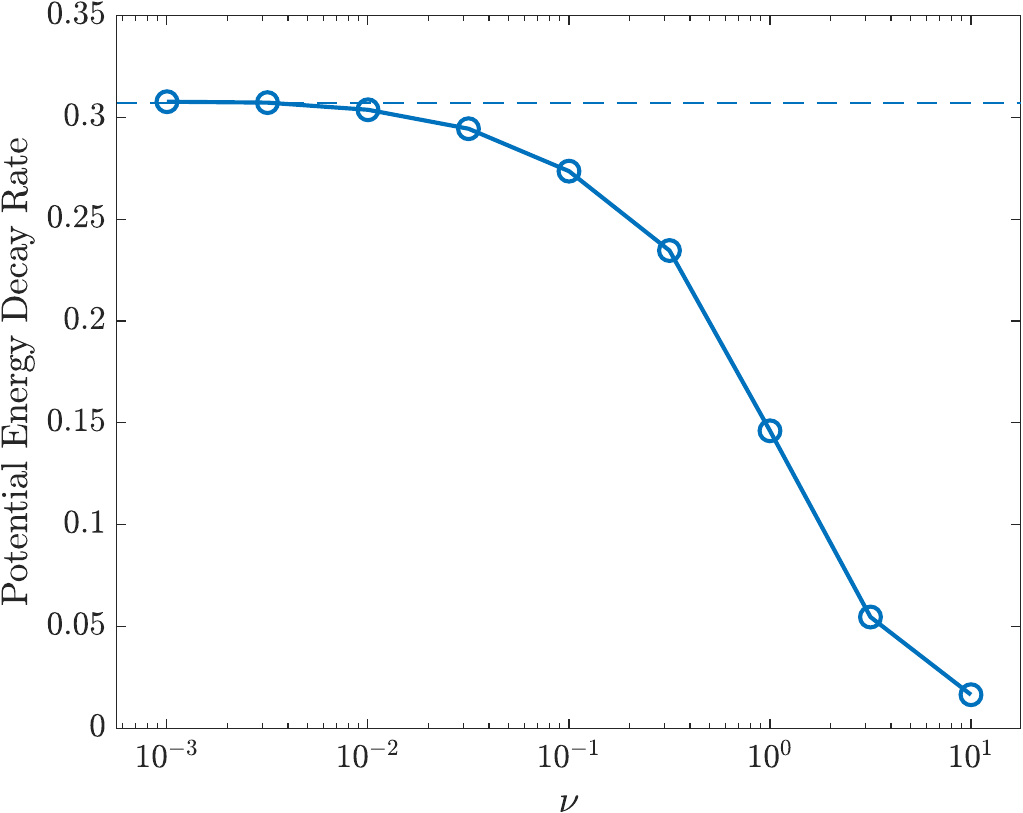}
        \caption{Exponential decay rate $\gamma$ for the potential energy as a function of $\nu$. The dashed line corresponds to $\gamma=0.307$. }
        \label{fig:landau:nu_plot:decay}
    \end{subfigure}
    \begin{subfigure}[b]{0.32\textwidth}
        \centering
        \includegraphics[width=\textwidth]{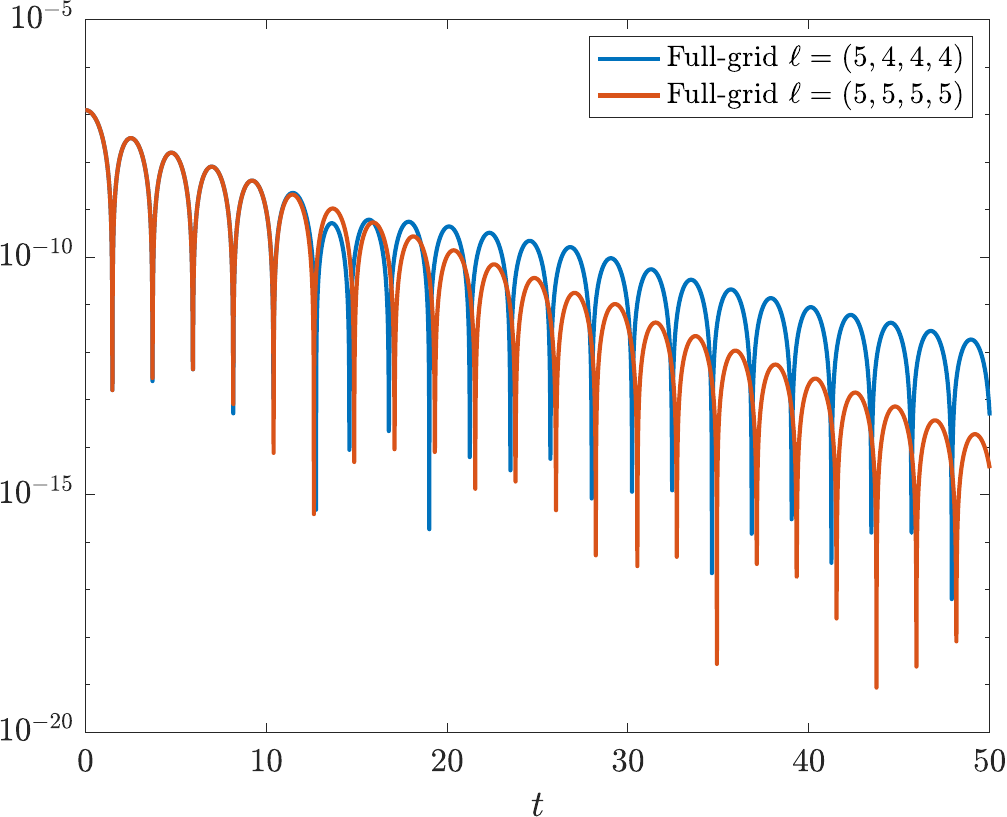}
        \caption{Plot of the electric potential for $\nu=10^{-2}$ and velocity resolutions $\ell_v=4$ and $\ell_v=5$.}
        \label{fig:landau:nu_plot:damping}
    \end{subfigure}
    \caption{Collisional Landau Problem -- \Cref{subsec:landau}: Plots demonstrating collisional Landau damping.  All runs use the Chu reduction method of \Cref{subsec:reduction}.  The levels set are $\ell_x=5$ and $\ell_v=6$ except in \Cref{fig:landau:nu_plot:damping}.}
    \label{fig:landau:nu_plot}
\end{figure}

In the collisionless case, the Landau damping problem is characterized by exponential decay of the potential energy with time, $E_{\rm Pot}(t)=\frac{1}{2}\int_{\Omega_{x}}E^{2}\dx x\propto\exp(-\gamma t)$, where the damping rate is $\gamma\approx0.307$ \cite{arberVann_2002}.  
Moreover, with evolving time, the solution will exhibit increasingly smaller-scale structures about the Maxwellian that eventually become unresolved with fixed or finite resolution (see \cite{endeveHauck_2022}).
With collisions, the damping rate decreases with increasing collision frequency (e.g., \cite{crestetto_etal_2012,hakim_etal_2020}), tending to zero in the Euler--Poisson limit ($\nu\to\infty$).  

\Cref{fig:landau:nu_plot:FG_runs} shows the potential energy versus time, as obtained with the full-grid method, for $\nu=10^{-2}$ (blue) and $\nu=1$ (red).  
\Cref{fig:landau:nu_plot:decay} shows numerically determined damping rates as a function of collision frequency.  
These results were obtained with the full-grid method using the Chu reduction technique.  
The damping rate is determined by a least squares fit using the local maxima of the potential energy.  
For small collision frequencies, the damping rate tends to the expected result in the collisionless limit indicated by the horizontal dashed line.  
The damping rate drops rapidly for $\nu\gtrsim0.3$, and has dropped to about $0.01$ for $\nu=10$.  
\Cref{fig:landau:nu_plot:damping} compares the evolution of the potential energy versus time for the $\nu=10^{-2}$ case with two different velocity resolutions; $\ell_{v}=4$ (blue) and $\ell_{v}=5$ (red).  
For the simulation with the coarser velocity resolution, the damping rate is consistent with the analytic prediction until $t\approx10$.  
For $t\gtrsim10$, the potential energy increases briefly with time before decreasing again with a modified damping rate.
For the finer velocity resolution, the damping rate stays constant at the correct value for all times.  
Based on this observation, we consider $\ell_{v}=5$ the minimum resolution needed to perform satisfactory on this test when $\nu=10^{-2}$.  
We performed a similar comparison with $\nu=1$, which revealed that $\ell_{v}=4$ is sufficient for this case.
In the following, we consider the two cases: $\nu=10^{-2}$ (low collisionality) and $\nu=1$ (moderate collisionality), in more detail to compare the adaptive sparse-grid method against the full-grid method.  
Due to the embedding of the $1x1v$ full-grid into the $1x3v$ mixed-grid as discussed in \Cref{subsec:choice_of_grids}, the electric field $E$ is similar for the full- and mixed-grids of the same level.
For this reason, the mixed-grid results are omitted.  

\begin{figure}[!ht]
    \centering
    \includegraphics[width=0.6\textwidth]{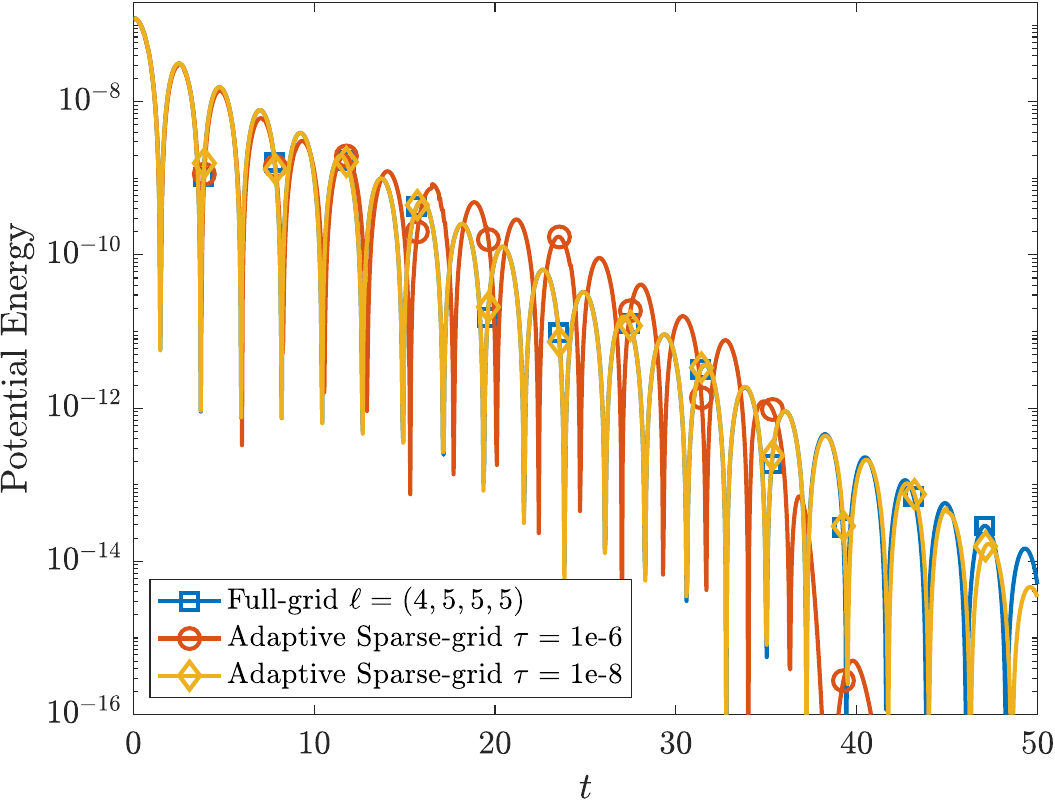}
    \caption{Collisional Landau Problem -- \Cref{subsec:landau}: Plot of the potential energy with $\nu=10^{-2}$.    
    The adaptive sparse grid is not allowed to refine past level $\bm{\ell}=(4,5,5,5)$, and the GMRES tolerance is set to $10^{-14}$.  A tolerance of $\tau=10^{-6}$ is not sufficient to capture the proper decay.  The tolerance of $\tau=10^{-8}$ agrees quite well with the full-grid solution except for a slight deviation at longer times.}
    \label{fig:landau:nu1e-2:asg_vs_fg}
\end{figure}

\Cref{fig:landau:nu1e-2:asg_vs_fg} compares adaptive sparse-grid against full-grid for the low collisionality case by plotting the potential energy versus time.  
The full-grid run with $\bm{\ell}=(4,5,5,5)$, used as reference in \Cref{fig:landau:nu1e-2:asg_vs_fg}, is in close agreement with the full-grid run with $\bm{\ell}=(5,5,5,5)$ plotted in the right panel of \Cref{fig:landau:nu_plot:damping}.
When the tolerance for refinement is $\tau=10^{-6}$, the adaptive results agree with the full-grid results up to about $t=10$.
For later times, the resolution allowed by the threshold is not sufficient to capture the correct damping of the potential energy.  
Past $t=35$, the solution coarsens to only global elements in $x$, i.e.~$\ell_x=0$, which forces the electric field to zero before refinement, and causes unreliable behavior in the potential energy.
When the tolerance is reduced to $\tau=10^{-8}$, the adaptive sparse-grid is in better agreement with the full-grid throughout the simulation, although some deviations near the end are observed.  
For the $\tau=10^{-6}$ case, the number of active elements stays around $1.1\times10^{4}$ throughout, while for the $\tau=10^{-8}$ case the number of active elements starts out around $3\times10^{4}$, which drops steadily to about $2.5\times10^{4}$ at the end of the simulation.  
In comparison, the full-grid with $\bm{\ell}=(4,5,5,5)$, the maximum allowed for the adaptive spares-grid, consists of about $5.2\times10^{5}$ elements.  
Thus, the adaptive grid provides significant savings in terms of the number of degrees of freedom.

\begin{figure}[!ht]
    \centering
    \includegraphics[width=0.6\textwidth]{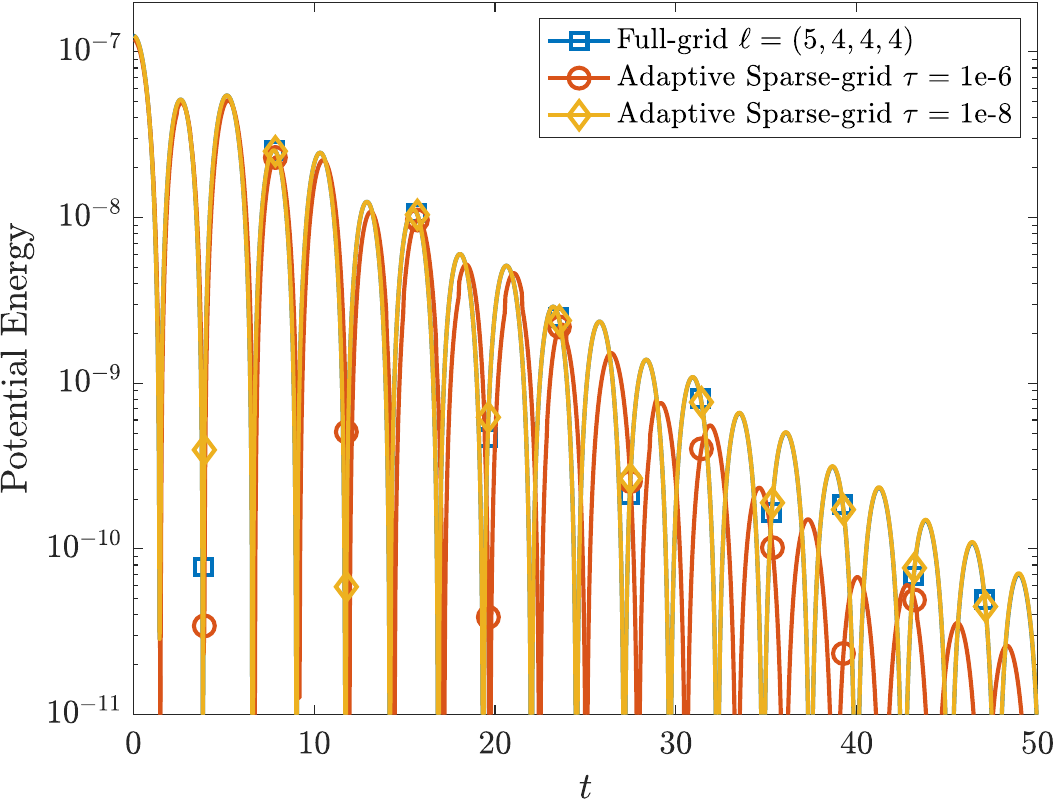}
    \caption{Collisional Landau  Problem -- \Cref{subsec:landau}: Plot of the potential energy with $\nu=1$.    
    The adaptive sparse grid is not allowed to refine past level $\bm{\ell}=(5,4,4,4)$, and the GMRES tolerance is set to $10^{-11}$.  A tolerance of $\tau=10^{-6}$ is not sufficient to capture the proper decay. The tolerance of $\tau=10^{-8}$ agrees quite well with the full-grid solution at all times plotted.}
    \label{fig:landau:nu1e0:asg_vs_fg}
\end{figure}

\begin{figure}[!ht]
    \centering
    \begin{subfigure}[b]{0.49\textwidth}
        \centering
        \includegraphics[width=\textwidth]{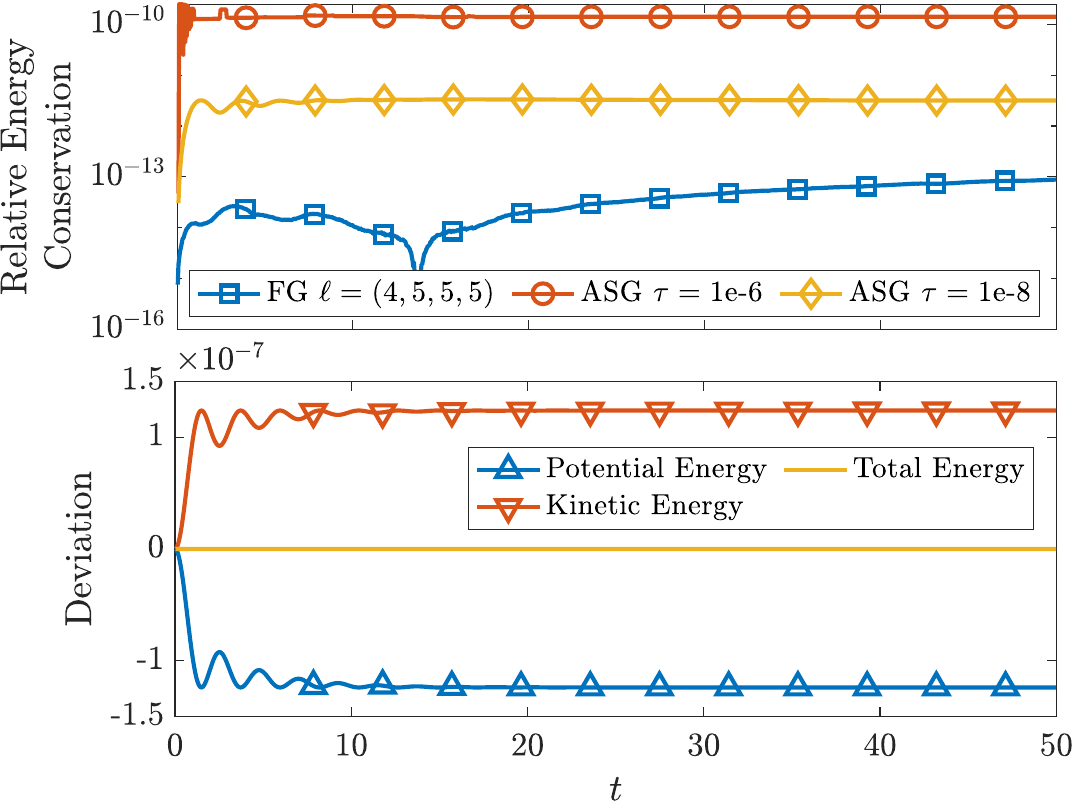}
        \caption{$\nu=10^{-2}$.  The GMRES tolerance is set to $10^{-14}$}
        \label{fig:landau:conservation:a}
    \end{subfigure}
    \begin{subfigure}[b]{0.49\textwidth}
        \centering
        \includegraphics[width=\textwidth]{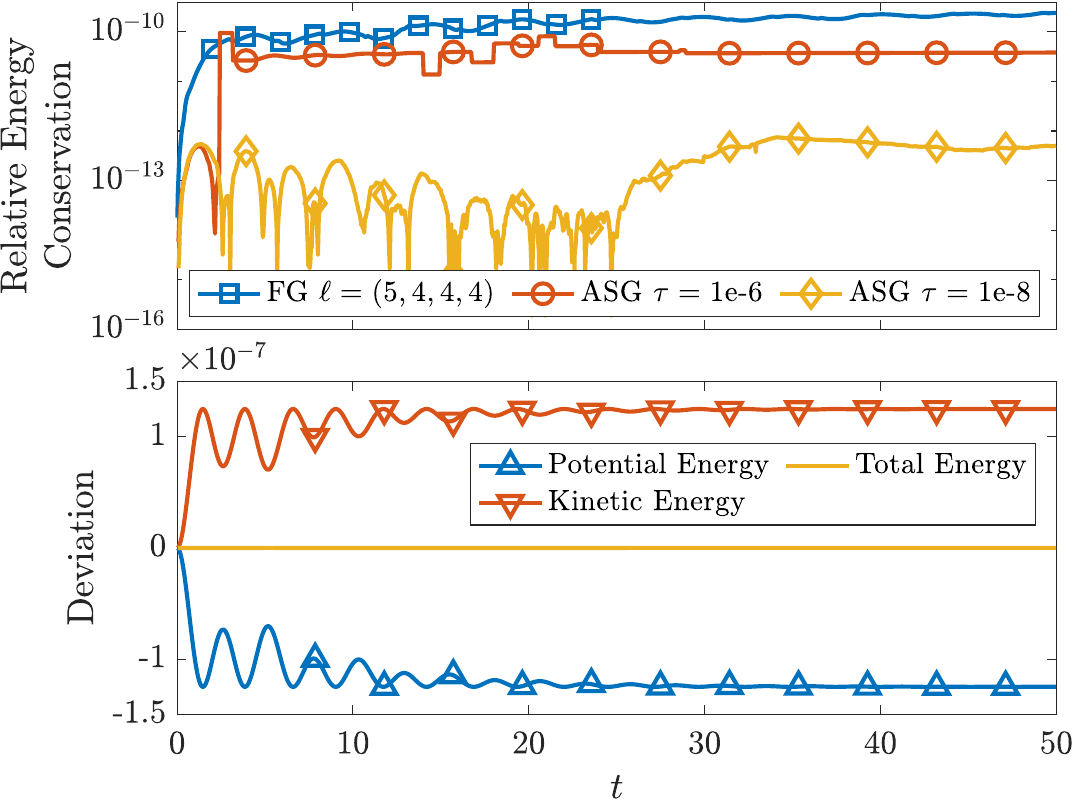}
        \caption{$\nu=1$.  The GMRES tolerance is set to $10^{-11}$}
        \label{fig:landau:conservation:b}
    \end{subfigure}
    \caption{Collisional Landau Problem -- \Cref{subsec:landau}: Top: Relative change in total energy versus time. Bottom:  Deviation of the potential, kinetic, and total energy from the initial condition for the adaptive sparse-grid method with $\tau=10^{-8}$.}
    \label{fig:landau:conservation}
\end{figure}

\Cref{fig:landau:nu1e0:asg_vs_fg} compares adaptive sparse-grid against full-grid for the moderate collisionality case, where we plot the same quantities as in \Cref{fig:landau:nu1e-2:asg_vs_fg}.  
For this collisionality, we have determined that a full-grid resolution of $\bm{\ell}=(5,4,4,4)$ is sufficient to accurately capture the evolution of the potential energy.
Similar to the low collisionality case, the potential energy evolution obtained with adaptivity threshold $\tau=10^{-6}$ is not in satisfactory agreement with the full-grid and analytic results.  
However, we find that the adaptive spare-grid and full-grid results are indistinguishable when the adaptivity threshold is reduced to $\tau=10^{-8}$.  
For $\tau=10^{-6}$, the number of active elements stays roughly constant at about $4.8\times10^{3}$, while for the case with $\tau=10^{-8}$, the number of active elements starts out around $1.1\times10^{4}$, and decreases to about $7\times10^{3}$ at the end of the simulation.  
For comparison, the full-grid with $\bm{\ell}=(5,4,4,4)$ consists of about $1.3\times10^{5}$ elements.  
Thus, the adaptive sparse-grid with $\tau=10^{-8}$ is as accurate as the full-grid solution, but with substantially fewer degrees of freedom.

In \Cref{fig:landau:conservation} we plot the the relative change in total energy for both collisionalities discussed above.  The relative change in the total energy is at the level of GMRES tolerance for the full-grid simulation.  For the adaptive sparse-grid methods, the relative change in the total energy decreases with the size of the threshold $\tau$ used; we expect this trend to continue until the GMRES tolerance pollutes the energy conservation.  We hypothesize the improvement in the relative energy conservation of the adaptive sparse-grid with $\tau=10^{-8}$ when compared with the full-grid (as seen in \Cref{fig:landau:conservation:b}) is due to the multiwavelets not being used in the Chu reduction discretization. 

The number of GMRES iterations varies between three and five for the sparse-grid runs. 
\section{Summary and Outlook}
\label{sec:conclusions}

In this work, we presented an adaptive sparse-grid DG method for the the VPLB model on a slab geometry.  The results of this project utilized the Adaptive Sparse-Grid Discretization (\asgard) codebase.
As demonstrated in \Cref{sec:numerical}, the adaptive sparse-grid method significantly decreases the storage cost of DG numerical approximations without compromising accuracy. 
Moreover, the adaptive sparse-grid method was able to capture physically relatively features of the distribution without the use of model specific error indicators.
The results also indicate that standard sparse-grids in velocity space, i.e.~the mixed-grid formulation, accurately captures low-order moments of the distribution, but are only slightly better when compared asymptotically against the full-grid for higher-order moments.  This necessitates further research into using in a coordinate system that more beneficially captures the radial behavior of the Maxwellian, e.g.~spherical-polar coordinates, or allowing some form of adaptivity in the mixed-grid.   
Other future plans include the expansion of the adaptive sparse-grid tests to full $3x3v$ phase-space simulations, efficient implementations of PDE operators on a sparse-grid basis, and the preservation of key quantities such as positivity of the discrete distribution in the multiwavelet basis.

\bibliography{references}

\end{document}